\documentclass[11pt]{amsart}
\usepackage[a4paper,margin=30mm]{geometry}
\usepackage[T1]{fontenc}
\usepackage{amsmath,amssymb,amsthm,mathtools,graphicx,tikz,booktabs,longtable,array}
\usepackage[hidelinks]{hyperref}
\IfFileExists{microtype.sty}{\usepackage{microtype}}{}
\newtheorem{theorem}{Theorem}[section]
\newtheorem{proposition}[theorem]{Proposition}
\newtheorem{lemma}[theorem]{Lemma}

\theoremstyle{definition}

\theoremstyle{remark}

\numberwithin{equation}{section}

\DeclareMathOperator{\per}{per}
\DeclareMathOperator{\Hypergeom}{Hypergeom}
\newcommand{\E}{\mathbb E}
\newcommand{\PP}{\mathbb P}
\newcommand{\ind}{\mathbf1}

\providecommand{\tightlist}{\setlength{\itemsep}{0pt}\setlength{\parskip}{0pt}}
\newcommand{\bigdotcup}{\mathop{\dot\bigcup}}

\hypersetup{pdftitle={The Borodin--Kostochka inequality from degree 750},pdfauthor={Lizhong Chen}}
\title{The Borodin--Kostochka inequality from degree 750}
\author{Lizhong Chen}
\date{}
\makeatletter
\let\bk@setauthors\@setauthors
\renewcommand{\@setauthors}{%
  \bk@setauthors
  \begingroup
  \centering\footnotesize
  Department of Mathematics\\
  Hong Kong University of Science and Technology\\
  Clear Water Bay, Hong Kong\\
  \href{mailto:lchendh@connect.ust.hk}{\texttt{lchendh@connect.ust.hk}}\par
  \endgroup
}
\makeatother

\begin{document}
\begin{abstract}
We prove that every finite simple graph \(G\) with maximum degree
\(\Delta\ge750\) satisfies the Borodin--Kostochka inequality
\(\chi(G)\le\max\{\omega(G),\Delta-1\}\), conjectured for all
\(\Delta\ge9\). Reed (1999) established the inequality for
\(\Delta\ge10^{14}\), and Liu, Sun, Wang and Zeng (2026) lowered the
threshold to \(5.3\times10^6\). To extend a random partial
colouring over large cliques, we find repeated colours using adaptively
chosen vertex-disjoint stars and apply Hall's theorem to the available
colour lists. We verify the required probability bounds by rational
and interval arithmetic with directed rounding.
\end{abstract}
\maketitle
\section{Introduction}

For a graph \(G\), let \(\chi(G)\), \(\omega(G)\) and \(\Delta(G)\) denote
its chromatic number, clique number and maximum degree, respectively.
A proper colouring assigns different colours to adjacent vertices.

By Brooks' theorem, every connected graph \(G\) other than a complete
graph or an odd cycle admits a proper colouring with at most \(\Delta(G)\)
colours \cite{Brooks}. Borodin and Kostochka conjectured that when
\(\Delta(G)\ge9\), one fewer colour suffices unless \(G\) contains a clique
that requires more \cite{BK}. In other words,
\[
 \chi(G)\le\max\{\omega(G),\Delta(G)-1\}
 \qquad\text{whenever }\Delta(G)\ge9.
\]
The degree restriction is necessary. Replacing each vertex of a five-cycle
by a triangle and joining consecutive triangles completely gives a graph
of maximum degree eight, clique number six and chromatic number eight.

Reed proved the conjecture for sufficiently large maximum degree
\cite{Reed}. Liu, Sun, Wang and Zeng gave the explicit bound
\(5.3\times10^6\) \cite{LSWZ}, and an earlier preprint of the author
lowered it to \(1000\) \cite{Chen1000}. We lower the threshold further.

\begin{theorem}\label{thm:main}
For every finite simple graph \(G\),
\[
 \Delta(G)\ge750
 \quad\Longrightarrow\quad
 \chi(G)\le\max\{\omega(G),\Delta(G)-1\}.
\]
\end{theorem}

The theorem applies to disconnected graphs as well as connected ones, with
no restriction on their order. The finite computation in the proof concerns
local probabilities. We reduce these estimates to inequalities in bounded
integer parameters whose ranges are independent of the order of \(G\).
The computational appendix gives the ranges and the evaluations.

We follow the approach of Reed and of Liu, Sun, Wang and Zeng, combining
a structural decomposition with a random partial colouring \cite{Reed,LSWZ}.
Rabern's decomposition \cite{RabernDense} provides the clique parts.
We use list-colouring results of Cranston and Rabern \cite{CR} to
distinguish the large clique parts that we colour as a whole from the
vertices that we can colour greedily. Independently at each vertex, we
assign a random colour and a random priority. A vertex keeps its colour
if it precedes every neighbour assigned the same colour. Pemmaraju and
Srinivasan studied this priority rule \cite{PS}, and Hurley, de Joannis
de Verclos and Kang used it in their colouring argument \cite{HDK}.
When several neighbours keep the same colour, they leave more colours
available at the vertex.

The degree-\(1000\) proof \cite{Chen1000} uses the same priority rule.
Its Hall condition uses both colours present in the initial clique
palette (represented colours) and colours absent from it (fresh colours).
It also uses matching arguments, star families chosen from fixed
endpoint sets, sampling without replacement and rank-sensitive colour
estimates. To reach
degree \(750\), we strengthen the lower bound on the initial number of
stars and estimate how many stars can be recovered in a larger endpoint
set once the available centres are known.

To see how a repeated colour helps on a clique \(C\), consider an
uncoloured vertex of \(C\) with a coloured neighbour outside it. If
the neighbour's colour is already used on \(C\), it removes no further
colour from the vertex's list. Two distinct repetitions of this kind at one vertex
and one at another suffice to colour the rest of \(C\) by Hall's theorem.
A second sufficient condition combines repetitions at different clique
vertices with fresh colours on their outside neighbours. Both conditions
persist as vertices outside \(C\) are coloured, so the clique parts can
be coloured one at a time.

We find these colours on vertex-disjoint stars. A \emph{\(d\)-star}
centred at \(c\in C\) consists of the edges from \(c\) to \(d\) distinct
neighbours outside \(C\), called its endpoints. The endpoints may be
adjacent in \(G\). We fix the set of possible endpoints in advance. The
initial colours on \(C\) and the relative priority order within each
colour class then tell us which vertices have already lost their colours
to clique neighbours. Choosing the stars at this stage lets us use
centres missed by the original family. The fixed endpoint set determines
which other events the clique event can depend on. The choice of stars
within that set improves our bound on its probability.

Combining a paired form of the Tutte--Berge formula with the bound
on common neighbours of endpoints gives a lower bound on the initial
number of disjoint two-stars. We can then enlarge the endpoint set to include all
eligible neighbours of unused centres. For a given number of available
centres, we obtain an affine lower bound on the expected maximum packing
size as a function of the original packing size. A bounded-difference
estimate controls the lower tail. Together with an exponential bound on
the conditional failure probability, this gives a log-convex upper bound
on the ratio of the failure probability to the Local Lemma allowance.
The allowance accounts for the enlarged endpoint set, and checking the
two ends of an interval suffices throughout that interval. In estimating
colour probabilities, we keep track of how many represented colours occur
on the endpoints, the order in which their labels are assigned and how
many neighbours the endpoints have in the clique.

For a vertex to be coloured greedily, the size of a largest clique
containing it controls how many edges are missing in its neighbourhood.
A weighted matching estimate then bounds the probability that too few
pairs of neighbours retain a common colour. We use the same random
partial colouring for these vertices and the clique parts, so a single
application of the asymmetric Local Lemma gives the required outcome.

Under the stronger clique restriction \(\omega(G)\le\Delta(G)-4\),
Cranston and Rabern proved \(\chi(G)\le\Delta(G)-1\) for
\(\Delta(G)\ge13\) \cite{CRBig}.
For the list-colouring version, Choi, Kierstead and Rabern
proved the conjectured bound at maximum degree at least \(10^{20}\)
\cite{CKR}. Dvo\v{r}\'ak, Kang and Mik\v{s}an\'ik proved the
correspondence-colouring version for maximum degree at least
\(3\cdot10^9\) \cite{DKM}. Theorem~\ref{thm:main} concerns the ordinary
chromatic number.

We first give sufficient conditions for completing a partial colouring.
They determine the bad events that the structural and probabilistic
arguments must address. After estimating the probabilities and the
available star packings, we combine the estimates in seven ranges
according to how far the clique order falls below the maximum degree. We include the proofs needed from the degree-\(1000\)
argument. Appendices~\ref{app:structure}--\ref{app:verification} give
the structural and small-deficit arguments, the probability recurrences
and the details of the finite computation.

\section{From a partial colouring to a colouring}
\label{sec:extension}

We complete the colouring in three steps. First we colour the exceptional
vertices attached to clique parts. We then colour the designated clique
parts one at a time and finish by colouring the remaining vertices
greedily. Section~\ref{sec:structure} gives the partition. Here we
determine what the initial partial colouring must provide at each step.

Put \(D=750\) and \(q=D-1\), and let \(G\)
have maximum degree at most \(D\). A \emph{partial \(q\)-colouring}
\(\varphi\) assigns colours from \([q]=\{1,\ldots,q\}\) to some vertices
and is proper on its coloured vertices.
Every extension considered below preserves colours already assigned.

For a vertex \(v\), let \(c_\varphi(v)\) be its number of coloured neighbours
and let \(a_\varphi(v)\) be the number of distinct colours on those neighbours.
The quantity
\[
 s_\varphi(v)=c_\varphi(v)-a_\varphi(v)
\]
counts the colours saved by repetitions. When an additional neighbour is
coloured, this quantity either stays unchanged or increases by one.
Consequently, if an uncoloured vertex already satisfies
\begin{equation}
s_\varphi(v)\ge d(v)-q+1,
 \label{eq:pilot-saving}
\end{equation}
then it can be coloured at any later stage, since its neighbours use at
most \(d(v)-s_\varphi(v)\le q-1\) distinct colours. Thus a vertex of
degree \(D\) needs two savings, and a vertex of degree \(D-1\) needs one.

On a large clique, savings at just a few vertices can suffice. Hall's
theorem lets us use their larger lists to colour the whole clique.

\subsection{The list-size deficit}
\begin{figure}[htbp]
\centering
\begin{tikzpicture}[x=1cm,y=1cm,
  vertex/.style={circle,draw,minimum size=6mm,inner sep=0pt},
  coloura/.style={vertex,fill=blue!18},
  colourb/.style={vertex,fill=orange!30}]
 \draw[rounded corners,draw=black!55,fill=black!2]
   (0,-1.55) rectangle (4,2.8);
 \node at (2,2.45) {clique \(C\)};
 \draw[rounded corners,dashed,draw=black!40]
   (0.2,-1.3) rectangle (3.8,0.76);
 \node[font=\small] at (2.45,-0.92) {uncoloured \(U\)};
 \node[coloura,label=above:\(a\)] (a) at (1.3,1.35) {\(\alpha\)};
 \node[colourb,label=above:\(b\)] (b) at (2.7,1.35) {\(\beta\)};
 \node[vertex,fill=white] (u) at (1,0.35) {\(u\)};
 \node[vertex,fill=white] (w) at (1,-0.45) {\(w\)};
 \node[coloura,label=below:\(x\)] (x) at (6,-0.45) {\(\alpha\)};
 \node[colourb,label=above:\(y\)] (y) at (6,0.65) {\(\beta\)};
 \draw (u)--(x);
 \draw (u)--(y);
 \draw (w)--(x);
 \node[align=left,anchor=west,font=\small] at (6.65,-0.45)
   {already excluded\\from every list on \(U\)};
 \node[align=left,anchor=west,font=\small] at (6.65,0.65)
   {already excluded\\from every list on \(U\)};
\end{tikzpicture}
\caption{The colours \(\alpha\ne\beta\) retained by \(a,b\in C\) have
already been removed from every list on \(U\). Thus \(x,y\) forbid no
further colours at \(u\), and \(x\) forbids no further colour at \(w\).
The list at \(u\) has size at least \(k\), and the list at \(w\) has
size at least \(k-1\). Only the outside edges used in this argument are
shown.}
\label{fig:repeated-colours}
\end{figure}
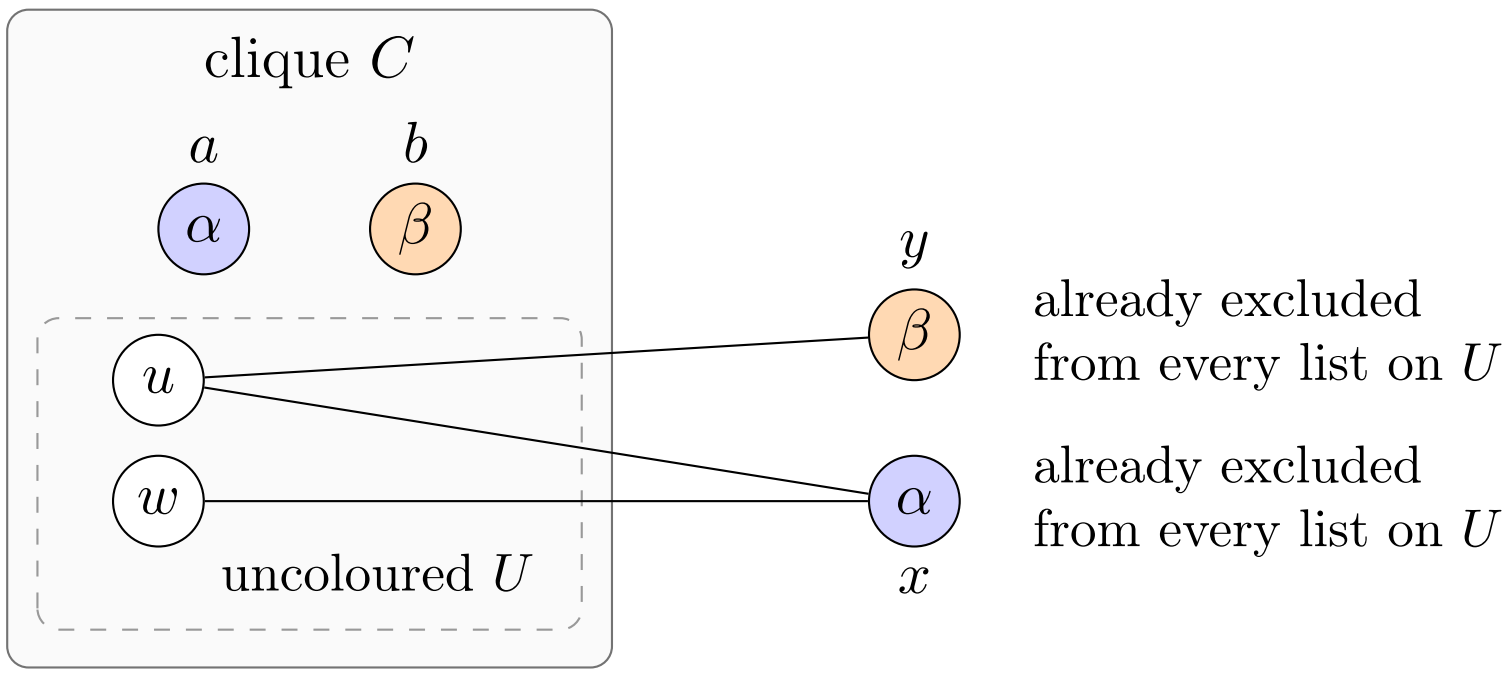

For a clique \(C\), call \(p=D-|C|\) its \emph{clique deficit}.
Suppose that \(k\) vertices of \(C\) remain uncoloured.
The coloured vertices of \(C\) have distinct colours.
Deleting those colours from \([q]\) leaves a palette of size \(k+p-1\).
Every vertex of \(C\) has at most \(p+1\) neighbours outside \(C\), so the list
\(L(v)\) of colours still available at an uncoloured vertex satisfies
\[
 |L(v)|\ge k-2.
\]
To colour the \(k\) vertices of the clique, we need \(k\) distinct colours
between their lists. Repetitions can help close the gap. If a coloured
neighbour of \(v\) outside \(C\) has a colour already used on \(C\), that
neighbour removes no further colour from \(L(v)\). Two distinct such
colours at one vertex \(u\) give \(|L(u)|\ge k\), and one at a different
vertex \(w\) gives \(|L(w)|\ge k-1\). Thus
\[
 |L(u)|\ge k,\qquad |L(w)|\ge k-1,\qquad
 |L(v)|\ge k-2\quad(v\in U),
\]
where \(U\) is the set of uncoloured clique vertices. The lists need not
be nested. Figure~\ref{fig:repeated-colours} illustrates
how the repeated colours give these larger lists.

These list sizes imply Hall's condition. For any nonempty set of at most
\(k-2\) vertices, even one of its lists has at least that many colours.
A set of \(k-1\) vertices contains at least one of \(u,w\), and the
whole set \(U\) contains \(u\). In each case the union of the lists has
at least as many colours as there are vertices to colour.

The probability argument begins with an initial colour assignment
\(\sigma:V(G)\to[q]\) and distinct priorities at all vertices. We retain
\(\sigma(v)\) precisely when no neighbour with that initial colour has
smaller priority. Let \(\varphi\) be the resulting proper partial colouring.
Let \(A_C\) consist of the vertices of \(C\) preceded by another vertex
of the same initial colour in \(C\). Every vertex of \(A_C\) is
uncoloured, even if the earlier vertex also loses its colour.

The set \(\sigma(C)\) includes every colour initially assigned to \(C\),
even if no vertex retains it. Colours absent from \(\sigma(C)\) give a
second way to satisfy Hall's condition. We state it in a form that lets
us colour other parts of the graph before returning to \(C\).

\subsection{A condition that survives outside extensions}
Fix an initial assignment and priorities as above, and a clique \(C\) of
order \(D-p\), with \(1\le p\le115\). Choose distinct
\(v_1,\ldots,v_t\in A_C\). Record the following colours in \(\varphi\).
For each \(i\), let \(r_i\) count distinct colours that occur both on a
coloured neighbour of \(v_i\) outside \(C\) and on a coloured vertex of \(C\).
Let \(T_i\) be a set of colours that occur on coloured outside neighbours of
\(v_i\) and are absent from \(\sigma(C)\). Put \(f_i=|T_i|\).
Assume that the sets \(T_i\) are pairwise disjoint.
Their disjointness will help rule out a common set of \(p\) forbidden
colours, even when they do not enlarge the individual lists.
We write \(x_+=\max\{x,0\}\).

\begin{lemma}\label{lem:clique-extension}
Suppose that \(|\{i:r_i\ge1\}|\ge2\), and that at least one of the following
two inequalities holds:
\begin{equation}
\max_i r_i\ge2
 \qquad\hbox{or}\qquad
 \sum_{i=1}^t(r_i+f_i-1)_+>p.
 \label{eq:pilot-certificate}
\end{equation}
Every proper extension of \(\varphi\) obtained by colouring vertices outside
\(C\), while leaving the vertices of \(C\) unchanged, extends to colour \(C\).
\end{lemma}
\begin{proof}
After such an outside extension, let \(U\) be the set of uncoloured
vertices of \(C\), put \(k=|U|\), and let \(L(v)\) be the list at
\(v\in U\). At least two selected vertices have \(r_i\ge1\), and both
remain in \(U\), so \(k\ge2\).
Delete the colours used in \(C\) from \([q]\), leaving a palette \(R_C\) of
size \(k+p-1\). The available lists satisfy
\[
 |L(v)|\ge k-2\quad(v\in U),\qquad
 |L(v_i)|\ge k-2+r_i.
\]
Every \((k-1)\)-set contains one of the two selected vertices with
\(r_i\ge1\). Hence every set of at most \(k-1\) vertices satisfies Hall's condition.
If some \(r_i\ge2\), the full set does too.

It remains to consider the case in which every \(r_i\le1\). If the union
of the lists had size at most \(k-1\), at least \(p\) colours of \(R_C\)
would be forbidden at every vertex of \(U\). Choose a set \(Q_0\) of
\(p\) such colours. At \(v_i\) at most \(p+1-r_i\) colours of \(R_C\) are
forbidden, and these include both \(Q_0\) and \(T_i\). Therefore
\[
 |T_i\cap Q_0|\ge(r_i+f_i-1)_+.
\]
The sets \(T_i\cap Q_0\) are disjoint subsets of the \(p\)-element set
\(Q_0\). Their sizes cannot sum to more than \(p\), contrary to
\eqref{eq:pilot-certificate}. Hall's theorem now gives a colouring of \(U\).
\end{proof}

Let \(E\subseteq N(C)\setminus C\). A \emph{Hall certificate in \(E\)}
for \(C\) consists of distinct centres \(v_i\in A_C\) and selected colours
on retained neighbours of each \(v_i\) in \(E\). At each centre, \(r_i\)
counts the distinct selected colours also retained on \(C\), and \(T_i\)
is the set of selected colours absent from \(\sigma(C)\). The \(T_i\) must
be pairwise disjoint, at least two \(r_i\) must be positive, and
\eqref{eq:pilot-certificate} must hold with \(f_i=|T_i|\). These selected
colours suffice for Lemma~\ref{lem:clique-extension}, since using all
retained neighbours outside \(C\) can only increase the \(r_i\). The selected repetitions
persist under outside extensions, and each colour in \(T_i\) remains absent
from \(C\).

\subsection{Completing the colouring}
In the partition of a minimal counterexample, each clique part has at
most one exceptional vertex attached to it. We designate certain clique
parts, all of the size required by Lemma~\ref{lem:clique-extension}, to
be coloured as a whole. Once the exceptional vertices and the designated
parts have been coloured, we colour every remaining vertex greedily.
These include the remaining vertices of the non-designated clique parts.

For this partition, it suffices to find a proper partial \(q\)-colouring
with the following properties:
\begin{enumerate}
 \item every exceptional vertex has two uncoloured neighbours in its own
 clique part;
 \item every designated clique part \(C\) has a Hall certificate in
 \(N(C)\setminus C\);
 \item every uncoloured vertex in the final greedy phase satisfies
 \eqref{eq:pilot-saving}.
\end{enumerate}
First colour the uncoloured exceptional vertices. The two reserved
neighbours of each such vertex lie in its own clique part and are still
uncoloured, so at most \(D-2=q-1\) neighbours already have colours.
Next colour the designated clique parts one at a time. The Hall
certificate for each part persists while the preceding parts are
coloured, by Lemma~\ref{lem:clique-extension}. The savings at the
remaining vertices also persist, so we may colour them in any order.

To obtain such a partial colouring, choose a uniform colour in \([q]\)
and a uniform priority in \((0,1)\) independently at each vertex, with all
choices independent. Ties have probability zero. The retention rule uses
the initial assignments, so even a vertex that loses its own colour can
cause a neighbour to lose that colour.

We define a bad event for each vertex to be coloured greedily, each
designated clique part and each exceptional vertex. The asymmetric Local
Lemma will give an outcome in which all three requirements hold.

\section{The structural reduction and the bad events}
\label{sec:structure}

We obtain the partition used in Section~\ref{sec:extension} by first
reducing to a critical graph of degree \(750\).
A graph is \emph{vertex-critical} if deleting any vertex lowers its
chromatic number. By Brooks' theorem, a counterexample at maximum degree \(k\ge750\)
contains a vertex-critical induced subgraph with \(\chi=\Delta=k\) and
\(\omega<k\). The Kostochka--Catlin reduction
\cite[Lemma~1.9]{CR} gives a graph with these properties at degree \(k-1\)
whenever one exists at degree \(k>750\). It suffices to exclude a
counterexample at \(D=750\).

Among all counterexamples of maximum degree \(D\), choose \(G\) of
minimum order and, subject to this, with the most edges. Then \(G\) is
connected and vertex-critical, and
\begin{equation}
\chi(G)=D,\qquad \omega(G)\le q,\qquad q\le d(v)\le D.
 \label{eq:critical-range}
\end{equation}
In the terminology of Cranston and Rabern, a \emph{child} of \(G\) is a
graph \(J\) nonisomorphic to \(G\) such that some induced subgraph of \(G\)
admits a vertex-surjective homomorphism to \(J\). A \emph{\(D\)-mule}
is a vertex-critical graph with \(\chi=\Delta=D>\omega\) having no
child in that same class. A child of the same order must have more
edges, since the homomorphism is then a bijection. Our extremal choice
therefore makes \(G\) a \(D\)-mule \cite[Lemma~2.1]{CR}, and
Corollary~2.11 of the same paper gives
\begin{equation}\label{eq:maximum-clique-attachment}
 |N(x)\cap K|\le1
 \qquad\text{for every \((D-1)\)-clique \(K\) and every \(x\notin K\)}.
\end{equation}
Both exceptions to that corollary have maximum degree seven.

Minimal order also gives the following colouring property for smaller graphs:
\begin{equation}\label{eq:smaller}
 |V(J)|<|V(G)|,\quad \Delta(J)\le D,\quad\omega(J)\le q
 \quad\Longrightarrow\quad \chi(J)\le q.
\end{equation}
Minimality gives a \(q\)-colouring of each component of \(J\) with
maximum degree \(D\). For components of smaller maximum degree, Brooks'
theorem suffices, since \(J\) contains no complete graph on \(D\)
vertices and an odd cycle needs only three colours.

\subsection{Clique parts}\label{sec:clique-parts}
A graph \(H\) is \(d_1\)-choosable if every list assignment satisfying
\(|L(v)|\ge d_H(v)-1\) admits a proper colouring. Such a nonempty induced
subgraph cannot occur in \(G\). Indeed, \eqref{eq:smaller} gives a
colouring of \(G-V(H)\). Each vertex of \(H\) then has a list of size
at least \(q-d_G(v)+d_H(v)\ge d_H(v)-1\), so the colouring would extend to \(G\).

Write \(I_t\) for the edgeless graph on \(t\) vertices.
We use two results of Cranston and Rabern \cite[Lemmas~3.3 and~3.11]{CR}.
For \(t\ge4\), the join \(K_t\vee B\) is \(d_1\)-choosable unless \(B\)
is almost complete, meaning \(\omega(B)\ge|B|-1\), or one of the following
exceptions occurs: \(t=4\) and \(B\in\{I_3,K_{1,3}\}\), or \(t=5\) and
\(B=I_3\). Also, \(H\vee I_2\) is \(d_1\)-choosable whenever \(H\) is
connected, incomplete and has at least four vertices.

A consequence of the first result will be used several times.
Let \(C\) be a clique of order \(m\ge7\). If two distinct vertices
\(x,y\) outside \(C\) have four common neighbours in \(C\), let \(A\)
consist of four such neighbours. Then
\[
 G[C\cup\{x,y\}]=K_4\vee G[(C\setminus A)\cup\{x,y\}].
\]
The second factor has order \(m-2\ge5\), so none of the small exceptions
applies. It must be almost complete. Hence \(C\cup\{x,y\}\) contains a
clique of order at least \(m+1\) containing \(A\).

Put \(h_0=444\). Form the intersection graph of the maximal cliques of
order at least \(h_0+1\), joining two cliques when they meet. Its
components determine our partition.
\begin{proposition}\label{prop:decomposition}
There is a partition
\[
 V(G)=L\mathbin{\dot\cup}\mathop{\dot\bigcup}_i S_i,\qquad
 S_i=C_i\quad\hbox{or}\quad S_i=C_i\cup\{u_i\},
\]
where \(C_i\) is a maximum-order representative of one component of the
intersection graph of the maximal cliques of order at least \(445\).
Writing \(|C_i|=D-p_i\), we have \(1\le p_i\le305\), and:
\begin{enumerate}
\renewcommand{\labelenumi}{(\roman{enumi})}
 \item an exceptional vertex \(u_i\), when present, has at least \(444\)
 neighbours in \(C_i\);
 \item every vertex of \(C_i\) has at least \(p_i\) outside neighbours
 whose number of neighbours in \(C_i\) is at most \(p_i+3\);
 \item two distinct outside vertices have at most three common
 neighbours in \(C_i\);
 \item if edges of maximum added degree at most \(7\) are added to
 \(G-C_i\), and total maximum degree remains at most \(D\), the resulting
 graph is \(q\)-colourable.
\end{enumerate}
Every vertex of \(L\) belongs to no clique of order greater than \(444\).
\end{proposition}

The description of each part as a clique with at most one additional
vertex is a special case of Rabern's decomposition
\cite[Lemma~4.1]{RabernDense}. Appendix~\ref{app:structure} proves this
and the additional properties in the proposition. The main observation is
that the exclusion of \(d_1\)-choosable induced subgraphs restricts the
form of every induced subgraph on at most \(D+1\) vertices with minimum
degree at least \(444\). Each consists of a clique and at most one
additional vertex. Applying this to unions of intersecting large cliques
gives the disjoint sets \(S_i\).
The common-neighbour consequence above then proves (iii), since a
larger clique meeting \(C_i\) would contradict its choice as representative.
Finally, the same dense-subgraph conclusion prevents a new \(D\)-clique
when edges are added as in (iv), so \eqref{eq:smaller} applies.

The \emph{attachment} of an outside neighbour \(z\) to \(C_i\) is
\(|N_G(z)\cap C_i|\), measured in the original graph. Call \(z\)
\emph{good for \(C_i\)} if its attachment is at most \(p_i+3\).
The \emph{good-incidence graph} has parts \(C_i\) and its good outside
neighbours, with all edges of \(G\) between the parts. A good \(d\)-star
is a \(d\)-star in this graph with its centre in \(C_i\).
In the good-incidence graph the centre degrees lie between \(p_i\)
and \(p_i+1\), the endpoint degrees are at most \(p_i+3\), and (iii)
bounds the number of common centres of any two endpoints by three.

We designate the parts with \(p_i\le115\) for Hall extension. Define
\[
 \mathcal V_A=L\cup\bigcup_{p_i>115}C_i.
\]
These are the vertices assigned to the final greedy phase. The sets
\(\mathcal V_A\), the designated clique parts, and the exceptional vertices
are disjoint. For \(v\in\mathcal V_A\), choose a largest clique \(C(v)\)
containing \(v\), and put \(c(v)=|C(v)|\). Then
\begin{equation}\label{eq:root-clique}
 c(v)\le634.
\end{equation}
For \(v\in L\) this follows from the definition of \(L\). Otherwise every
large clique through \(v\) lies in the same intersection component as
its non-designated part, whose chosen representative has order at most
\(D-116=634\). The cliques \(C(v)\) need not be disjoint.

\subsection{One event for each required extension}
Use the random assignment and priorities of Section~\ref{sec:extension},
and let \(\varphi\) be the retained proper partial colouring.
For \(v\in\mathcal V_A\), define
\[
 A_v=\{v\text{ is uncoloured and }s_\varphi(v)\le d(v)-q\}.
\]
We call \(v\) the root of \(A_v\). If \(A_v\) is avoided, then \(v\)
either already has a colour or has the saving required for greedy colouring.

For each designated clique \(C\), fix a set
\(E(C)\subseteq N(C)\setminus C\) before assigning any random colours
or priorities. Section~\ref{sec:closure} specifies the choice from the
graph. Let \(F_C\) be the event that \(C\) has no Hall certificate in
\(E(C)\). After seeing the clique colours and their within-class priority
orders, we will choose stars in this set to bound the probability of \(F_C\).

For each exceptional vertex \(u_i\), fix a set
\(R_i\subseteq N(u_i)\cap C_i\) of \(h_0=444\) vertices and put
\[
 W_i=\{h_0-|\sigma(R_i)|\le1\}.
\]
Outside \(W_i\), at least two vertices of the clique \(R_i\) are not first
in their initial colour class within \(R_i\). They are therefore
uncoloured neighbours of \(u_i\). This remains true even if their
earlier same-coloured neighbours also discard their colours.

By Section~\ref{sec:extension}, any outcome avoiding all the events
\(A_v,F_C,W_i\) gives a \(q\)-colouring of \(G\). It remains to show that
the probability of such an outcome is positive.

\subsection{Dependencies and the probability targets}
Write \(B_r(S)\) for the closed radius-\(r\) ball about \(S\). For each
bad event, the random marks on the following set determine whether it occurs.
\[
 X(A_v)=B_2(v),\qquad
 X(F_C)=B_1(C\cup E(C))\subseteq B_2(C),\qquad
 X(W_i)=R_i.
\]
Both colour and priority marks are included except for \(W_i\), which
uses only colours. Join two events when their supports intersect.
An event is independent of the sigma-field of all its non-neighbours,
because their combined mark coordinates are disjoint from its own.

Assign activity \(x_0\) to each event \(A_v\) or \(F_C\), and activity
\(x_W\) to each \(W_i\), where
\[
 x_0=D^{-4},\qquad x_W=2\cdot10^{-37},\qquad
 \eta_0=\frac{x_0}{1-x_0},\qquad \eta_W=\frac{x_W}{1-x_W},
 \qquad S_3=D^2-D+1.
\]
If \(|C|=m=D-p\) and \(|E(C)|\le e\), define
\begin{align}
 J_1(p,e)&=m(p+2)+e(D-1),\label{eq:J1}\\
 J_3(p,e)&=m[1+(p+1)S_3]+e(D-1)^3.\label{eq:J3}
\end{align}
These quantities bound \(|X(F_C)|\) and \(|B_2(X(F_C))|\), respectively.
The clique has at most \(m(p+1)\) outgoing edges. Since each endpoint lies
in \(B_1(C)\), a vertex within distance \(r\) of that endpoint but outside
\(B_r(C)\) must be at distance exactly \(r\) from the endpoint. A shortest
path to that vertex cannot start with the chosen edge back to \(C\).
There are at most \((D-1)^r\) such paths per endpoint. Summing over the
fixed endpoint set for \(r=1,3\) gives \eqref{eq:J1}--\eqref{eq:J3}.

For an \(A\)-root whose chosen clique has order \(c\), write \(b=D-c\) and
define
\[
 Q_4=1+D^4-2D^3+2D^2,\qquad
 V(c)=\min\{Q_4,\ D+1+(b+1)(2D-b-2)S_3\}.
\]
This bounds \(|B_4(v)|\). The first bound counts nonbacktracking paths.
For the second, \(N[v]\) is the union of the clique \(C(v)\) and at most
\(b+1\) further vertices, each adjacent to \(v\); its edge boundary is at
most \((b+1)(2D-b-2)\). There are at most \(S_3\) further vertices along
paths of length at most two beyond each boundary edge. For \(c<250\)
we will simply use \(Q_4\).

To count neighbouring events of types \(A\) and \(F\), assign each one
a vertex in \(B_2(X)\), where \(X\) is the support under consideration.
If \(X(A_v)\) meets \(X\), use its root \(v\). If \(X(F_C)\) meets
\(X\), choose any vertex of \(C\) in \(B_2(X)\). The chosen vertices
are distinct because the designated clique parts are disjoint and
contain no \(A\)-roots. There are therefore at most \(|B_2(X)|\)
neighbouring events of these two types combined. Since the \(R_i\) are
disjoint, at most \(|X|\) events \(W_i\) meet \(X\), and no two distinct
\(W\)-events are neighbours.

These counts reduce the asymmetric Local Lemma inequalities to
\begin{align}
 \PP(A_v)&\le x_0\exp[-\eta_0 V(c(v))-\eta_W(1126125/2)],
       &&c(v)\ge250,\label{eq:A-target}\\
 \PP(A_v)&\le x_0\exp[-\eta_0 Q_4-\eta_W(1126125/2)],
       &&c(v)<250,\label{eq:A-mass-target}\\
 \PP(F_C)&\le x_0\exp[-\eta_0J_3(p,e)-\eta_WJ_1(p,e)],\label{eq:F-target}\\
 \PP(W_i)&\le x_W\exp[-\eta_0 V_W],
       &&V_W=444[1+(D-443)D].\label{eq:W-target}
\end{align}
Here \(1126125/2>D^2+1\ge|B_2(v)|\). For \(V_W\), each of the 444
vertices of the clique \(R_i\) has at most \(D-443\) outgoing edges,
each reaching at most \(D\) vertices within distance two along that edge.
For each activity \(x\), the inequality \(-\log(1-x)\le x/(1-x)\)
gives a lower bound on the product of the factors \(1-x\) contributed
by neighbouring events.

In this notation, the finite asymmetric Local Lemma states that if
\begin{equation}\label{eq:LLL}
 \PP(B_i)\le x_i\prod_{j\sim i}(1-x_j),\qquad 0\le x_i<1,
\end{equation}
and each event is independent of the sigma-field of its non-neighbours,
then all the events can be avoided with positive probability. Indeed,
conditional on avoiding any set of other events, the probability of
\(B_i\) is at most \(x_i\). This follows by induction on the size of the
conditioning set, using independence for non-neighbours and the inductive
bound for neighbours. Multiplying the conditional probabilities of
avoidance proves the claim \cite{EL}.

We can verify \eqref{eq:W-target} directly. The probability of at most
one repeated assignment among the \(444\) colours is at most
\[
 \left(1+\frac{444\cdot443}{2q}\right)
 \exp\left(-\frac{443\cdot442}{2q}\right)
 <133e^{-130}<10^{-37}.
\]
The series for \(e\) gives \(e>8/3\) and hence \(e^5>100\), which verifies
the last comparison. On the other hand,
\(\eta_0 V_W<1/2\) and \(e^{1/2}<2\), so the right side of
\eqref{eq:W-target} exceeds \(x_W/2=10^{-37}\).

In \eqref{eq:F-target}, \(J_1\) and \(J_3\) depend on the entire endpoint
set fixed before the experiment. Choosing a subfamily of stars later
can improve the bound on \(\PP(F_C)\), while the Local Lemma allowance
continues to use that entire set.

\section{Useful colours and complete probability bounds}
\label{sec:probability}

We estimate a clique event in two stages. First expose enough information
inside the clique to identify a useful star family. Then estimate the chance
that this family supplies the colours required by
Lemma~\ref{lem:clique-extension}. The estimate will depend on the size of the
available family, but not on the labels of its vertices. This permits us
to use the largest family available after the first exposure.

\subsection{Choosing stars after a partial exposure}
Fix a designated clique \(C\) of order \(m=D-p\) and its endpoint set
\(E=E(C)\). Expose the initial colours on \(C\) and the relative priority
order within each colour class. Write \(\mathcal C\) for this
information. The numerical priorities and the order of the minima of
different classes remain unexposed. The earliest vertex of each class
in \(C\) is its \emph{clique representative}.

Write \(O=|\sigma(C)|\) for the number of initial clique colours and put
\(h=m-O\). The exposure determines \(A_C\), the set of the \(h\) vertices
preceded by a vertex of their own colour in \(C\). Every vertex of
\(A_C\) is uncoloured. A representative can lose its colour to an outside
neighbour, but its colour still contributes to \(O\). We evaluate the
colour laws below for \(1\le O<q\) and use the upper bound one on
the failure probability at any other value.

Choose, as a deterministic function of \(\mathcal C\) and the graph, an
ordered family of \(K\) vertex-disjoint \(d\)-stars with centres in \(A_C\)
and endpoints in \(E\), where \(2\le d<q\). A maximum family, with ties
broken by a fixed graph order, is an allowed choice. So are its prefixes.
The choice uses neither endpoint marks nor numerical clique minima.
Conditional on \(\mathcal C\), the selected vertices are fixed, and all
marks outside \(C\) retain their independent original laws.
Vertices outside the selected family, including unused vertices of \(E\),
retain their original adjacencies and marks.

For a centre set \(A\subseteq C\), the \emph{available rank} \(R_E(A)\)
is the maximum number of vertex-disjoint \(d\)-stars with centres in \(A\)
and endpoints in \(E\).
We seek a nonincreasing function \(\psi_{p,O}\) for which
\(\PP(F_C\mid O)\le\E[\psi_{p,O}(R_E(A_C))\mid O]\).
We first bound the failure probability for a fixed exposure
\(\mathcal C\), then average over exposures with the same \(O\).
This section constructs \(\psi_{p,O}\). Section~\ref{sec:rank}
estimates its argument, the largest packing available in the fixed set.

Each colour occurring on the \(n=dK\) selected endpoints has a unique
earliest endpoint, its \emph{winner}. That colour is a \emph{target colour} of the
star containing its winner. Call a colour \emph{represented} if it belongs
to \(\sigma(C)\), and \emph{fresh} otherwise. For a represented colour
to count, both its winner and its clique representative must retain it,
and none of the winner's neighbours in \(C\) may have been assigned that
colour initially. A fresh colour counts whenever its winner retains it.
Each colour counted in this way supplies a colour for the Hall condition.

For this comparison, write \(r_i,f_i\) for the counted represented and
fresh successes at the \(i\)th centre. These represented counts are lower
bounds for the counts in Lemma~\ref{lem:clique-extension}. Take \(T_i\)
there to be the fresh colours counted here, and put
\(J=\sum_i(r_i+f_i-1)_+\). Call a star \emph{strong} when \(r_i\ge2\)
and \emph{weak} when \(r_i\ge1\). A strong star is therefore also weak.
Each fresh colour has only one winner, so the sets of counted fresh
colours at different centres are disjoint. By Lemma~\ref{lem:clique-extension},
\(F_C\) can occur only in the union
\begin{equation}\label{eq:F-obstruction}
 \{\text{every }r_i\le1,\ J\le p\}\ \cup\
 \{|\{i:r_i\ge1\}|\le1\}.
\end{equation}
The first event rules out both alternatives in \eqref{eq:pilot-certificate}.
The second says that at most one star is weak. Under any given clique
exposure, the probability of the union in \eqref{eq:F-obstruction}
therefore bounds the probability of \(F_C\).

\subsection{A product comparison}
To bound \eqref{eq:F-obstruction}, we use products of nonnegative
functions \(g_i(r_i,f_i)\), one for each star. Each function is
nonincreasing in both success counts. We compare the expected product
with one in which different stars face independent blockers. Here a
blocker is a vertex whose initial colour and priority prevent another
vertex from retaining its colour.

Suppose that in any realisation each independent random coordinate can
supply a blocker to at most one of several groups. For a product of
nonnegative functions that increase as blockers are added, giving each
group its own independent copy of every coordinate can only increase
the expectation. To see this, fix all but one coordinate and write each
factor as \(a_i+\delta_i\), where \(a_i\ge0\) is its value when this
coordinate supplies no blocker and \(\delta_i\ge0\). At most one
increment is nonzero. Independent copies preserve the expected constant
and linear terms in the product and add only nonnegative terms of higher
order. Now replace the coordinates one at a time.

In this comparison, an exterior vertex lies outside \(C\) and the
selected endpoints. Unused vertices of \(E\) are included.
Partition the endpoints by initial colour and distinguish represented from
fresh classes. Given these unlabelled classes, the represented labels are
sampled uniformly without replacement from \(\sigma(C)\). For the exterior
comparison, fix the clique exposure \(\mathcal C\), the actual endpoint
colours and priorities, and the minimum priority of each clique colour
class. All exterior vertices retain their independent original marks.
Each target colour has a unique winner, so the stars have disjoint sets
of target colours. An exterior vertex can block only its own colour.
The preceding comparison therefore bounds the expected product by the
product of the expectations for the individual stars.

Within one star there are at most \(d\) target colours. For each exterior
vertex and target colour, the blocker coordinate is nonzero if the vertex
has that colour and priority below \(t\). Here \(t\) is the largest priority
among its neighbours that are the winner or, for a represented colour, the
clique representative of that colour, with \(t=0\) if there are none.
Conditional on
preceding blocker coordinates being zero, a new coordinate with threshold
\(0\le t\le1\) has probability at most \(t/(q-d+1)\); after a nonzero
coordinate it has probability zero. Since
\[
 1-e^{-t/(q-d)}\ge\frac{t}{q-d+1},
\]
we can use independent uniform random variables to generate blockers of
probability \(1-e^{-t/(q-d)}\) that dominate the corresponding blockers
within a star. Together with the comparison between stars, this gives
an upper bound on the product expectation in terms of separate
experiments for the individual stars.

One dependence remains. The label assigned to a represented endpoint
class determines the size of its clique class, which in turn determines
the law of the representative's minimum priority. We express that
minimum using a uniform variable independent of the choice of label.

Consider a represented class of \(\ell\) selected endpoints. Let \(X_\ell\)
be their minimum priority. For a clique colour class of size \(s\), with
minimum priority \(M_s\), the variable
\(V=1-(1-M_s)^s\) is uniform on \((0,1)\), independent of its relative order,
and \(M_s\le V\). Independent transformed uniforms can therefore be assigned
to represented endpoint classes before their labels are injected into the
clique palette. The actual clique minimum still depends on the size of
the class to which the label is assigned.

Let \(w\) be the fixed winner of the endpoint class, and let \(Q_w\)
consist of \(w\) together with the other endpoints in its class whose
priorities are below \(V\). Put \(D_\ell=|Q_w|\).
Avoiding the initial clique colours appearing on neighbours of \(Q_w\)
makes both the winner and its clique representative free of exposed
same-colour blockers.

Fix the clique exposure \(\mathcal C\), the unlabelled endpoint classes
and their palette types, the order and numerical values of the priorities
in each endpoint class, and the transformed minima assigned to the
represented classes. Also fix an independent uniform ordering of those
classes. Call these data \(\mathcal D\). At step \(j\), condition on
\(\mathcal D\) and all label draws and comparison uniforms from steps
\(1,\ldots,j-1\). The current label is uniform among the \(O-j+1\)
remaining clique colours.

Let \(a_w=|N(w)\cap C|\) and
\(A_\circ=\sum_{z\in Q_w\setminus\{w\}}|N(z)\cap C|\).
At most \(a_w+A_\circ\) labels occur on clique neighbours of \(Q_w\).
Write \(b_j\) for the fraction of the remaining labels that are forbidden
in this way, and put
\(\bar b_j=\min\{1,(a_w+A_\circ)/(O-j+1)\}\).
An independent uniform \(U_j\) chooses a forbidden label if
\(U_j\le b_j\), and an allowed one otherwise. Choose uniformly within
the resulting category. Thus \(U_j>\bar b_j\) guarantees an allowed
label after every history.

Put \(Q=q-d\). For each allowed label, the exterior comparison above
gives a conditional survival probability at least
\(s_j=\exp(-(D-a_w)X_\ell/Q-\kappa V)\), where
\(\kappa\ge(p_*+1)/Q\) and \(p_*\ge p\). Indeed, the winner has at most
\(D-a_w\) exterior neighbours with threshold \(X_\ell\), and its
clique representative has at most \(p+1\) with threshold at most \(V\).
At a shared neighbour the maximum threshold is at most their sum.
In the comparison experiment the exterior coordinates have already
been separated by target colour. We can therefore generate survival
for this class with a fresh uniform \(W_j\), independently of the
label draw, using the survival probability for the chosen label as
threshold. Mark this class if \(U_j>\bar b_j\) and \(W_j\le s_j\).
A mark guarantees success in the exterior comparison experiment.
Conditional on the history at step \(j\), its probability is
\begin{equation}\label{eq:joint-floor}
 (1-\bar b_j)s_j
 =e^{-(D-a_w)X_\ell/Q-\kappa V}
 \left(1-\frac{a_w+A_\circ}{O-j+1}\right)_+ .
\end{equation}
Given \(\mathcal D\), both thresholds are fixed before we draw any labels.
The independent pairs \((U_j,W_j)\) thus give independent marks, each
sufficient for success in the comparison experiment, even though the
forbidden fraction \(b_j\) can depend on earlier draws. Counting only
these marks can only increase the expected failure products.

When \(O\le Q\), the function
\[
 a\longmapsto e^{aX_\ell/Q}
       \left(1-\frac{a+A_\circ}{O-j+1}\right)_+
\]
is nonincreasing for \(0\le X_\ell\le1\). On its positive part, its
derivative has the sign of
\(X_\ell(O-j+1-a-A_\circ)/Q-1\le0\); the clipped function is continuous.
Increasing any other attachment also decreases the expression.
If every endpoint has at most \(P\) clique neighbours, replacing
\(a_w\) by \(P\) and \(a_w+A_\circ\) by \(P D_\ell\) lowers the mark
probability. The bound now depends only on the class priorities and its
rank \(j\). Averaging over those priorities at fixed class size and rank gives
\begin{align}
 \Gamma_{\ell,j}
 &=\E\left[e^{-\alpha_R X_\ell-\kappa V}
       \left(1-\frac{P D_\ell}{O-j+1}\right)_+\right],
 &\alpha_R&=\frac{D-P}{q-d}.\label{eq:Gamma}
\end{align}
We call \((D-P)/(q-d)\) the compensated rate, since it accounts for
the winner's clique neighbours when bounding the exterior blockers.
For general \(O\), \(a_w\ge1\) gives the uncompensated rate
\(\alpha_R=q/(q-d)\). In both cases the factor involving \(P D_\ell\)
is the same. For fresh winners, put \(\alpha_F=q/(q-d)\). Their
probability of success is at least
\begin{equation}\label{eq:fresh-floor}
 \Phi_\ell=J_\ell(\alpha_F),\qquad
 J_\ell(t)=\ell\int_0^1e^{-tx}(1-x)^{\ell-1}\,dx.
\end{equation}

Conditional on the endpoint partitions and rank order, the numerical
order statistics are independent across classes, and the transformed
clique uniforms are independent of them. The averaging in
\eqref{eq:Gamma} and \eqref{eq:fresh-floor} therefore preserves
independence of the auxiliary marks. Any deterministic lower bounds
\[
 0\le\gamma_{\ell,j}\le\Gamma_{\ell,j},\qquad
 0\le\phi_\ell\le\Phi_\ell
\]
can therefore be used. Lowering a success probability can only increase
a failure product. In the auxiliary experiment, each class places a
Bernoulli mark at its uniformly chosen winning endpoint position, with
success probability given by its size and, for a represented class, its
rank in the label order. These marks describe successes in the comparison
experiment. The expectation inequality persists when we average over the
partitions and the rank order, giving the following result.

\begin{lemma}\label{lem:product}
Conditional on \(\mathcal C\), the expectation of every product of
nonnegative local functions \(g_i(r_i,f_i)\) that are nonincreasing in
both counts is at most its expectation in the auxiliary class-mark
experiment above. The same holds for nonnegative sums of these products.
The selected family and the functions are fixed at the exposed clique
information.
\end{lemma}

We apply Lemma~\ref{lem:product} by expressing bounds for the two events
in \eqref{eq:F-obstruction} as nonnegative sums of the products in the
lemma. The same comparison applies if a fixed deletion leaves some groups
with fewer than \(d\) endpoints. Each group still has at most \(d\) target
colours, and the colour laws use the actual number of remaining endpoints.

\subsection{From the comparison to a failure curve}
For the computations below, a \emph{colour profile} is a fixed choice of
the success bounds \(\gamma_{\ell,j},\phi_\ell\), together with the survival
parameters and evaluation ranges used to compute them. We also refer to
these success bounds as \emph{floors}.
To find the distribution of the marks we first count the occupied colour
classes. If \(S\) is the number of the \(n\) endpoint positions that
receive represented colours, then
\[
 \PP(S=s)=\binom ns(O/q)^s((q-O)/q)^{n-s}.
\]
Given \(S=s\), the probability that exactly \(r\) represented classes
are occupied is
\[
 \frac{(O)_r{s\brace r}}{O^s},\qquad 0\le r\le\min(s,O),
\]
where \({s\brace r}\) is a Stirling number of the second kind and
\((a)_r=a(a-1)\cdots(a-r+1)\), with \((a)_0=1\). For nonnegative
integers \(n\) we take \(\binom nk=0\) when \(k<0\) or \(k>n\).
If \(s=r=0\), there are no represented marks. Given \(s\) and \(r\ge1\), a uniformly chosen occupied represented
class has size \(\ell\) with probability
\[
 \frac{\binom s\ell{s-\ell\brace r-1}}{r{s\brace r}}.
\]
In a uniform order of the \(r\) classes, the last class has this size
distribution. Using \(\gamma_{\ell,r}\) for its mark and then removing it
gives a recurrence with nonnegative terms for the number of marks. There may be more than \(O\) represented
positions, but at most \(O\) distinct represented classes.

For clique deficits \(p\ge3\), the represented success bounds depend on the
injection rank. In these profiles, marks are set to zero above class size four,
but those classes still consume labels and ranks. Omitting represented
partitions containing a class of size greater than \(24\) incurs the error
\[
 \tau_{n,O}=\min\left\{1,\frac{O\binom n{25}}{q^{25}}\right\}.
\]
The fresh recurrence is evaluated for at most \(q-O+1\) fresh positions.
The probability of larger fresh counts is included in a separate term
\(\nu_{n,O}\). Appendix~\ref{app:colour-laws} specifies both recurrences,
the positive remainders for their truncated sums, and their endpoint
conventions.
For deficits one and two we use a more conservative two-star profile whose
represented success bounds are independent of the label order. The two
recurrences in Appendix~\ref{app:constant-rank-laws} give the joint
distribution of the represented and fresh marks and the represented
marginal needed to bound the probability of at most one weak star.

We now use these distributions to bound the two Hall obstructions.
Take \(K\ge2\) stars, each with \(d\) endpoints, and put
\(L=(d-1)K\). On \(r_i\le1\), define
\[
 B_{i,t}=\binom{d-1-(r_i+f_i-1)_+}{t}\quad(0\le t\le d-1),
\]
and set \(B_{i,t}=0\) when \(r_i\ge2\).
These functions are nonnegative and nonincreasing in both success counts.
The coefficient of \(z^k\) in
\(\prod_i\sum_t B_{i,t}z^t\) is
\[
 \ind_{\{\forall i:r_i\le1\}}\binom{L-J}{k}.
\]
On the first event in \eqref{eq:F-obstruction}, this is at least
\(\binom{L-p}{k}\), for \(0\le k\le L-p\).
Lemma~\ref{lem:product} therefore bounds the probability of the first
event by the expected auxiliary coefficient divided by \(\binom{L-p}{k}\).
Every permitted \(k\) bounds the same probability, so we may take the
minimum. We use this factorial bound when \(p<L\). When \(p\ge L\), we
bound the first event by the probability that no star is strong.

To bound the probability that no star is strong, or that at most one
star is weak, it suffices to count placements of represented marks.
After averaging over the colour partitions and the independent rank
order, the \(b\) marked represented positions form a uniform
\(b\)-subset of the \(n=dK\) positions, conditional on \(b\).
Let \(\widetilde r_b\) be the probability that there are \(b\)
represented marks and that the counts and partition lie in the retained
ranges. The omitted mass is bounded by \(\nu_{n,O}+\tau_{n,O}\).
On the retained states the two bounds are
\begin{align}
 NS_K&=\sum_{b=0}^K\widetilde r_b
              \frac{d^b\binom Kb}{\binom nb},\label{eq:NS}\\
 WU_{K,d}&=\sum_{b=0}^d\widetilde r_b
              \frac{K\binom db}{\binom nb}.\label{eq:WU}
\end{align}
In the first sum the \(b\) marks must lie in distinct stars. In the
second they must all lie in a single star, and we sum over its \(K\)
possible choices. Thus the zero-mark term has coefficient \(K\), and
terms with \(K<b\le d\) still contribute.

Let \(M_{K,p}\) be the factorial bound for \(p<L\), and omit it from the
minimum when \(p\ge L\). A complete conditional bound is
\begin{equation}\label{eq:complete-U}
 U_{p,O}(K)=\min\{1,\nu_{n,O}+\tau_{n,O}
                 +\min(NS_K,M_{K,p})+WU_{K,d}\}.
\end{equation}
A further valid exponential-product bound for the first obstruction can
also be included in the inner minimum. It has the form
\(z^{-p}\E_{\rm aux}[\ind_{\{\forall i:r_i\le1\}}z^J]\), \(0<z<1\).
Here \(\nu_{n,O}\) and \(\tau_{n,O}\) bound the omitted states,
\(WU_{K,d}\) bounds the second Hall obstruction, and the inner minimum
compares bounds for the first obstruction. The oversized-class event is
removed before the two obstruction bounds are applied.

For a fixed profile with prefix cutoff \(B\), define
\begin{equation}\label{eq:psi}
 \psi_{p,O}(r)=\min_{0\le j\le\min(r,B)}U_{p,O}(j),
 \qquad U_{p,O}(0)=U_{p,O}(1)=1.
\end{equation}
If several profiles bound the same conditional failure event under the
same exposure, we take the minimum of their complete bounds. The curve
in \eqref{eq:psi} is nonincreasing and is constant from \(B\) onwards.
Applying the bounds to prefixes of a maximum family in the fixed set
gives, for every permitted exposure,
\[
 \PP(F_C\mid\mathcal C)\le\psi_{p,O}(R_E(A_C)).
\]
Conditional on \(O\) alone, symmetry makes \(A_C\) a uniform
\((m-O)\)-subset of \(C\). Hence
\begin{equation}\label{eq:rank-interface}
 \PP(F_C\mid O)
 \le\E_{A\in\binom C{m-O}}\psi_{p,O}(R_E(A)).
\end{equation}
To remove the conditioning on \(O\) in \eqref{eq:rank-interface},
we use the occupancy probabilities
\begin{equation}\label{eq:occupancy}
 w_{m,O}=\PP(|\sigma(C)|=O)=\frac{(q)_O{m\brace O}}{q^m}.
\end{equation}
For each value of \(O\) outside the range of the estimate, we include
its full probability with the failure bound one.

\subsection{The vertices left for greedy colouring}
Fix \(v\in\mathcal V_A\) and write \(n=d(v)\in\{749,750\}\).
If \(v\) is uncoloured, one saving suffices at degree 749 and two
suffice at degree 750. To find these savings, let
\(Q=\overline{G[N(v)]}\). In each initial neighbourhood colour class
of size at least two, take the first two vertices in priority order. Call their
pair \emph{admissible} if it is an edge of \(Q\). An admissible pair
\emph{survives} if both vertices retain their colour. Each surviving
admissible pair contributes one saving to \(s_\varphi(v)\).

The degree of an endpoint in \(Q\) controls its number of neighbours
outside \(N[v]\). Writing \(\delta_x=d_Q(x)\), put
\[
 t_x=\frac{D-n+\delta_x}{D-2}.
\]
First ignore \(v\) as a possible blocker. We use the symmetric edge
weight \(b_{xy}\) defined in \eqref{src:SG12} to bound the probability
that a selected admissible pair \(xy\in E(Q)\) survives. This probability
is at least \((9/10)b_{xy}\) for a class of size two and at least
\(b_{xy}\) for a larger class. On the range of degrees considered here,
\[
 b_{xy}\ge A_n-\frac{\gamma}{2}(\delta_x+\delta_y),\qquad
 A_n=\frac{49}{50}-\frac{3(750-n)}{5\cdot748},\quad
 \gamma=\frac3{5\cdot748}.
\]
The excluded joins from Section~\ref{sec:clique-parts} restrict the
degrees in \(Q\). If the endpoints of an edge \(xy\) of \(Q\) have at
least three common \(G\)-neighbours in \(N(v)\), those neighbours form a
clique. Otherwise they and \(v\) induce a connected incomplete graph on
at least four vertices joined to \(I_2\), contrary to the second
Cranston--Rabern theorem. A triangle of \(Q\) has at most four common
\(G\)-neighbours in \(N(v)\). Five such neighbours would form a clique
and, together with \(v\), give \(K_6\vee I_3\).
Finally, two disjoint \(Q\)-edges have at most two common \(G\)-neighbours in \(N(v)\).
Three, together with \(v\), would give a \(K_4\) joined to four vertices
whose complement has a matching of size two, again excluded by the
first theorem.

Suppose first that \(c(v)<250\), and put \(a=248\). An independent set
in \(Q\) forms a clique with \(v\) in \(G\), so \(\alpha(Q)\le a\).
The common non-neighbours of \(x,y\) in \(Q\) are their common
\(G\)-neighbours in \(N(v)\). For \(xy\in E(Q)\), these form an
independent set in \(Q\) whenever there are at least three. There are
therefore at most \(a\) of them, and inclusion--exclusion gives
\begin{equation}\label{src:SG14}
 \delta_x+\delta_y\ge n-a+|N_Q(x)\cap N_Q(y)|.
\end{equation}
For a three-element set \(X\), consider the vertices outside \(X\)
adjacent in \(Q\) to none of its members. The graph they induce in \(Q\)
has matching number at most one by the restriction on disjoint edges.
It consists of a star or a triangle, together with isolated vertices
(the edgeless case is included). In the star case all but at most one
vertex are independent, so the order is at most \(a+1\). In the triangle
case \(X\) and the isolated vertices are all common \(G\)-neighbours of
the triangle, so the preceding triangle bound again gives order at most
\(a+1\). Consequently
\begin{equation}\label{src:SG15}
 \sum_{x\in X}\delta_x\ge n-a-4+2e_Q(X).
\end{equation}

Put \(H=\lceil(n-a)/2\rceil\) and \(L=\lfloor(n-a-4)/3\rfloor\).
Vertices of degree below \(H\) are independent by
\eqref{src:SG14}, so there are at most \(a\). By
\eqref{src:SG15}, at most two vertices have degree below \(L\).
The sorted degrees thus dominate
\[
 (0,0,\underbrace{L,\ldots,L}_{a-2},
       \underbrace{H,\ldots,H}_{n-a}).
\]
Let \(B_1=(a-2)L+(n-a)H\), \(B_2=(a-2)L^2+(n-a)H^2\),
and \(M_n=n-1+H\). Let \(b_i\) be the displayed lower bound for the
\(i\)th smallest degree \(\delta_i\). Then
\(\delta_i^2\le b_i^2+M_n(\delta_i-b_i)\). If \(e=e(Q)\), it follows that
\[
 \sum_i\delta_i^2\le B_2+M_n(2e-B_1),\qquad 2e\ge B_1.
\]
Summing the linear kernel bound, and using
\(A_n-\gamma M_n=6652/(50\cdot748)>0\), gives
\begin{equation}\label{src:SG16}
 \sum_{xy\in E(Q)}b_{xy}\ge
 (A_n-\gamma M_n)e+\frac{\gamma}{2}(M_nB_1-B_2)
 \ge S_n:=\frac{A_nB_1-\gamma B_2}{2}.
\end{equation}
The lower bound \(S_n\) thus applies at every permitted edge density.
The appendix gives its exact value for each of the two root degrees.

When \(250\le c=c(v)\le634\), put \(T=C(v)\setminus\{v\}\),
\(a_C=c-1\), \(S=N(v)\setminus C(v)\), and \(s=n-a_C\).
Two vertices of \(S\) have at most two common neighbours in \(T\).
Otherwise, together with their common neighbour \(v\), the
four-common-neighbour consequence would give a larger clique through \(v\).
If \(L_C=e_G(S,T)\), pair counting followed by Cauchy's inequality gives
\[
 L_C^2-a_CL_C-2a_Cs(s-1)\le0.
\]
Each vertex of \(S\) also misses a vertex of \(T\), because \(C(v)\)
is a largest clique through \(v\). Thus the number of missing cross edges
is at least
\[
 E(n,c)=a_Cs-\min\left\{s(a_C-1),
 \left\lfloor\frac{a_C+\sqrt{a_C^2+8a_Cs(s-1)}}2\right\rfloor\right\}.
\]
Select \(E(n,c)\) of these missing cross edges in the fixed graph order.
In \(Q\) we give each selected edge the lower weight \(b\) from
\eqref{src:SG17} and every other edge weight zero. The sum of the
weighted degrees on each side is then \(bE(n,c)\).

Choose \((e_1,\ldots,e_r)\) uniformly from all ordered sequences of
\(r\) disjoint unordered pairs of the \(n\) neighbours of \(v\). A pair
is tested for membership in \(Q\) only after it is chosen. Write
\[
 Z_r=\E\prod_{j=1}^r(1-b_{e_j}\ind_{\{e_j\in Q\}}).
\]
It remains to account for the root \(v\), which can both lose its own
colour and prevent a pair from retaining that colour. Temporarily remove
it as a blocker. For each neighbourhood colour class \(\gamma\), let
\(J_\gamma\) indicate that its selected admissible pair survives with
the root removed. Put \(J_\gamma=0\) if there is no such pair. Write
\(S=\sum_\gamma J_\gamma\), let \(B_v\in\{0,1\}\) count the successes
lost when we restore the root, and let \(I_v\) indicate that the root
is uncoloured.
For \(z_*=2843/62500\),
\[
 \PP(A_v)\le z_*^{-(n-q)}\E[z_*^{S-B_v}I_v].
\]
Since \(B_v\in\{0,1\}\), expanding
\(z_*^{-B_v}=1+(z_*^{-1}-1)B_v\) gives a first term
\(\E[z_*^S I_v]\) and a correction proportional to
\(\E[z_*^S B_v I_v]\). The correction accounts for a pair whose
success is lost even though the restored root is itself uncoloured.

Conditional on the marks away from the root, the probability that the
root receives a given occupied colour and loses it is \((1-U_1)/q\),
where \(U_1\) is the minimum priority in that colour class.
If the class has at least two vertices, let \(U_2\) be its second-smallest
priority. For a pair that succeeds with the root removed, restoring the root
destroys that success while leaving the root uncoloured precisely when
the root receives this class colour and has priority between \(U_1\)
and \(U_2\). The conditional loss probability is therefore
\(J_\gamma(U_2-U_1)/q\). For a class of size \(k\), the mean of
\(1-U_1\) is \(k/(k+1)\), and for \(k\ge2\) the mean of \(U_2-U_1\)
is \(1/(k+1)\). We first apply the exterior product comparison and
then average over the priorities. In the first term, the pair-failure
factor increases with the pair priorities and \(1-U_1\) decreases, so
the product-space covariance inequality bounds the mean of their product
by the product of their means. In the correction term we bound the
distinguished pair's survival indicator by its admissibility indicator
before averaging the gap \(U_2-U_1\). This gives a nonnegative sum of
matching expectations. The distinguished pair is tested for admissibility
in the full \(Q\), with weight one if it passes. Dropping this test bounds
its matching expectation by the corresponding \(Z_r\). Hence
\[
 \PP(A_v)\le\sum_{r=0}^{\lfloor n/2\rfloor}\Gamma^A_{n,r}Z_r,
 \qquad \Gamma^A_{n,r}\ge0.
\]
The coefficients \(\Gamma^A_{n,r}\) collect the colour-class contributions
from the first term and the correction, after averaging over class sizes
and assignments to the neighbourhood labels. Both terms use matchings on
all \(n\) labels. The appendix derives these coefficients and their
truncation bounds, retaining the joint dependence of the root and the
surviving repetitions.

The matching quantities admit bounds from the binary Brégman inequality
\cite{Bregman,AF}. If
\(\beta(j)=(j!)^{1/j}\) for positive integers \(j\), \(\beta(0)=0\), and
\(\beta\) is linearly interpolated, then for a matrix with entries in
\([0,1]\) and row sums \(r_i\),
\[
 \per A\le\prod_i\beta(r_i).
\]
For the weighted extension, maximise the permanent row by row while
preserving row sums until each row has at most one fractional entry.
Independent rounding of those entries preserves the expected permanent,
and the binary inequality gives the displayed bound.
The square of a weighted perfect-matching sum is at most the permanent
of its symmetric adjacency matrix: the union of two perfect matchings
decomposes into doubled edges and alternating even cycles.
Completing a matching by labelled dummy vertices therefore converts
the bounds on total or cross kernel mass into bounds on \(Z_r\).

Appendix~\ref{app:vertices} gives the normalisation for the matchings,
the coefficients for both root terms and the bounds on the errors from
truncating the class-size sums. Together they give the upper bound
\(U_{n,c}\) used in Section~\ref{sec:closure}. For each of the two root
degrees we make one comparison using the total weight when \(c<250\),
and one using the cross edges for each of the \(385\) values
\(250\le c\le634\). By \eqref{eq:root-clique}, these cases cover every
possible \(c\).

\section{Star supply and the available packing size}
\label{sec:rank}

The failure bound in \eqref{eq:rank-interface} decreases as more stars
become available. We first use the graph structure to find disjoint
stars before any colours are assigned. Once the available centres are
known, we can keep the stars whose centres are available and look for
additional stars using the same fixed endpoint set.

\subsection{A paired matching formula}
Let \(H\) be a simple bipartite graph with centre side \(C\) and endpoint
side \(E\), and let \(r_2(H)\) be its maximum number of vertex-disjoint
two-stars centred in \(C\). For \(X\subseteq C\) and \(Y\subseteq E\),
let \(o_E(H-X-Y)\) count components with an odd number of endpoint
vertices. The Tutte--Berge formula gives
\begin{equation}\label{eq:paired}
 2r_2(H)=\min_{X\subseteq C,\ Y\subseteq E}
       \bigl(|E|+2|X|+|Y|-o_E(H-X-Y)\bigr).
\end{equation}
To prove this, replace each centre by two adjacent copies, with the same
endpoint neighbours. A packing of \(t\) two-stars gives a matching of size
\(|C|+t\), using the internal edge at every other centre. Conversely,
a matching with \(t\) centres using two external edges and \(w\) centres
using no edge has size \(|C|+t-w\le|C|+r_2(H)\).

In the graph of twin vertices, choose a deletion set that minimises the
Tutte--Berge expression and, among these, has smallest cardinality. It
contains both copies of a centre or neither. Indeed, if just one copy
were deleted, restoring it would enlarge its twin's component by one
vertex without merging any components. The expression could not increase,
contrary to our choice of a smallest minimiser. Thus the chosen set
deletes the pairs of copies of some \(X\subseteq C\) together with some
\(Y\subseteq E\). In each remaining component the number of centre
copies is even, so its order has the parity of its endpoint count. This
proves \eqref{eq:paired}, including isolated vertices and empty sides.

\begingroup\samepage
We use \eqref{eq:paired} in two degree estimates. In the application
of the first, \(p=D-|C|\) is the integer clique deficit.
\begin{lemma}\label{lem:codegree-supply}
Suppose the centre degrees \(k_c\) satisfy \(p\le k_c\le U\), where
\(p\ge2\), every endpoint has degree at most \(P\), and every pair of
endpoints has at most \(\rho\ge1\) common centres. If \(L\ge3\) is odd,
\(\theta=(L-1)/L\), and
\[
 \theta U\le2P,\qquad
 \rho(L-1)(L-2)\le P(p-1)L,
\]
then
\[
 r_2(H)\ge
 \left\lceil\frac{\theta}{2P}\sum_{c\in C}k_c\right\rceil
 \ge\left\lceil\frac{\theta p|C|}{2P}\right\rceil.
\]
\end{lemma}
\endgroup

\begin{proof}
Choose \(X,Y\) attaining \eqref{eq:paired}.
In a component \(H_i\) of \(H-X-Y\), let \(t_i\) be its endpoint count,
let \(r_c\) be a centre's degree into that component, and put
\[
 g_i=\sum_{c\in C(H_i)}(2r_c-k_c),\qquad
 \epsilon_i=\ind_{\{t_i\text{ odd}\}}.
\]
We claim that \(\theta g_i\le P(t_i-\epsilon_i)\).
For even \(t_i\), use \(g_i\le\sum_c r_c\le Pt_i\).
For odd \(t_i\), we have both
\[
 g_i\le Pt_i,\qquad
 g_i\le\frac{\rho t_i(t_i-1)}{p-1}.
\]
For the second bound, use
\((2r-k)_+\le r(r-1)/(p-1)\) for \(k\ge p\) and \(0\le r\le k\),
then count endpoint pairs. When \(r\le p\), the required difference is
\((r-p)(r-p+1)\ge0\); when \(r\ge p\), use \(2r-k\le r\).
If \(t_i\le L-2\), the second bound proves the claim, since
\(\theta\rho t_i\le P(p-1)\). If \(t_i\ge L\), the first does, since
\(\theta t_i\le t_i-1\). These exhaust odd \(t_i\).

Let \(e_Y\) count incidences from \(C\setminus X\) to \(Y\). Then
\[
 \sum_{c\in C}k_c=\sum_{c\in X}k_c+2e_Y+\sum_i g_i.
\]
Multiplying by \(\theta\), and using the degree bounds and the claim, gives
\[
 \theta\sum_c k_c
 \le2P|X|+2P|Y|+P\sum_i(t_i-\epsilon_i)
 =2P r_2(H).
\]
Integrality finishes the proof.
\end{proof}

\begin{lemma}\label{lem:irregular}
If every endpoint has degree at most \(b\ge1\) and every centre has degree
\(k_c\le2b+1\), then
\[
 r_2(H)\ge
 \left\lceil\frac{\sum_c(k_c-1)_+}{2b}\right\rceil.
\]
\end{lemma}
\begin{proof}
Delete all centres of degree at most one. None can be the centre of a
two-star, and each contributes zero to the sum in the bound. Choose
\(X,Y\) as in \eqref{eq:paired}.
For a component \(H_i\), let \(a_i,t_i,e_i\) count its centres, endpoints
and edges. As above, let \(\epsilon_i\) indicate that \(t_i\) is odd and
let \(e_Y\) count edges from \(C\setminus X\) to \(Y\). Then
\[
 e_i-a_i\le b(t_i-\epsilon_i).
\]
For even \(t_i\), this follows from \(e_i\le bt_i\).
For odd \(t_i\) and \(a_i\ge b\), subtract \(a_i\) from the same inequality.
For odd \(t_i\) and \(a_i\le b\), simplicity gives
\(e_i-a_i\le a_i(t_i-1)\le b(t_i-1)\).
Partitioning the original edge surplus now gives
\[
 \sum_c(k_c-1)
 =\sum_{c\in X}(k_c-1)+e_Y+\sum_i(e_i-a_i)
 \le2b\,r_2(H)-b|Y|\le2b\,r_2(H).\qedhere
\]
\end{proof}

Lemma~\ref{lem:codegree-supply} uses the common-neighbour bound to
guarantee stars in the original graph. Once we sample the available
centres and cap degrees at the endpoints, the remaining centre degrees
can vary. Lemma~\ref{lem:irregular} applies to that residual graph.

\subsection{Recovering stars in a fixed endpoint set}
Fix a family of \(T\) disjoint \(d\)-stars, \(1\le T\le m=|C|\), with
centre set \(I\) and original endpoint set \(E_0\). We use either \(E_0\)
itself or the set \(E\) obtained by adjoining to \(E_0\) all good neighbours
of unused centres. If the original family is maximal in the full good-incidence
graph, each unused centre has at most \(d-1\) good neighbours outside \(E_0\).
Consequently the enlarged set has
\begin{equation}\label{eq:expanded-size}
 |E|\le dT+(d-1)(m-T)=(d-1)m+T.
\end{equation}

Label each original endpoint by its centre and give each added endpoint
a fixed label from \(I\). We call the label of \(e\) its \emph{owner}
and denote it by \(\iota(e)\). An added endpoint can receive any label
from \(I\), whether or not it is adjacent to that centre. We will use
the labels when averaging over the available centres.

Let \(A\) be a uniform \(h\)-subset of \(C\). Keep the
\(Z=|A\cap I|\) available original stars. We look for additional stars at
centres in \(A\setminus I\), using only endpoints whose owners are absent
from \(A\). This endpoint set is disjoint from the kept stars.
Write \(N=|A\setminus I|\), and let \(\ell_e\) be the number of
neighbours of endpoint \(e\) in \(A\setminus I\). Put
\[
 Y_b=\sum_{\iota(e)\notin A}\min(\ell_e,b),
\]
where \(b\ge1\) is fixed. Retain exactly \(\min(\ell_e,b)\) incidences at
each such endpoint, breaking ties from the graph.

Suppose centre degrees in the fixed pool are at most \(U\).
A maximal residual \(d\)-star packing of size \(t\) gives
\[
 Y_b\le Ut+dt(b-1)+(d-1)(N-t)
       =(d-1)N+(U-2d+1+db)t.
\]
The three terms count edges at selected centres, edges from unselected
centres to selected endpoints, and edges between unselected centres and
endpoints. The last count uses maximality of this residual packing.
For \(d=2\) and \(U\le2b+1\), Lemma~\ref{lem:irregular} gives the sharper
denominator \(2b\), by choosing a maximum residual two-star packing.
Thus, with
\[
 \ell_*=
 \begin{cases}
 2b,&d=2,\ U\le2b+1,\\
 U-2d+1+db,&\text{otherwise},
 \end{cases}
\]
the actual maximum \(R_E(A)\) satisfies
\begin{equation}\label{eq:residual-rank}
 R_E(A)\ge Z+\frac{(Y_b-(d-1)N)_+}{\ell_*}.
\end{equation}
Keeping the \(Z\) original stars and adding a residual packing gives
at least this many stars. A maximum packing may use a different family.

Let \(P\) bound the endpoint degrees, and let every nonoriginal centre
have degree at least \(a_0\) in the same pool. The number of nonoriginal
neighbours of an endpoint is bounded by
\[
 s=
 \begin{cases}
 P-1,&E=E_0,\\
 P,&\text{in the enlarged pool}.
 \end{cases}
\]
An original endpoint already has one neighbour in \(I\), its owner. An
added endpoint need not be adjacent to its assigned owner, so for it we
use the full degree bound \(P\). In the applications \(1\le s\le m-1\).
Write \(\Hypergeom(M,h,t)\) for the size of the intersection of a uniform
\(h\)-subset of an \(M\)-set with a fixed \(t\)-subset. For \(h<m\), define
\[
 f_h(t)=\E\min\{\Hypergeom(m-1,h,t),b\}.
\]
\[
 \beta_h=
 \frac{\bigl[(1-h/m)(a_0/s)f_h(s)-(d-1)h/m\bigr]_+}{\ell_*},
\]
and put \(\beta_m=0\). We claim that
\begin{equation}\label{eq:affine-mean}
 \E R_E(A)\ge\mu_h(T):=\frac{hT}{m}+(m-T)\beta_h.
\end{equation}
To prove this, note that \(f_h(0)=0\) and \(f_h\) is discretely
concave. For \(h>0\) its next increment is
\[
 \frac h{m-1}
 \PP\{\Hypergeom(m-2,h-1,t)\le b-1\},
\]
which decreases with \(t\). For \(h=0\) the function is zero.
Therefore \(f_h(r)\ge(r/s)f_h(s)\) for \(r\le s\).
Each owner's label lies outside the set of nonoriginal neighbours.
Conditional on that owner being absent, the sample is a uniform
\(h\)-subset of the other \(m-1\) labels. Linearity of expectation gives
\[
 \E Y_b\ge(1-h/m)\frac{a_0(m-T)}s f_h(s).
\]
Apply the convexity of the positive part in \eqref{eq:residual-rank},
together with \(\E Z=hT/m\) and \(\E N=h(m-T)/m\), to obtain
\eqref{eq:affine-mean}. No independence between owners or endpoints is
needed.

Since \(f_h(s)\le hs/(m-1)\), for \(1\le h<m\) we have
\[
 0\le\beta_h\le\frac hm\frac{(a_0-d+1)_+}{\ell_*}\le\frac hm.
\]
The last inequality follows from \(a_0\le U\) and the choice of
\(\ell_*\). Thus \(\mu_h(T)\) is affine and nondecreasing in \(T\)
when its other parameters are fixed throughout the interval. The
boundary cases follow directly. If \(h=0\), the rank is zero. If
\(h=m\), every original star is available, so \(R_E(C)\ge T\) without
conditioning on an absent owner. If \(T=m\), every available centre
has its original star and \(R_E(A)=h\).

If \(hT/m+(m-T)\gamma_h\) is a second lower bound for the same
\(\E R_E(A)\) throughout a fixed interval in \([0,m]\), then
\(m-T\ge0\) allows us to replace \(\beta_h\) in
\eqref{eq:affine-mean} by \(\max\{\beta_h,\gamma_h\}\).
The small-deficit and original-pool constructions in Appendix~\ref{app:small}
combine lower bounds in this way.

\subsection{From the mean to an interval of probability bounds}\label{sec:rank-envelopes}
Exchanging one available centre changes \(R_E(A)\) by at most one.
Deleting the centre from a maximum packing loses at most one star.
Interchanging the roles of the two centre sets gives the reverse
inequality. In both comparisons the endpoint set is fixed.

Put \(k_h=\min(h,m-h)\) and
\[
 V_h=
 \begin{cases}
 0,&k_h=0,\\
 (m-k_h)^2\displaystyle\sum_{j=m-k_h}^{m-1}j^{-2},&k_h>0.
 \end{cases}
\]
Of \(A\) and \(C\setminus A\), reveal the one of size \(k_h\) in a
uniformly random order. After \(i\) labels have been revealed, compare
two choices for the next label by transposing them in the two
completions. The completed sets agree unless the other label is absent
from the rest of the sample. This has probability
\((m-k_h)/(m-i-1)\), and in that case the sets differ by one exchange.
The conditional means for the next step therefore lie in an interval of
length at most this probability. Applying Hoeffding's scalar bound
\cite[equation~(4.16)]{Hoeffding} to each centred increment gives
\begin{equation}\label{eq:rank-moment}
 \E e^{-\lambda R_E(A)}
 \le\exp\left(-\lambda\mu_h(T)+\frac{\lambda^2V_h}{8}\right),
 \qquad\lambda\ge0.
\end{equation}
When the sample is deterministic, \(V_h=0\).

Fix an interval \([T_0,T_1]\). For each \(h\), use the same curve
\(\psi_h\) at occupancy \(O=m-h\), cutoff \(B\), and coefficients
\(A_h,C_h,\lambda_h\) throughout the interval. Choose these coefficients
so that
\[
 A_h\ge\psi_h(B),\quad C_h\ge0,\quad\lambda_h\ge0,\qquad
 C_h\ge\max_{0\le j<B}
           (\psi_h(j)-A_h)e^{\lambda_hj}.
\]
Here the maximum is taken together with zero and equals zero if there
is no index \(j<B\). The resulting bound
\(\psi_h(j)\le A_h+C_he^{-\lambda_hj}\) and
\eqref{eq:rank-moment} give
\begin{equation}\label{eq:outer-envelope}
 \PP(F_C\mid O=m-h)
 \le A_h+C_h\exp\left(-\lambda_h\mu_h(T)
                          +\frac{\lambda_h^2V_h}{8}\right).
\end{equation}
We may take \(A_h\) larger than the constant value \(\psi_h(B)\). For
each smaller prefix the displayed condition on \(C_h\) can be checked
directly.

Average \eqref{eq:outer-envelope} with the weights in
\eqref{eq:occupancy}, assigning failure probability at most one to each
unsupported occupancy. Write \(P(T)\) for the bound before taking its
minimum with one. Since \(\mu_h(T)\) is affine, \(P(T)\) is a
nonnegative sum of exponentials in \(T\). If the upper bound \(e(T)\)
on the size of the endpoint set is also affine, the allowance in
\eqref{eq:F-target} has the form \(x_0e^{-L_0-L_1T}\). The ratio
\(Q(T)=P(T)/(x_0e^{-L_0-L_1T})\) is again a nonnegative sum of
exponentials. When the sum is positive, \((\log Q)''\) is the
variance of their slopes under the weights proportional to the terms of
the sum. Thus \(Q\) is log-convex and
\begin{equation}\label{eq:endpoints}
 \frac{P(T)}{x_0e^{-L_0-L_1T}}
 \le
 \max\left\{\frac{P(T_0)}{x_0e^{-L_0-L_1T_0}},
            \frac{P(T_1)}{x_0e^{-L_0-L_1T_1}}\right\}
 \quad(T_0\le T\le T_1).
\end{equation}
The endpoint inequality also holds for zero sums and singleton intervals.
It is an analytic inequality throughout \([T_0,T_1]\), irrespective of
which packing sizes are realised by graphs.

\subsection{A pointwise alternative}
For the larger deficits we use a deterministic lower bound on the
available rank. Suppose every two endpoints have at most \(\rho\)
common centres, and every nonoriginal centre has at least \(k\) neighbours
in the fixed pool. Keep the available original stars and extend them
maximally using other available centres. If the resulting family has
\(t\) stars, its endpoint set has size \(dt\). Every unselected centre
has at least \(a=k-d+1\) neighbours there, by maximality.
Each selected centre has its \(d\) assigned endpoints. When \(k\ge d+1\),
pair counting gives
\[
 t\binom d2+(h-t)\binom a2\le\rho\binom{dt}{2},
\]
or equivalently
\begin{equation}\label{eq:rank-polynomial}
 \rho d^2t^2+
 [a(a-1)-d(d+\rho-1)]t-ha(a-1)\ge0.
\end{equation}
For \(h>0\), let \(L(h,k)\) be the least positive integer satisfying
\eqref{eq:rank-polynomial}; put \(L(0,k)=0\).
The polynomial is negative at zero and nonnegative at \(h\), so this
integer exists. The actual maximum rank is at least \(L(h,k)\), as well
as the number of available originals.

For a maximal family in the full good-incidence graph, maximality gives
each unused centre at least \(p-d+1\) neighbours in the original
endpoint set. Thus, when we look for additional stars within that set,
the pair count uses \(a=p-2d+2\). When we count the initial stars
in the full graph, all good neighbours are available and we use
\(a=p-d+1\).

When only a fixed family of \(T_*\) stars is retained, its number of
available centres has the exact law
\[
 Z\sim\Hypergeom(m,h,T_*).
\]
Averaging \(\psi_{p,O}(Z)\) against this full distribution gives the
bound directly. We use this calculation when the original packing is
large enough to retain a fixed number of stars. The retained family
need not be maximal.

\section{The finite ranges and completion of the proof}
\label{sec:closure}

It remains to choose the endpoint set for each designated clique and
check the Local Lemma bound. Write \(m=750-p\), where \(1\le p\le115\).
Whenever we use an initial maximum packing, we choose it in the full
good-incidence graph, breaking ties in a fixed order before assigning
any random marks. Its size \(T\) is therefore determined by the graph. Table~\ref{tab:cases}
specifies the choice of endpoint set, and hence the event \(F_C\), for
each possible deficit. At the smallest deficits, criticality gives the
required two-stars and we use both alternatives in the Hall condition.
As \(p\) increases we can use four-stars and then eight-stars. For
\(41\le p\le115\), removing endpoints with many neighbours in \(C\)
allows us to use a single bound based on represented colours alone.

\begin{table}[htbp]
\centering\small
\setlength{\tabcolsep}{4pt}
\begin{tabular}{@{}p{12mm}p{75mm}p{24mm}r@{}}
\toprule
\(p\)&Construction and size range&Endpoints&Ratio\\
\midrule
1&373 fixed disjoint two-stars
 &746&\(<0.164\)\\
2&Expanded pool for \(187\le T\le374\).
 375 fixed stars for \(T\ge375\)&\(m+T\), 750&\(<0.454\)\\
3&Free pairs in the expanded pool for \(187\le T\le737\).
 738 fixed stars for \(T\ge738\)&\(m+T\), 1476&\(<0.734\)\\
4--10&Expanded two-star pool for \(T_0(p)\le T\le m\)&\(m+T\)&\(<0.642\)\\
11&Expanded four-star pool for \(99\le T\le739\)&\(3m+T\)&\(<0.000886\)\\
12--40&All four-stars of the packing for \(T<250\).
 250 fixed stars for \(T\ge250\)&\(4T\), 1000&\(<0.744\)\\
41--115&Eight-stars filtered by attachment, using at most 250 stars
 &2000&\(<0.353\)\\
\bottomrule
\end{tabular}
\caption{The endpoint sets used for clique events. The third column
bounds the number of endpoints in each case. The last column bounds the ratio of the estimated failure probability to the allowance
in \eqref{eq:F-target}, rounded outwards. The allowance uses the entire
endpoint set. At \(p=2\) the two intervals are 187--229 and 230--374.}
\label{tab:cases}
\end{table}

\subsection{The smallest deficits}
For \(p=1,2,3\), we use criticality to obtain enough disjoint two-stars.
\begin{lemma}\label{lem:small-supply}
At deficit one a designated clique has at least 745 vertex-disjoint
good two-stars. At deficits two and three it has at least 187.
\end{lemma}
For \(p=1\), Lemma~\ref{app:lemma:2.1} gives at least \(747\) degree-\(D\)
centres in the \(749\)-clique \(C\). Each has two outside neighbours.
By \eqref{eq:maximum-clique-attachment} no outside vertex occurs at two
centres, so the outside pairs for different centres are disjoint. The
packing bounds at deficits two and three follow from the recolouring and matching
arguments in Appendix~\ref{app:small}.

For \(p=1\), fix 373 of the pairs before sampling. Write \(I\) for their centres
and \(E\) for their 746 endpoints. Each endpoint has only its own centre
as a neighbour in \(C\), so \(R_E(A)=|I\cap A|\) for every \(A\subseteq C\).
Given \(O=o\), the available set \(A_C\) is a uniform \(h\)-subset of
\(C\), where \(m=749\) and \(h=m-o\). Thus \(K_0=|A_C\cap I|\)
has law \(\Hypergeom(m,h,373)\). For \(0<t<1\),
\begin{equation}\label{eq:p1-exact-moment}
 \E[t^{K_0}\mid O=o]
 =\frac{[u^h](1+u)^{376}(1+ut)^{373}}{\binom mh}.
\end{equation}
Each selected centre in \(I\) contributes a factor \(ut\), and each
selected centre outside \(I\) contributes \(u\).

The factor \((1+u)^{376}\) from the centres outside \(I\) gives a
sharper bound on the coefficient. For a polynomial \(B(u)\) with
nonnegative coefficients and a positive argument \(v\), put
\[
 M_k(\theta)=\max_{0\le r\le k}\binom kr\theta^r(1-\theta)^{k-r},
 \qquad \theta=\frac v{1+v}.
\]
Direct multiplication gives
\[
 v^h[u^h](1+u)^kB(u)
 =\sum_j [u^j]B(u)\,v^j\binom{k}{h-j}v^{h-j}
 \le M_k(\theta)(1+v)^kB(v).
\]
Apply this inequality with \(k=376\), \(B(u)=(1+ut)^{373}\) and
\(v=h/o\). Write \(x=h/m\) and \(y=o/m\). At every supported occupancy
we have \(h,o>0\), and \((1+vt)/(1+v)=y+xt\). The numerical
envelope uses \(r(x,t)=(y+\sqrt{y^2+4xt})/2\). Since \(y+xt\le1\),
we have \((y+xt)^2-y(y+xt)=xt(y+xt)\le xt\), so \(y+xt\le r(x,t)\).
Consequently \eqref{eq:p1-exact-moment} gives
\begin{equation}\label{eq:p1-smoothed}
 \E[t^{K_0}\mid O=o]\le W_o^{\rm sm}:=
 \min\left\{1,M_{376}(x)
       \frac{r(x,t)^{373}}{\binom{749}{h}x^h y^o}\right\}.
\end{equation}
At \(t=e^{-1/3}\), \eqref{eq:p1-smoothed} bounds the exponential
term in the upper bound for the failure curve. We retain the constant
term and assign the failure bound one to every unsupported occupancy.
Appendix~\ref{app:p1-probability} gives the conservative colour profile,
the exact rational smoothing factor and the evaluation details. The
endpoint set has size 746 throughout.

For \(p=2\), choose a maximum packing of good two-stars. At an unused
centre with at least two good neighbours in the original endpoint set,
fix two of them. Every other unused centre has exactly one good
neighbour in the set and one outside. Indeed, it has at least two in
total, and two outside the set would extend the packing. At each such
centre, add the outside endpoint and fix this pair. The enlarged set
has size at most \(m+T\). An original endpoint is free if its owner was not sampled, and
an added endpoint is always free of this restriction. Among sampled
centres with two free endpoints, select a subfamily whose pairs are
disjoint. The selection argument in the appendix gives the coefficient
\[
 \delta_h=\frac{\binom{m-2}{h}-\binom{m-9}{h}}
                 {7\binom mh},\qquad 0\le\delta_h\le h/m.
\]
Thus \(\E R_E(A)\ge hT/m+(m-T)\delta_h\).
We use the same conditional failure curve on \(187\le T\le229\) and
\(230\le T\le374\), with a fixed exponential envelope on each interval.
For \(T\ge375\), keep 375 original stars and average the complete curve
over the hypergeometric distribution of the number of available centres.

\paragraph{The deficit-three mean.}\label{sec:p3-mean}
Let \(p=3\), \(m=747\), \(P=6\), and \(187\le T\le737\).
Start with a maximum two-star packing in the full good-incidence graph,
write \(I\) for its centre set, and add all good outside neighbours of
unused centres, giving \(|E'|\le m+T\).
At each unused centre at most one neighbour lies outside the original
endpoint set. Throughout this argument we use the good-incidence graph,
and \(P\) bounds each endpoint's full degree into \(C\). An original
endpoint is owned by its original centre. Added endpoints are unowned
in this calculation. The auxiliary owner labels used for averaging in
Section~\ref{sec:rank} are unnecessary here.

Fix \(0\le h\le m\), let \(A\) be a uniform \(h\)-subset of \(C\),
and expose \(U=A\cap I\). Put \(N=h-|U|\) and \(V_0=A\setminus I\).
Delete the pairs assigned to \(U\) from the endpoint set. At every unused
centre that still has two free neighbours, select a pair by a fixed
endpoint order. Let \(\mathcal L(U)\) be the set of these centres, and
join two of them in a conflict graph when their selected pairs intersect.
Each closed neighbourhood in the conflict graph has size at most ten.
An original endpoint meets at most \(P-2\) other unused centres, and an
unowned endpoint meets at most \(P-1\). A pair contains at most one
unowned endpoint.
Write \(M=m-T\ge10\). Conditional on \(U\), the set \(V_0\) is a uniform
\(N\)-subset of the unused centres. Generate it as the first \(N\) labels
of a uniformly random ordering of those centres. Select \(c\in\mathcal L(U)\) if it is first in its
closed conflict neighbourhood and belongs to \(V_0\). The random-order
calculation in Appendix~\ref{app:selection}, with no forbidden labels, gives
\[
 \E[R_{E'}(U\cup V_0)\mid U]\ge |U|+
 |\mathcal L(U)|\,\alpha(M,N),\qquad
 \alpha(M,N)=\frac{1-\binom{M-10}{N}/\binom MN}{10}.
\]

Before sampling, choose three good outside neighbours at each unused centre.
Its type \(\tau=(r_0,k_2,k_1)\) records the number of unowned endpoints
and the numbers of original owners occurring twice and once.
The four possibilities are
\[
 (0,0,3),\quad(0,1,1),\quad(1,0,2),\quad(1,1,0).
\]
For a uniform \(z\)-subset of \(T\) original owners, put
\(H_j=(z)_j/(T)_j\).
The probabilities that fewer than two chosen endpoints remain free are,
in this order,
\[
 3H_2-2H_3,\qquad H_1,\qquad H_2,\qquad H_1.
\]
If all three endpoints are owned, fewer than two remain free precisely
when all owners are present or the only absent owner occurs just once.
If one endpoint is unowned, failure requires every owner to be present.
Since \(H_2\le H_1\), the largest of the four failure probabilities is
\(\max\{H_1,3H_2-2H_3\}\). Thus at least two endpoints remain free with
probability at least \(e_3(T,z)\) for every type, where, for \(n\ge3\),
\[
 e_3(n,s)=1-\max\left\{\frac sn,\,
       \frac{s(s-1)(3n-2s-2)}{n(n-1)(n-2)}\right\}.
\]

Set \(Z=|U|=|A\cap I|\). Conditional on \(Z=z\), the sets \(U\) and
\(V_0\) are independent uniform samples from \(I\) and \(C\setminus I\),
of sizes \(z\) and \(h-z\), respectively.
For a centre of type \(\tau\), let \(e_\tau(T,z)\) be the conditional
probability that at least two of its three chosen neighbours remain free.
This bounds its probability of membership in \(\mathcal L(U)\) from below.
Both \(e_\tau(T,z)\) and \(\alpha(M,h-z)\) decrease with \(z\).
Their covariance is nonnegative, as follows by comparing two independent
copies of \(Z\). Moreover,
\[
 \E e_\tau(T,Z)=e_\tau(m,h),\qquad
 \E\alpha(M,h-Z)=\nu_h:=
 \frac{1-\binom{m-10}{h}/\binom mh}{10}.
\]
The first identity counts the fixed owner labels in \(A\) directly. The
second uses the probability that \(A\) meets a fixed set of ten unused
labels. Applying the covariance inequality to each centre and then taking
the minimum over the four types gives
\[
 \E R_{E'}(A)\ge hT/m+(m-T)\beta_h^{\rm free},\qquad
 \beta_h^{\rm free}=e_3(m,h)\nu_h,\qquad
 0\le\beta_h^{\rm free}\le h/m.
\]
We take the largest of three coefficients for the term \(m-T\). The
first is \(\beta_h^{\rm free}\), the second is \(\beta_h\) from
\eqref{eq:affine-mean} for the original endpoint set with
\(a_0=2,s=5,b=3,\ell_*=6\), and the third, obtained by selecting
disjoint pairs of original endpoints, is
\[
\frac{\binom{m-2}{h}-\binom{m-11}{h}}{9\binom mh}.
\]
For the third coefficient, a pair has at most two original owners, and
its closed neighbourhood in the conflict graph has size at most nine.
Appendix~\ref{app:selection} gives the random-order calculation. Each
coefficient gives a lower bound for the expected maximum packing in the
same enlarged set and is fixed throughout \(187\le T\le737\). The choice
of available stars uses only the exposed clique information. For
\(T\ge738\), retain 738 vertex-disjoint two-stars from the original packing.

\subsection{Two-stars and four-stars}
For \(4\le p\le10\), apply Lemma~\ref{lem:codegree-supply} with
\[
 U=p+1,\quad P=p+3,\quad\rho=3,\quad L=4p-7.
\]
Both inequalities in the lemma hold for each integer in this range.
The guaranteed number of initial stars is
\begin{equation}\label{eq:T0-low}
 T_0(p)=
 \left\lceil\frac{2p(p-2)(750-p)}{(4p-7)(p+3)}\right\rceil
 =190,215,234,248,260,269,276.
\end{equation}
Enlarge the original pool by all good neighbours of unused centres.
The full endpoint bound is \(m+T\).
Use \eqref{eq:affine-mean} with
\[
 d=2,\quad a_0=p,\quad s=p+3,\quad
 b=\max\{3,\lceil p/2\rceil\},\quad \ell_*=2b.
\]
Here \(p+1\le2b+1\) throughout the range, as required for \(\ell_*=2b\).

At \(p=4\) and \(T=190\), recovering stars at unused centres makes the
difference between satisfying the Local Lemma inequality and exceeding
its allowance. With recovery, the ratio of the upper bound to the
allowance is \(0.641384\) to six decimal places. If we set \(\beta_h=0\)
in \eqref{eq:affine-mean}, we count only the available original stars.
Keeping the endpoint bound \(m+T=936\), the colour curves and every
envelope parameter fixed, the ratio then lies between \(217716\) and
\(217717\).

The generic two-star bound \cite[Lemma~B.1]{Chen1000} guarantees
\(T_{\rm old}=\lceil(p-1)m/(2(p+3))\rceil\).
With recovery retained, the left-end ratio at \(p=4\) is \(8.086\)
at \(T_{\rm old}=160\) and \(0.821\) at \(T=187\), to three decimal
places. For \(5\le p\le10\), all six left-end ratios at
\(T=T_{\rm old}\) are below \(0.212\) with the present curves.
The right endpoints \(T=m\) are unchanged, so \eqref{eq:endpoints}
covers these larger intervals as well. For these six deficits Lemma~\ref{lem:codegree-supply}
improves the guaranteed supply even though the old supply suffices
for the new probability estimate.

For \(p=11\), take a maximum packing of good four-stars.
The elementary incidence count
\[
 (p-d+1)(m-T)\le dT(p+2)
\]
gives \(T\ge\lceil m(p-d+1)/((d+1)(p+1))\rceil=99\).
Enlarge its pool as before and use
\[
 d=4,\quad U=12,\quad P=s=14,\quad a_0=11,\quad b=5,\quad\ell_*=25.
\]
Its full endpoint bound is \(3m+T\), and
\eqref{eq:endpoints} covers the whole interval 99--739.

For \(12\le p\le40\), again take a maximum packing of good four-stars.
Let \(T_{\rm pair}(p)\) be the least positive integer \(t\) satisfying
\[
 48t^2+[(p-3)(p-4)-24]t-m(p-3)(p-4)\ge0,
\]
and put
\[
 T_0(p)=\max\left\{
 \left\lceil\frac{m(p-3)}{5(p+1)}\right\rceil,\,
 T_{\rm pair}(p)\right\}.
\]
The two bounds come from the incidence count above and the initial
pair count in \eqref{eq:rank-polynomial}. If \(T<250\), retain the
original pool, of size \(4T\), and take
\[
 a_0=p-3,\quad s=p+2,\quad
 b=\max\{5,\lceil(p+3)/3\rceil\},\quad
 \ell_*=p-6+4b.
\]
For the coefficient of \(m-T\) we may also use the priority bound in
Appendix~\ref{app:small}, taking the larger of the two coefficients.
That bound comes from choosing four neighbours in \(E_0\) at each
unused centre and selecting disjoint four-sets. The degree bound in
this calculation uses maximality when \(T<250\). When \(T\ge250\), we
keep 250 original stars and average over the full hypergeometric law of
their available centres. This uses 1000 endpoints.

Table~\ref{tab:profiles} gives the probability parameters for stars with
equal numbers of endpoints. The larger ceilings in the table apply only
to the colour probabilities. In the Hall condition we use the actual
deficit \(p\), in the occupancy law the actual order \(m=750-p\), and in
the mean estimates the degree bounds for the chosen construction.
\begin{table}[t]
\centering\small
\begin{tabular}{@{}ccrrrrccc@{}}
\toprule
Actual \(p\)&\(d\)&\(P\)&\(p_*\)&\(N\)&\(B\)&\(\alpha_R\)&\(\alpha_F\)&\(\kappa\)\\
\midrule
3--5&2&8&5&374&187&\(742/747\)&\(749/747\)&\(6/747\)\\
6--10&2&13&10&400&200&\(737/747\)&\(749/747\)&\(11/747\)\\
11--20&4&23&20&384&96&\(727/745\)&\(749/745\)&\(21/745\)\\
21--40&4&43&40&384&96&\(707/745\)&\(749/745\)&\(41/745\)\\
\bottomrule
\end{tabular}
\caption{Complete equal-star colour profiles on \(375\le O\le558\).
\(N\) is the maximum number of endpoint positions stored, and \(B\)
is the maximum number of stars in a prefix. Each probability law uses
the actual endpoint count \(n=dK\). The small deficits 1 and 2 use the
separate conservative data described in Appendix~\ref{app:small}.}
\label{tab:profiles}
\end{table}

\subsection{One filtered eight-star construction}
For \(41\le p\le115\), choose a maximum packing of good eight-stars and retain
its first \(\min(T,250)\) stars. Write \(E_0\) for their endpoints.
The event \(F_C\) uses the fixed set
\(E(C)=\{z\in E_0:|N(z)\cap C|<60\}\).
At most 27 endpoints are deleted. Indeed, if \(H>0\) endpoints have
attachment at least 60, let \(r_c\) count their neighbours at centre \(c\),
let \(B=\sum_c r_c\), and set \(u=B/H\ge60\). Pair counting and Cauchy's
inequality give
\[
 B^2\le m\sum_c r_c^2\le m[B+3H(H-1)],
 \qquad
 H\le\frac{m(u-3)}{u^2-3m}
 \le\frac{57m}{3600-3m}<28
\]
for \(635\le m\le709\). The middle expression decreases for \(u\ge60\);
the last increases with \(m\) and equals \(40413/1473<28\) at 709.

After the exposure \(\mathcal C\), choose \(K\) vertex-disjoint
eight-stars with centres in \(A_C\) and endpoints in \(E_0\), using
only the graph and \(\mathcal C\). For \(K\ge4\), at least
\(n_K=8K-27\) of their endpoints belong to \(E(C)\). Keep the first
\(n_K\) of these in fixed label order, each with its selected centre.
Each of the \(K\) groups now has between zero and eight endpoints. For
\(K\le3\) we use the failure bound one. We evaluate prefixes up to \(K=50\).

All remaining endpoints have attachment at most 59. Use the single profile
\[
 P=59,\quad p_*=115,\quad d_{\max}=8,\quad
 \alpha_R=691/741,\quad\kappa=116/741,\quad
 300\le O\le600.
\]
We count only represented colours. Their distribution uses the actual
number \(n_K\le373\) of endpoints, with the rank-dependent lower bounds
and the full error term for omitted partitions. Among groups with given
total size, the coefficient for no strong star is largest when the group
sizes differ by at most one. Moving an endpoint from a larger group to a
smaller one when their sizes differ by at least two can only increase the
relevant elementary symmetric coefficients. The union bound for at most
one weak star is largest when as many groups as possible have size eight.
Appendix~\ref{app:colour-laws} proves both estimates, including empty
groups and terms of degree above \(K\) in the weak-star bound. Together
they give a complete failure curve \(\psi_O\) for the prefixes.

If \(T<250\), maximality in the full good graph gives every unused centre
at least \(p-7\) neighbours in \(E_0\). Hence the available rank
\(R_{E_0}(A_C)\) is at least the least positive integer satisfying
\[
 192t^2+[(p-14)(p-15)-80]t-h(p-14)(p-15)\ge0,
\]
with value zero when \(h=0\). If \(T\ge250\), we instead use the exact
hypergeometric distribution of the number of available centres among the
250 retained stars. Both bounds use at most 2000 endpoints. For each
\(p\), the graph determines which case occurs, so the larger of the two
probability bounds, each averaged over the occupancies, covers both
cases in a single comparison.

\subsection{The verified inequalities}
The constructions above reduce the probability requirements to finite
calculations, whose exact specification is in
Appendix~\ref{app:verification}.
\begin{proposition}\label{prop:finite}
The complete bounds satisfy \eqref{eq:A-target} and
\eqref{eq:A-mass-target} for both root degrees and all permitted root
clique sizes, and satisfy \eqref{eq:F-target} for every construction in
Table~\ref{tab:cases} throughout its full size range.
\end{proposition}
For the \(A\)-events, the computation checks 772 inequalities using
289886 matching records. For each degree \(n=749,750\), one comparison
using total edge weight covers all \(c<250\). A separate comparison for
each \(250\le c\le634\) uses the cross edges. For the \(F\)-events,
14 comparisons check the two ends of the intervals for \(p=4,\ldots,10\),
and 75 comparisons for \(p=41,\ldots,115\) use 150 averages over the two
packing cases. The other constructions require 98 comparisons, giving
187 in all. Equation~\eqref{eq:endpoints} covers each whole interval.
For a fixed family of stars we average over the full hypergeometric
distribution.

We round survival lower bounds and Local Lemma allowances down, and
probability coefficients and omitted positive sums up. Final comparisons
use the resulting rational bounds. We check each proposed exponential
envelope at every prefix, so its validity is independent of how the
parameters were optimised. The appendix gives the inputs, finite ranges,
rounding rules and programs, and identifies the reference outputs.

Proposition~\ref{prop:finite} and the bound on \(W_i\) give all the
inequalities for the asymmetric Local Lemma in Section~\ref{sec:structure},
with the supports and activities specified there. With positive
probability all bad events are avoided. Section~\ref{sec:extension} then
gives a \(q\)-colouring, first of the uncoloured exceptional vertices,
then of the designated cliques and finally of the remaining vertices.
This contradicts the choice of \(G\) at degree \(750\). The reduction in
Section~\ref{sec:structure} gives Theorem~\ref{thm:main} at every larger
maximum degree as well.

\label{end-main}
\clearpage
\section*{AI assistance}
The author used GPT-6 Astra to assist with exploring candidate estimates,
numerical parameter trials, coding, checks of local arguments, and the
writing and editing of this manuscript. AI assistance does not constitute
independent verification of the proofs or computations. The author takes
responsibility for the mathematical arguments, computational results and
final text.

\clearpage
\appendix
\section{Structural details}\label{app:structure}

We first record a property of the largest possible cliques in the
critical graph \(G\). We then prove the decomposition used for the other
large cliques.

\begin{lemma}\label{app:lemma:2.1}
The \((D-1)\)-cliques of \(G\) are pairwise disjoint. Each has at most
two vertices of degree \(D-1\).
\end{lemma}
\begin{proof}
If two distinct \((D-1)\)-cliques intersect, a common vertex sees their
entire union except itself. Their intersection therefore has at least
\(D-3\) vertices. A vertex of one clique outside the other sees this
intersection, contrary to \eqref{eq:maximum-clique-attachment}.

For the degree assertion we use Reed's recolouring argument \cite{Reed}.
Suppose that \(v,w,s\) have degree \(q=D-1\) in \(G\) and belong to a
\(q\)-clique \(K\). Let \(z,u\) be the unique outside neighbours of
\(v,w\), respectively. They are distinct by \eqref{eq:maximum-clique-attachment}.
Add \(zu\) to \(G-K\) if it is missing. The maximum degree remains at most
\(D\), since both endpoints lost a clique neighbour. A new \(D\)-clique
would contain \(z,u\) and \(D-2\) other vertices. Those other vertices
together with \(z\) already form a \(q\)-clique in \(G\) to which \(u\)
has at least \(D-2\) neighbours, again contrary to
\eqref{eq:maximum-clique-attachment}. Thus \eqref{eq:smaller} colours the
augmented graph with \(q\) colours. Give \(w\) the colour of \(z\), which
differs from that of \(u\). Colour \(K\setminus\{v,w,s\}\) greedily,
leaving \(v,s\) uncoloured. We can now colour \(s\) because it has degree
\(q\) and its neighbour \(v\) is still uncoloured. Finally \(v\) has
degree \(q\) and sees the same colour on \(w,z\), so it too can be coloured.
\end{proof}

\subsection{Dense subgraphs and the clique decomposition}\label{dense-subgraphs-and-the-clique-decomposition}

Set \(h_0=444\), \(M=D-h_0=306\), \(s_0=2h_0-D-3=135\) and
\(k=\lfloor M/(s_0-1)\rfloor=2\). The table gives the numerical
inequalities used in the decomposition.

\begin{longtable}[]{@{}lr@{}}
\toprule\noalign{}
Condition & Value at \(D=750\) \\
\midrule\noalign{}
\endhead
\bottomrule\noalign{}
\endlastfoot
\(M\ge6\), \(s_0>1\), \(k\le2\) & \(306,135,2\) \\
\(12h_0-7(D+1)>15\) & \(71>15\) \\
\(2h_0-D-k-3\ge2\) & \(133\ge2\) \\
\(q-2M>1\) & \(137>1\) \\
\(2h_0-D-19\ge2\) & \(119\ge2\) \\
\(h_0-11\ge2\) & \(433\ge2\) \\
\(h_0-10>5\), \(2(h_0-11)>q\) & \(434>5,866>749\) \\
\(2h_0>D-1\) & \(888>749\) \\
\end{longtable}

Direct substitution also gives the coarser bounds
\(h_0\le D-6\), \(3h_0\le2D\), \(3M\ge D\),
\(59D/100+1\le h_0\le3D/5+2\), and
\(2D/5-2\le M\le41D/100-1\).

\textbf{Join corollary.} For \(t\ge6\), the first Cranston--Rabern theorem
makes \(K_t\vee I_3\) \(d_1\)-choosable: \(I_3\) is not almost complete
and neither exceptional value of \(t\) occurs. Thus lists of sizes at least
\(t-1\) on the independent part and \(t+1\) on the clique suffice.
In particular no such induced join occurs in \(G\).

\textbf{Cut lemma.} Let \(H\) be an induced subgraph on at most \(D+1\)
vertices with minimum degree at least \(h_0\), and let \(S\subseteq V(H)\)
have size \(s\ge s_0\). Let \(T\) consist of the vertices of
\(V(H)\setminus S\) with at least two neighbours in \(S\), and put
\(W=V(H)\setminus(S\cup T)\), \(|W|=w\).
Counting the cut between \(S\) and \(W\) gives
\begin{equation}\label{src:SG5}
s(w-M)\le e_H(S,W)\le w,
\qquad w\le M+k.
\end{equation}
Indeed every vertex of \(S\) has at least \(w-M\) neighbours in \(W\),
and every vertex of \(W\) has at most one neighbour in \(S\). Thus each
vertex of \(W\) has at least \(2h_0-D-k+1\) neighbours in \(S\cup T\).

\textbf{Dense-subgraph lemma.} Every induced \(H\) with
\(|V(H)|\le D+1\), \(\delta(H)\ge h_0\), is a clique or a maximum clique
plus one vertex. In the nonclique case its clique part has size at most
\(D-2\).

For a proof, suppose that \(H\) has six vertex-disjoint nonedges and
let \(Q_i\) be their common-neighbour sets in \(H\). If \(n=|V(H)|\),
then \(|Q_i|\ge2h_0-n\), and counting memberships gives
\[
\sum_{i<j}|Q_i\cap Q_j|\ge\sum_i|Q_i|-n
\ge12h_0-7n>15.
\]
Two of the sets \(Q_i\) therefore share two distinct vertices \(a,b\).
Both are adjacent to all four endpoints of two disjoint nonedges
\(xy,vw\). Colour \(G-V(H)\). Give \(x,y\) a common colour and \(v,w\)
a different common colour, which is possible because the two lists in
each pair share at least \(q-2M>1\) colours. Apply the cut lemma with
\[
S=(\{a,b\}\cup(N_H(a)\cap N_H(b)))\setminus\{x,y,v,w\}
\]
Here \(|S|\ge s_0\), and all four coloured endpoints belong to \(T\).
Colour \(W\), then the uncoloured vertices of \(T\), and then
\(S\setminus\{a,b\}\). A vertex of \(W\) has at least
\(2h_0-D-k-3\ge2\) uncoloured neighbours in \(S\cup T\), even after
the four endpoints have been coloured. Each vertex of \(T\) has two
neighbours in \(S\), and each vertex of \(S\setminus\{a,b\}\) is
adjacent to both \(a\) and \(b\). Finally, \(a\) and \(b\) each see two
distinct repeated colours, so they too can be coloured. At every step
at most \(q-1\) colours are forbidden, a contradiction.

The complement therefore has a vertex cover of size at most 10, namely
the endpoints of a maximum matching. A maximum clique \(C_H\) leaves at
most 10 vertices outside it and has size at least \(h_0-9\). If two outside
vertices \(u,v\) remain, each has at least \(h_0-9\) neighbours in \(C_H\).
Since \(|C_H|\le D-1\), their common neighbourhood in \(C_H\) has size at least
\[
2(h_0-9)-|C_H|\ge2h_0-D-17=121\ge4.
\]
Choose four such neighbours. Since \(|C_H|\ge h_0-9\ge7\), the
common-neighbour argument for large cliques gives a larger clique
in \(C_H\cup\{u,v\}\subseteq V(H)\). This contradicts the choice
of \(C_H\), so at most one vertex of \(H\) lies outside \(C_H\).
In the nonclique case, a clique part of size \(D-1\) would have an
exceptional vertex with at least \(h_0\) neighbours in it, contrary to
\eqref{eq:maximum-clique-attachment}. This proves the dense-subgraph lemma.

Let \(\mathcal K\) be all maximal cliques of order at least \(h_0+1=445\).
If two intersect, their union lies in the closed neighbourhood of a
common vertex, hence has order at most \(D+1\) and minimum degree at
least \(h_0\). The dense lemma makes this union a clique plus one vertex;
its only maximal cliques are its clique part and that exceptional vertex
with its neighbours in the part. For every intersecting pair their
intersection has size at least \(h_0\). No member can meet two others:
\(2h_0>D-1\) forces a triple intersection, and applying the dense lemma
to the triple union contradicts the presence of three maximal cliques.

Choose a maximum-size representative \(C_i\) from each intersection
component. Component unions \(S_i\) are disjoint and satisfy
\[
S_i=C_i\quad\hbox{or}\quad S_i=C_i\cup\{u_i\},\qquad
|N(u_i)\cap C_i|\ge h_0,
\]
\begin{equation}\label{src:SG6}
V(G)=L\,\dot\cup\,\bigdotcup_iS_i,
\qquad L=V(G)\setminus\bigcup_{K\in\mathcal K}K.
\end{equation}
Write \(|C_i|=D-p_i\). Then \(1\le p_i\le305\), so
\(p_i<D/2\), \(p_i+3<h_0\), and \(D-2p_i-3>0\).

For distinct \(x,y\notin C_i\), the nonempty sets
\(C_i\setminus N(x)\), \(C_i\setminus N(y)\) have distinct
representatives. Otherwise both equal \(\{a\}\), and the cliques
\((C_i\setminus\{a\})\cup\{x\}\) and
\((C_i\setminus\{a\})\cup\{y\}\) extend to members of \(\mathcal K\)
meeting \(C_i\). These extensions must coincide, since the component
of \(C_i\) has at most two members and neither extension is \(C_i\).
Their common extension contains \(C_i\setminus\{a\}\) and both \(x,y\),
contrary to the choice of \(C_i\).

Call an outside neighbour good when its full attachment into \(C_i\)
is at most \(p_i+3\). Each vertex of \(C_i\) has at least \(p_i\)
good neighbours. To prove this, suppose a centre \(v\) has two outside
neighbours \(x,y\), each with attachment at least \(p_i+4\). Choose
distinct representatives \(a\in C_i\setminus N(x)\) and
\(b\in C_i\setminus N(y)\), and colour \(G-(C_i\cup\{x,y\})\).
Give \(x,a\) one common colour and \(y,b\) a different common colour.
The lists in each pair share at least \(q-(D-p_i-4)-(p_i+1)=2\)
colours. Reserve \(v\) and another uncoloured neighbour \(w\) of \(x\) in the
clique. After colouring the other uncoloured clique vertices, colour
\(w\), which still has \(v\) uncoloured and sees the same colour on
\(x,a\). Finally \(v\) sees both repeated colours. If \(d(v)=q\), even
one such \(x\) is impossible. The same argument with one repeated pair
and one reserved clique neighbour colours \(v\) last, since its degree
is \(q\).
A degree-\(q\) centre has exactly \(p_i\) outside neighbours; a
degree-\(D\) centre has \(p_i+1\), of which at most one is bad.

\textbf{Adding edges outside the clique.} Add edges to \(G-C_i\) so
that the added edges have maximum degree \(L_0\) and the resulting graph
has maximum degree at most \(D\). If \(D-1-L_0\ge h_0\), the resulting
graph is \(q\)-colourable. Indeed, a new \(D\)-clique would induce in
\(G\) a subgraph on \(D\) vertices of minimum degree at least \(h_0\).
By the dense-subgraph lemma, it would either be a clique, contrary to
\(\omega(G)\le q\), or have a clique part of order \(D-1\), which the
same lemma excludes in the nonclique case. Now apply \eqref{eq:smaller}.
For the small-deficit arguments we use
\(L_0=2,3,4,5,6,7\), for which \(D-1-L_0=747,746,745,744,743,742\),
respectively.

\section{The neighbourhood probability estimate}\label{the-neighbourhood-probability-estimate}\label{app:vertices}

Here \(D=750\), \(q=749\), and the structural cutoff is 444. We derive
the bound on \(\PP(A_v)\) by counting pairs of neighbours that retain a
common colour.

For the root \(v\) of an event \(A_v\), choose a largest clique \(C(v)\)
containing \(v\). Put \(c=|C(v)|\), \(N=N_G(v)\),
\(n=d(v)\in\{749,750\}\) and \(Q=\overline{G[N]}\). Then
\begin{equation}\label{src:SG8}
c=1+\alpha(Q),\qquad c\le634.
\end{equation}
For \(v\in L\), \(c\le444\). If \(v\) lies in a non-designated clique
part \(C_i\), then \(C(v)\) is a member of \(\mathcal K\) meeting
\(C_i\). Since \(C_i\) is largest in its component, \(c=|C_i|\le634\).
The cliques chosen for different roots may overlap. There is one event
\(A_v\) for each root.

We use three restrictions on common neighbourhoods. If the endpoints
of a \(Q\)-edge \(xy\) have at least three common \(G\)-neighbours in
\(N\), these neighbours form a clique in \(G\). Otherwise, adjoining
\(v\) to them gives a connected incomplete graph on at least four
vertices completely joined to the independent pair \(x,y\), contrary
to the second Cranston--Rabern theorem in Section~\ref{sec:clique-parts}.
A triangle of \(Q\) has at most four common \(G\)-neighbours in \(N\).
Five would form a clique by the preceding argument and, together with
\(v\), would give the forbidden \(K_6\vee I_3\). Finally, two disjoint
\(Q\)-edges have at most two common \(G\)-neighbours in \(N\). Three
such neighbours and \(v\) would give a \(K_4\) joined to the four
endpoints. The complement of the graph on those endpoints contains a
matching of size two, so the graph is neither almost complete nor a
claw, and its order rules out \(I_3\). The first Cranston--Rabern
theorem excludes this join.

For \(c\ge7\), put \(T=C(v)\setminus\{v\}\), \(a_C=c-1\),
\(S=N\setminus C(v)\), and \(s=n-a_C\). Distinct vertices \(x,y\in S\)
have at most three common G-neighbours in \(C(v)\), hence at most two
in \(T\). If four existed in \(C(v)\), choose four including \(v\).
The common-neighbour consequence in Section~\ref{sec:clique-parts} gives a larger clique containing \(v\), contrary
to the choice of \(C(v)\).

Let \(L_C=e_G(S,T)\). Counting pairs of vertices of \(S\) through
\(T\), then applying Cauchy's inequality, gives
\[
L_C^2-a_CL_C-2a_Cs(s-1)\le0.
\]
Each vertex of \(S\) also misses a vertex of \(T\). Therefore
\begin{equation}\label{src:SG9}
e_Q(S,T)\ge E(n,c):=a_Cs-
\min\left\{s(a_C-1),
\left\lfloor\frac{a_C+\sqrt{a_C^2+8a_Cs(s-1)}}2\right\rfloor\right\}.
\end{equation}
The bound counts only edges of \(Q\) between \(S\) and \(T\) and permits
arbitrary additional edges in \(Q\). We verify the floor in the formula
by evaluating the quadratic inequality at consecutive integers. The
argument requires only the edge count.

\subsection{Survival kernel and structural mass}\label{survival-kernel-and-structural-mass}

In the notation of Section~\ref{sec:extension}, write \(c_v\) for the
number of coloured neighbours of \(v\) and \(r_v\) for the number of
their distinct colours. Colouring one more neighbour increases \(c_v\)
by one and \(r_v\) by at most one. Thus the saving
\({s_\varphi(v)}=c_v-r_v\) persists, and at most
\(d(v)-{s_\varphi(v)}\) colours can be forbidden when we come to colour
\(v\). The bad event at this root is
\begin{equation}\label{src:SG10}
A_v=\{v\hbox{ is uncoloured in }\varphi\hbox{ and }{s_\varphi(v)}\le d(v)-q\}.
\end{equation}
Avoiding \(A_v\) gives a still-uncoloured root at least one saving when
\(n=749\) and at least two when \(n=750\). The marks on \(B_2(v)\)
determine this event.

In each neighbourhood colour class of size at least two, take the first
two vertices in priority order. The pair is admissible if it is an edge
of \(Q\), and it contributes one saving if both vertices retain their
colour. Let \(Y_{\rm full}\) count such surviving pairs and let
\(Y_{\rm cross}\) count only those on the selected cross edges, fixed
from the graph. Then
\begin{equation}\label{src:SG11}
Y_{\rm cross}\le Y_{\rm full}\le{s_\varphi(v)}.
\end{equation}
Either count can therefore be used to bound the probability of the same
event \eqref{src:SG10}.

Write \(\delta_x=d_Q(x)\) and
\(d_x=|N_G(x)\setminus N_G[v]|\), and define \(\delta_y,d_y\) similarly.
Then \(d_x\le D-n+\delta_x\). Put
\(t_x=(D-n+\delta_x)/(D-2)\), and define \(t_y\) similarly. Define
\[
g(t)=\int_0^1e^{-tu}\,du,
\]
\[
s_3(u,v)=3\int_0^1\!\int_0^y(1-y)
(e^{-ux-vy}+e^{-vx-uy})\,dx\,dy.
\]
\begin{equation}\label{src:SG12}
b_{xy}=\min\{g(t_x)g(t_y)/\tau,s_3(t_x,t_y)\},
\quad\tau=9/10.
\end{equation}
Then \(0<b_{xy}\le1\). With the root temporarily removed as a blocker,
a selected admissible pair \(xy\in E(Q)\) from a class of size two has survival at least
\(\tau b_{xy}\), and from a larger class at least \(b_{xy}\). For a given
assignment of the two priorities to \(x,y\), multiplying the logarithmic
chord bounds for the exterior vertices gives the survival factor
\(e^{-(d_xU_1+d_yU_2)/(D-2)}\) when \(x\) has the earlier priority. A common exterior
neighbour only improves this bound, since
\(1-\max(U_1,U_2)/q\ge(1-U_1/q)(1-U_2/q)\). The first two order
statistics of a larger class are coordinatewise dominated by those of
three uniforms, with uniform orientation on the two endpoint labels.

Here all rates are at most \(749/748<3/2\). The common kernel satisfies
\begin{equation}\label{src:SG13}
b_{xy}\ge\frac{49}{50}-\frac35\frac{t_x+t_y}{2}
=A_n-\frac\gamma2(\delta_x+\delta_y),
\quad A_n=\frac{49}{50}-\frac{3(750-n)}{5\cdot748},
\quad\gamma=\frac3{5\cdot748}.
\end{equation}
We first prove the bound when the two rates are equal. By convexity,
the tangents to \(g(t)^2/\tau\) at \(1/2\) and \(3/5\) give lower
bounds. Using the alternating-series estimates
\(6065/10000\le e^{-1/2}\le6066/10000\) and
\(5488/10000\le e^{-3/5}\le5489/10000\), we obtain the lines
\(1-79t/125\) and \(483/500-567t/1000\). The first exceeds
\(49/50-3t/5\) on \([0,1/2]\) and the second on \([1/2,\infty)\).
For \(s_3(t,t)\), the \(j\)th moment of the sum of the two smallest
of three independent uniform variables is
\(6(2^{j+1}-1)/((j+1)(j+2)(j+3))\). The cubic Taylor lower bound for
the exponential therefore gives \(1-3t/4+7t^2/20-t^3/8\). Subtracting
\(49/50-3t/5\) gives a polynomial that is positive on \([0,1/2]\) by
the quadratic discriminant bound and increasing on \([1/2,3/2]\).
For unequal rates, Cauchy--Schwarz gives
\(g((u+v)/2)^2\le g(u)g(v)\), and averaging the two orientations gives
\(s_3(u,v)\ge s_3((u+v)/2,(u+v)/2)\). This proves \eqref{src:SG13}.

For \(c<250\), use the degree floors and kernel-mass estimate
\eqref{src:SG14}--\eqref{src:SG16} from the main text. The exact scalar
values are as follows.

\begin{longtable}[]{@{}lrrrrl@{}}
\toprule\noalign{}
Root degree \(n\) & \(L\) & \(H\) & \(B_1\) & \(B_2\) & Exact mass \(S_n\) \\
\midrule\noalign{}
\endhead
\bottomrule\noalign{}
\endlastfoot
749 & 165 & 251 & 166341 & 38260851 & \(1235978643/18700\) \\
750 & 166 & 251 & 166838 & 38405278 & \(72982177/1100\) \\
\end{longtable}

The estimate applies to the full \(Q\) at every permitted density.

When \(250\le c\le634\), select exactly \(E(n,c)\) cross Q-edges
using the fixed graph order. For the two sides use exterior rates
\[
t_C=(D-c+1)/748,\qquad t_S=749/748.
\]
A constant lower kernel weight is
\begin{equation}\label{src:SG17}
b\le\min\{g(t_C)g(t_S)/\tau,
S_3((t_C+t_S)/2)\},
\qquad S_3(t)=s_3(t,t).
\end{equation}
Taylor polynomials of odd degree 31 give positive lower bounds for the
integrals. The certificate rounds their minimum down to a rational with
denominator \(2^{128}\). Give the selected cross edges weight \(b\) and
the other edges of \(Q\) weight zero. Each side then has total weighted
degree \(S=bE(n,c)\). The term in \eqref{src:SG21} for a success lost
to the root still tests whether the pair is an edge of the full \(Q\).

\subsection{Positive root coefficients and weighted matching bounds}\label{positive-root-coefficients-and-weighted-matching-bounds}

We apply the exterior product comparison before averaging over the
numerical priorities. A pair-failure factor increases when exterior
blockers are added. Within each colour class we may condition on the
relative priority order, since the numerical order statistics remain
independent across classes.

Fix the transform parameter
\[
a_*=59657/62500,\qquad z_*=1-a_*=2843/62500.
\]
\begin{equation}\label{src:SG18}
t_0=t_1=0,\quad t_2=9/10,\quad t_j=1\ (j\ge3).
\end{equation}
Fix deterministic edge weights \(0\le b_e\le1\) bounded above by the
survival kernel. Choose an ordered sequence of disjoint unordered pairs
uniformly from all such sequences on the \(n\) labels. Define
\[
Z_r=\mathbb E\prod_{i=1}^r(1-b_{e_i}\mathbf1_{e_i\in Q}),\quad Z_0=1,
\]
\begin{equation}\label{src:SG19}
H_r=\mathbb E\left[\mathbf1_{e_{r+1}\in Q}
\prod_{i=1}^r(1-b_{e_i}\mathbf1_{e_i\in Q})\right]\le Z_r
\quad(2r+2\le n).
\end{equation}
The distinguished edge in \(H_r\) is tested for membership in the full
\(Q\), with weight one if it belongs to \(Q\). This explains the
distinction from \(Z_r-Z_{r+1}\), which would weight that edge by
\(b_e\). In particular, the test still includes edges outside the
selected cross edges.

We next restore the root to the experiment. For each neighbourhood colour
class \(\gamma\), let \(J_\gamma\) indicate that its selected admissible
pair succeeds with the root removed. It is zero for a class of size less
than two or an inadmissible pair. Put \(S=\sum_\gamma J_\gamma\),
let \(B\in\{0,1\}\) count successes destroyed by restoring the root,
and let \(I_v\) mean that the root is uncoloured. Put \(Y=S-B\).
The actual event \eqref{src:SG10} implies \(I_v=1\), \(Y\le n-q\), so
\begin{equation}\label{src:SG20}
\Pr(A_v)\le z_*^{-(n-q)}\mathbb E[z_*^{S-B}I_v].
\end{equation}
Conditional on the marks away from \(v\), the probability that the root
receives a given occupied colour and loses it is \((1-U_1)/q\), where
\(U_1\) is the minimum priority in that class. If its pair succeeds
without the root, the root can destroy that success and still be
uncoloured only by receiving that colour with priority between \(U_1\)
and the second priority \(U_2\). This condition also suffices. The loss
contribution is thus
\(J_\gamma(U_2-U_1)/q\). A failure
factor increases with pair priorities, whereas \(1-U_1\) decreases.
The product-space covariance inequality, obtained by comparing independent
copies and then conditioning coordinate by coordinate, bounds their joint
mean by the product of means. For a class of size \(k\),
\(\E(1-U_1)=k/(k+1)\), and for \(k\ge2\),
\(\E(U_2-U_1)=1/(k+1)\). Finally expand
\(z_*^{-B}=1+(a_*/z_*)B\) and use
\[
 z_*^S J_\gamma
 =z_*J_\gamma\prod_{\delta\ne\gamma}z_*^{J_\delta}
 \le z_*\ind_{\{\text{selected pair of }\gamma\in Q\}}
           \prod_{\delta\ne\gamma}z_*^{J_\delta}
\]
in the loss term, dropping its distinguished survival indicator. The
result retains the full coefficient \(a_*\), including size-two classes.
The exterior product comparison precedes these numerical integrations.

To express the resulting positive sum, set
\[
c_j(y)=1-a_*t_j+a_*t_jy,\quad
\mathcal K(z,y)=\sum_{j\ge0}c_j(y)z^j/j!,
\]
\[
\mathcal L_u(z,y)=\sum_{j\ge1}\frac{j}{j+1}c_j(y)z^j/j!,
\qquad
\mathcal L_g(z)=\sum_{j\ge2}\frac{z^j}{(j+1)j!}.
\]
Let
\(G_u=(n!/q^n)[z^n]q\mathcal L_u\mathcal K^{q-1}\) and
\(G_g=(n!/q^n)[z^n]q\mathcal L_g\mathcal K^{q-1}\).
Expand \(1-a_*t_jb_e=(1-a_*t_j)+a_*t_j(1-b_e)\) into its two
nonnegative terms. The variable \(y\) records each choice of the second
term, so the power of \(y\) counts the pairs to be included in the
matching product. Averaging over the uniform permutation of the labels
gives
\begin{equation}\label{src:SG21}
\Pr(A_v)\le\frac{z_*^{-(n-q)}}q
\left(\sum_r[y^r]G_u\,Z_r+a_*\sum_r[y^r]G_g\,H_r\right)
\le\sum_{r=0}^{\lfloor n/2\rfloor}\Gamma^A_{n,r}Z_r.
\end{equation}
The matchings in both sums use all \(n\) neighbourhood labels. At root
degree 749 the prefactor contributes no power of \(z_*\), and at degree
750 it contributes \(z_*^{-1}\).

To evaluate the root coefficients and bound the error from truncating the
class sizes, put \(N_c=q-1=748\),
\(B_j(y)=j!N_c^{-j}[z^j]\mathcal K(z,y)^{N_c}\), \(B_0=1\). For
\(1\le j\le749\),
\begin{equation}\label{src:SG22}
B_j(y)=\sum_{t=1}^j
\frac{(N_c+1)t-j}{j}\binom jtN_c^{-t}c_t(y)B_{j-t}(y).
\end{equation}
All multipliers are nonnegative; at \(j=749,t=1\) the multiplier is
exactly zero. The root polynomial is
\begin{equation}\label{src:SG23}
\sum_r\Gamma^A_{n,r}y^r=
\frac{(N_c/q)^n}{z_*^{n-q}}
\sum_{t=1}^n\binom ntN_c^{-t}
\left[\frac{t}{t+1}c_t(y)+\frac{a_*\mathbf1_{t\ge2}}{t+1}\right]
B_{n-t}(y).
\end{equation}
We truncate only the sums over class size, at \(T=64\). For each row
\(j>T\) we add the following remainder to the upper bound on every
coefficient.
\begin{equation}\label{src:SG24}
\epsilon_j=\frac{N_c+1}{N_c}
\frac{((j-1)/N_c)^T}{T!}\frac1{1-(j-1)/(N_c(T+1))}.
\end{equation}
Use zero for \(j\le T\). This follows because each exact coefficient
of \(B_j\) is at most \(B_j(1)=1\), and the omitted shift multipliers
are bounded by a geometric remainder of the exponential series. The
omitted outer root vector has total mass at most
\begin{equation}\label{src:SG25}
\epsilon_n^{\rm root}=
\frac{(N_c/q)^n}{z_*^{n-q}}
\frac{(n/N_c)^{T+1}}{(T+1)!}
\frac1{1-(n/N_c)/(T+2)}.
\end{equation}
The bracket in \eqref{src:SG23} is at most one at \(y=1\), and all
the denominators above are positive at the stated parameters. We add
\(\epsilon_j\) to individual coefficients during the recurrence and
use \(\epsilon_n^{\rm root}\) to bound the total omitted mass in the
outer sum. We retain every matching size, including the largest.

We next prove the matching bounds used with \eqref{src:SG21}. Write
\(\beta(0)=0\), \(\beta(j)=(j!)^{1/j}\) for \(j\ge1\), and interpolate
linearly between integers. The table of factorial-root bounds in the computational data certifies that this
function is increasing and concave on \([0,1502]\): its rational brackets
\(\ell_j,u_j\) satisfy
\begin{equation}\label{src:SG26}
\ell_j^j\le j!\le u_j^j,\qquad
\ell_j>u_{j-1},\qquad
2\ell_j>u_{j-1}+u_{j+1}
\end{equation}
for the relevant indices, with \(\ell_0=u_0=0\). The power bounds
enclose the factorial roots, and the adjacent-slope inequalities give
concavity of the piecewise-linear interpolation on this finite interval.

\textbf{A weighted extension of Brégman's bound.} For an \(n\)-by-\(n\) array with entries \(a_{ij}\) in \([0,1]\) and row totals \(r_i\), we will prove
\begin{equation}\label{src:SG27}
\sum_{\pi\text{ a bijection}}\prod_{i=1}^n a_{i,\pi(i)}
\le\prod_{i=1}^n\beta(r_i).
\end{equation}
The finite concavity of \(\beta\), proved using \eqref{src:SG26},
enters later, in \eqref{src:SG28}--\eqref{src:SG29}.

For a zero-one array, \eqref{src:SG27} is Brégman's theorem in the
form given in \cite[equation~(2.1)]{AF}. If a row is zero, the permanent
and its corresponding factor on the right are both zero. Otherwise all
row sums are positive and the theorem applies.

For a general array, fix the row sums. If two entries in a row lie
strictly between zero and one, increase one by \(t\) and decrease the
other by \(t\). Each product in \eqref{src:SG27} uses only one entry of
the row, so the permanent is affine in \(t\). Move in a direction in
which it does not decrease until an entry reaches zero or one. Apply this
operation row by row until every row has at most one fractional entry.

Write \(r_i=k_i+t_i\) with \(k_i=\lfloor r_i\rfloor\), \(0\le t_i<1\). The adjusted row has \(k_i\) ones and, if \(t_i>0\), one remaining entry \(t_i\). Independently round that entry to one with probability \(t_i\) and to zero otherwise. Averaging preserves each product on the left of \eqref{src:SG27}, hence its sum. Apply the binary result to every rounded array and average. The rounded row totals are independent, so the right side averages to
\[
\prod_i\bigl[(1-t_i)\beta(k_i)+t_i\beta(k_i+1)\bigr]
=\prod_i\beta(r_i).
\]
This proves \eqref{src:SG27} for the original array as well. Integer and zero row totals are included.

\textbf{From bijections to pairs.} Give each unordered pair of distinct
labels a nonnegative symmetric weight. Let \(W\) be the sum over perfect
pairings of the products of their edge weights. In the expansion of
\(W^2\), an ordered pair of perfect pairings has a union consisting of
doubled edges and alternating even cycles of length at least four. For
a fixed union, each long cycle has two ways to assign its edges
alternately to the two pairings. Equally, each long cycle has two
orientations as a cycle of a permutation. A doubled edge has a unique
assignment and gives a transposition. Thus for each union the numbers
of ordered pairs of pairings and permutations with those even cycles
agree, and they contribute the same weight. Therefore
\[
W^2=\sum_{\pi\text{ with only even cycles}}\prod_i a_{i,\pi(i)}
\le\sum_{\pi\text{ a bijection}}\prod_i a_{i,\pi(i)},
\]
where diagonal weights are zero. This direct count works for arbitrary nonnegative weights.

Suppose \(\sum_{e\in Q}b_e\ge S_0\), with \(0<S_0<\binom{n}{2}\), and set
\[
d_0=n-1-2S_0/n,\qquad d_r=n-2r,\qquad
\alpha_r=\min(1,2r/d_0),\qquad R_r=d_r+\alpha_rd_0.
\]
For \(r\ge1\) append \(d_r\) labelled dummy vertices. Give real pairs weight \(\alpha_r(1-b_e\mathbf 1_Q)\), real-dummy pairs weight \(1\) and dummy-dummy pairs weight \(0\). Every nonzero perfect pairing uses exactly \(r\) real pairs and matches the remaining real labels bijectively to the dummies, in \(d_r!\) ways. Since \(Z_r\) averages over all ordered disjoint real pairs, their weighted sum after completion is
\[
W=\frac{\alpha_r^r n!}{2^r r!}Z_r.
\]
Each dummy label has incident weight total \(n\). The real labels have average total at most \(R_r\), and each has total at most \(d_r+\alpha_r(n-1)\le2n-1<1502\). Apply the preceding pair count and \eqref{src:SG27}. Arithmetic--geometric mean and the verified monotonicity and concavity of \(\beta\) give
\begin{equation}\label{src:SG28}
Z_r\le\mathcal B_r^{\rm mass}:=
\min\left\{1,\frac{2^rr!}{\alpha_r^rn!}
\beta(n)^{d_r/2}\beta(R_r)^{n/2}\right\},\qquad\mathcal B_0=1.
\end{equation}
Take \(S_0=S_n\) from \eqref{src:SG16}. This step uses only the lower
bound on the sum of the edge weights.

For the selected cross kernel take \(S_0=S=bE(n,c)\). Its two sides have sizes \(a_C,s\) and weighted-degree totals \(S\) on each side. Their separate average incident totals after adding dummies are
\[
R_C=d_r+\alpha_r(n-1-S/a_C),\qquad
R_S=d_r+\alpha_r(n-1-S/s).
\]
Apply the same arithmetic--geometric mean and \(\beta\) concavity separately on these two groups to obtain
\begin{equation}\label{src:SG29}
Z_r\le\mathcal B_r^{\rm cross}:=
\min\left\{1,\frac{2^rr!}{\alpha_r^rn!}
\beta(n)^{d_r/2}\beta(R_C)^{a_C/2}\beta(R_S)^{s/2}\right\}.
\end{equation}
Both groups are nonempty for \(n=749,750\) and \(250\le c\le634\).
Every individual incident weight total lies in the interval on which
concavity was verified. The pairing argument applies for any arrangement
of the selected cross edges. The finite calculation uses these bounds
for \(r\ge1\) and the separate value at \(r=0\).

\textbf{Finite Certificate A.} For both root degrees, take one row with
\(c<250\), using \eqref{src:SG16}, \eqref{src:SG28}, and one row for each integer
\(250\le c\le634\), using \eqref{src:SG9}, \eqref{src:SG17}, \eqref{src:SG29}. Let
\(\overline{\Gamma^A}_{n,r}\) denote the coefficient upper bounds from
\eqref{src:SG22}--\eqref{src:SG24}, and let \(\overline{\mathcal B}_{n,c,r}\)
bound the applicable matching expression from above. Then
\(\PP(A_v)\) is at most
\begin{equation}\label{src:SG30}
U_{n,c}=\min\left\{1,
\sum_{r=0}^{\lfloor n/2\rfloor}
\overline{\Gamma^A}_{n,r}\overline{\mathcal B}_{n,c,r}
+\epsilon_n^{\rm root}\max_r\overline{\mathcal B}_{n,c,r}\right\}.
\end{equation}
The certificate has 772 rows and 289886 matching records. Its exact
positive coefficient arrays include the full root-loss term and all
remainder bounds for the class-size truncation at \(T=64\). The A-comparison file (entry \texttt{A comparisons} in Appendix~\ref{app:verification}) records
for every row the strict rational inequality
\begin{equation}\label{src:SG31}
U_{n,c}<x_0\exp[-\eta_0 V(c)-\eta_W(1126125/2)],
\end{equation}
using a certified lower bound for the exponential on the right. For
\(c<250\) the comparison uses \(V=Q_4\). For \(250\le c\le634\)
it uses \(V(c)\) from Section~\ref{sec:structure}, with a separate
comparison for each of the 385 clique sizes at each root degree.
The largest ratio, displayed to 16 decimal places, is
0.9152870180309839 at \(n=750,c=634\). The computational appendix
verifies the strict inequalities using rational bounds. The bound on
the root clique and the two edge-weight estimates ensure that these
comparisons cover every permitted neighbourhood density.

\section{The colour laws and their complete bounds}
\label{app:colour-laws}

We derive the distributions of the numbers of represented and fresh
marks used in \eqref{eq:complete-U}. Throughout, \(O\) counts the
initial clique colours, \(A=q-O\) counts the other colours, and \(n\)
counts the selected endpoint positions. The profiles considered here
have \(O,A\ge1\). We fix all class success bounds before evaluating a
prefix. For the Stirling numbers of the second kind \({s\brace r}\), we
use \({0\brace0}=1\) and give impossible indices value zero.

\subsection{The represented recurrence}
Conditional on \(s\) represented positions and \(r\) occupied represented
classes, order the classes uniformly. Let \(V_{s,r}(x)\) be the generating
polynomial for the number of successful class marks, with floor
\(\gamma_{\ell,j}\) at size \(\ell\) and injection rank \(j\).
For positive-mass states,
\begin{equation}\label{eq:represented-recurrence}
 V_{s,r}(x)=
 \sum_{\ell=1}^{s-r+1}
 \frac{\binom s\ell{s-\ell\brace r-1}}{r{s\brace r}}
 (1-\gamma_{\ell,r}+\gamma_{\ell,r}x)V_{s-\ell,r-1}(x).
\end{equation}
Set \(V_{0,0}=1\) and give the other impossible states value zero. To
obtain the size distribution in the recurrence, choose the positions of
the last class and partition the remaining positions among \(r-1\)
classes. The probabilities sum to one because
\[
 \sum_\ell\binom s\ell{s-\ell\brace r-1}=r{s\brace r}.
\]
Removing the last class leaves the first \(r-1\) ranks and their
success bounds unchanged. This step requires bounds that are independent
of the total number of endpoint positions and occupied classes at the
end of the recurrence.

Restrict the sum to class sizes \(\ell\le24\), keeping the original
Stirling denominator. Denote the resulting polynomial by
\(V^H_{s,r}\), and put
\[
 r^H_{s,b}=[x^b]\sum_{r=0}^{\min(s,O)}
                    \frac{(O)_r{s\brace r}}{O^s}V^H_{s,r}(x).
\]
The numbers \(r^H_{s,b}\) give the exact probabilities of \(b\) marks
and no represented class larger than 24, conditional on \(s\)
represented positions. Classes of size 5--24 receive no mark but still
count towards \(r\). In the experiment on \(n\) positions, the
probability of omitting a partition is bounded by
\[
 \tau_{n,O}=\min\{1,O\binom n{25}/q^{25}\},
\]
by choosing a monochromatic 25-set and one of the \(O\) represented labels.
No division by a zero Stirling number is performed.

\subsection{The fresh recurrence}
Put \(\phi_0=0\), \(a_0(y)=1\), and
\(a_\ell(y)=1-\phi_\ell+\phi_\ell y\) for \(\ell\ge1\).
The full fresh marked law on \(t\) positions is
\[
 F_t(y)=\frac{t!}{A^t}[z^t]
       \left(\sum_{\ell\ge0}a_\ell(y)\frac{z^\ell}{\ell!}\right)^A,
 \qquad f_{t,a}=[y^a]F_t(y).
\]
It satisfies \(F_t(1)=1\). Differentiating the power of the generating
series and comparing coefficients gives
\begin{equation}\label{eq:fresh-recurrence}
 F_j(y)=\sum_{\ell=1}^j
   \frac{(A+1)\ell-j}{j}\binom j\ell A^{-\ell}
   a_\ell(y)F_{j-\ell}(y).
\end{equation}
The coefficients of the recurrence are nonnegative for \(j\le A+1\).
At \(j=A+1\) the first coefficient vanishes, so this endpoint is included.

When the sum is cut at \(\ell=64\), an upper error for each coefficient
at \(j>64\) is
\[
 \epsilon_{A,j}=
 \frac{A+1}{A}\,
 \frac{((j-1)/A)^{64}}{64!}\,
 \frac1{1-(j-1)/(65A)}.
\]
For \(j\le64\) set the error to zero.
To obtain it, bound the omitted coefficient multipliers by
\((A+1)A^{-1}\binom{j-1}{\ell-1}A^{-(\ell-1)}\), use that every
coefficient of \(F_{j-\ell}\) is at most one, and bound the remaining
exponential-series tail geometrically. The denominator is positive on
the indicated domain. Thus the truncated computation, with these error bounds, encloses the
full law.

Let
\[
 w_{n,s}=\binom ns(O/q)^s(A/q)^{n-s}.
\]
\[
 G_n=\{s:0\le s\le n,\ n-s\le A+1\},\qquad
 \nu_{n,O}=\sum_{s\notin G_n}w_{n,s}.
\]
The supported joint and represented marginal masses are
\begin{align}
 \widetilde p_{b,a}
  &=\sum_{s\in G_n}w_{n,s}r^H_{s,b}f_{n-s,a},\label{eq:joint-law}\\
 \widetilde r_b&=\sum_{s\in G_n}w_{n,s}r^H_{s,b}.\label{eq:rep-marginal}
\end{align}
The two palettes are independent conditional on the set of positions
with represented colours. We use this conditional independence before
averaging over the size of that set. The marginal identity follows
from \(F_t(1)=1\) for the exact law. Rounding coefficients upwards may
make their sum exceed one.

\subsection{Grouping the marked positions}
For \(K\) equal \(d\)-stars, \(n=dK\), set
\[
 P_d(x,y,z)=
 \sum_{u=0}^1\sum_{v=0}^{d-u}
 \frac{d!}{u!v!(d-u-v)!}\,
 x^uy^vz^{(u+v-1)_+}.
\]
After averaging over the position set and the independent rank order,
the marked represented and fresh positions, conditional on their counts
\(b,a\), are uniform disjoint subsets of the \(n\) positions. In the
auxiliary experiment, let \(c_j\) be the probability mass on the supported
states in which no star is strong and the sum \(J\) equals \(j\). Then
\[
 c_j=\sum_{\substack{b,a\ge0\\b+a\le n}}\widetilde p_{b,a}
     \frac{[x^by^az^j]P_d(x,y,z)^K}
          {\binom n{b,a,n-b-a}}.
\]
All coefficients are nonnegative. The factorial argument in the main
text therefore gives, for \(p<L=(d-1)K\),
\[
 M_{K,p}=
 \min_{0\le k\le L-p}
 \frac{\sum_{j=0}^{L-k}c_j\binom{L-j}{k}}{\binom{L-p}{k}}.
\]
Equations~\eqref{eq:NS}--\eqref{eq:WU} give the bounds for no strong
star and for at most one weak star. We first set aside partitions with
an oversized class, so their probability enters only once when we add
\(\nu_{n,O}+\tau_{n,O}\) to complete the bound. Some such partitions
may also have unsupported counts. Counting this overlap twice still gives
an upper bound.

Now allow the \(K\ge1\) groups to have different sizes
\(a_1,\ldots,a_K\), with total \(n>0\). For \(0\le b\le n\), the
coefficient for placing \(b\) represented marks with at most one in each
group is \(e_b(a_1,\ldots,a_K)/\binom nb\), where \(e_b\) is the
elementary symmetric polynomial of degree \(b\). Suppose the group sizes
are integers in \([0,d]\). If two sizes \(a,b'\) satisfy \(a\ge b'+2\),
changing them to \(a-1,b'+1\) preserves the linear coefficient of
\((1+ax)(1+b'x)\) and increases the quadratic coefficient. The other
factors have nonnegative coefficients, so this change cannot decrease
any elementary symmetric coefficient. A maximum is thus attained when
the group sizes differ by at most one. Writing \(n=Ku+r\), \(0\le r<K\),
we obtain
\[
 NS_b^*(K,n)=
 \frac{[x^b](1+ux)^{K-r}(1+(u+1)x)^r}{\binom nb}.
\]
For the union bound over the star containing all the represented marks,
\(WU_0^*=K\). For \(b\ge1\), the coefficient is
\(\sum_i\binom{a_i}b/\binom nb\). The discrete differences of
\(\binom ab\) are nondecreasing, so transferring a position from a
smaller nonempty group to a larger group of size less than \(d\) cannot
decrease the sum. Filling as many groups as possible to size \(d\) gives,
for \(n=dv+w\) with \(0\le w<d\),
\[
 WU_b^*(K,n)=\frac{v\binom db+\binom wb}{\binom nb}.
\]
For the eight-star application, the nonzero terms lie in
\(0\le b\le\min(n,\max(K,8))\). In particular the weak coefficient at
\(K=4,n=5,b=5\) is one, although the no-strong coefficient is zero.
Zero-arity groups remain present in the coefficient \(WU_0^*=K\).
If \(n=0\), the family supplies no certificate and we use bound one.

\subsection{Positive formulas for the class floors}
We give the rational formulas behind \eqref{eq:Gamma}.
Let \(U_{(1)}=X<\cdots<U_{(\ell)}\) be the endpoint order statistics and
let \(V\) be the independent uniform variable. Put \(U_{(0)}=0\),
\(U_{(\ell+1)}=1\), and
\[
 B_u=\{U_{(u-1)}<V<U_{(u)}\},\qquad 1\le u\le\ell+1.
\]
On \(B_u\), the set \(Q_w\) contains the winner and all other endpoints
of its class with priority below \(V\), so \(D_\ell=\max(1,u-1)\).
For nonnegative integers \(s,t\), define
\[
 M_{\ell,u}(s,t)=\E[(1-X)^s(1-V)^t\ind_{B_u}].
\]
The exact positive formulas are
\begin{align}
 M_{\ell,1}(s,t)
 &=\frac{\ell\sum_{i=0}^t
        \binom ti i!(\ell-1+s+t-i)!}{(\ell+s+t+1)!},\label{eq:moment-u1}\\
 M_{\ell,u}(s,t)
 &=\frac{\ell!\sum_{i=0}^s
       \binom si (u-2+i)!(\ell+1-u+s-i+t)!}
 {(u-2)!(\ell+1-u)!(\ell+s+t+1)!},
 \quad 2\le u\le\ell+1.\label{eq:moment-u2}
\end{align}
For the first, integrate \(0<V<X<1\), expand
\(1-V=(1-X)+(X-V)\), and use the beta integral.
For the second, the joint density on \(0<X<V<1\) is
\[
 \frac{\ell!}{(u-2)!(\ell+1-u)!}
 (V-X)^{u-2}(1-V)^{\ell+1-u}.
\]
Expand \(1-X=(1-V)+(V-X)\) and apply the same integral.
As a check, summing over \(u\) gives
\(\ell/[(\ell+s)(t+1)]\).

For rational nonnegative rates \(\alpha,\kappa\), choose a positive
rational lower bound \(E_0\le e^{-\alpha-\kappa}\), and put
\[
 A^-_{\ell,u}
 =E_0\sum_{s=0}^{12}\sum_{t=0}^{12}
       \frac{\alpha^s\kappa^t}{s!\,t!}M_{\ell,u}(s,t).
\]
The identity
\(e^{-\alpha X-\kappa V}
 =e^{-\alpha-\kappa}e^{\alpha(1-X)+\kappa(1-V)}\)
has a positive expansion in the last two factors.
It follows that the dyadic number
\begin{equation}\label{eq:dyadic-floor}
 \gamma_{\ell,j}=2^{-40}
 \left\lfloor2^{40}
 \sum_{u=1}^{\ell+1}
 \left(1-\frac{P\max(1,u-1)}{O-j+1}\right)_+
 A^-_{\ell,u}\right\rfloor
\end{equation}
is a lower bound on the probability of a mark. For \(E_0\) we use the
Taylor polynomial of odd degree 31, which gives a positive lower bound
at these rates. We evaluate \eqref{eq:dyadic-floor} for
\(\ell=1,\ldots,4\) and use zero for larger classes.

The two-star profiles may also use the uncompensated rate and the simpler
product floor
\[
 J_\ell(\alpha)g(\kappa)
 \left(1-\frac{P(\ell+1)}{2(O-j+1)}\right)_+,\qquad
 g(\kappa)=\int_0^1e^{-\kappa v}\,dv,\quad
 \alpha=q/(q-d).
\]
Indeed, conditional on \(X=x,V=v\),
\(\E D_\ell\le1+(\ell-1)v\), and
\(\int_0^1ve^{-\kappa v}\,dv\le g(\kappa)/2\).
These estimates bound the weighted contribution of \(D_\ell\) in the
product formula. For the two integrals we use the rational lower bounds
\[
 J^-_\ell(\alpha)=\sum_{j=0}^{31}(-\alpha)^j
                        \frac{\ell!}{(\ell+j)!},\qquad
 g^-(\kappa)=\sum_{j=0}^{31}\frac{(-\kappa)^j}{(j+1)!}.
\]
These follow by integrating the odd Taylor lower bound for the exponential.
Both are positive for every specified profile. For fresh colours we
round the corresponding integral bound down to a dyadic rational. For
the equal two-star profiles in Table~\ref{tab:profiles}, we take the
largest of three lower bounds at each class size and rank. Two come
from the same moment formula with \(\alpha_R=(D-P)/(q-d)\) and
\(\alpha_R=q/(q-d)\), respectively. The third is the product bound
rounded down.

The four-star and filtered eight-star profiles use the moment formula
with \(\alpha_R=(D-P)/(q-d)\). For the conservative two-star profile at
deficits one and two we replace \(O-j+1\) by \(O-192+1\). This gives
a lower bound valid at every rank up to the storage limit of 192.
Here \(n\le192\le\min(O,A)+1\), and the joint computation below
includes every partition, so \(\nu_{n,O}=\tau_{n,O}=0\).
The input files specify the profiles and the downward rounding. The
bounds depend on \(O\), class size and rank. Since they are independent
of the final number of endpoints, the same bounds apply at every prefix, regardless of
how many more endpoints are used at later prefixes.

For the filtered eight-star profile, fresh colours are inactive.
Thus \(F_t(y)=1\) for all unrepresented counts and no fresh-domain
restriction is needed. Use the represented law with all binomial counts,
the actual \(n=8K-27\), and the complete partition remainder.
This gives the single high-deficit curve used in the main table.

\subsection{Joint and marginal laws at deficits one and two}
\label{app:constant-rank-laws}
At fixed occupancy \(O\), put \(A=q-O\). A prefix of \(K\le96\) two-stars
has \(n=2K\le192\) endpoint positions. Let \(x,y\) count represented and
fresh marks respectively, and let \(z\) count endpoint positions.
Write \(\rho_\ell,\phi_\ell\) for the stored success floors of a represented and a fresh class of
size \(\ell\), respectively. They are zero for \(\ell>4\).
Let \(r_0=s_0=1\) and, for \(\ell\ge1\), put
\[
 r_\ell(x)=\frac{1-\rho_\ell+\rho_\ell x}{\ell!},\qquad
 s_\ell(y)=\frac{1-\phi_\ell+\phi_\ell y}{\ell!}.
\]
Write \(R(z)=\sum_{\ell\ge0}r_\ell(x)z^\ell\) and
\(S(z)=\sum_{\ell\ge0}s_\ell(y)z^\ell\). Each labelled colour contributes
one class, possibly empty. Thus the probability generating polynomial for
the two mark counts on \(j\) endpoint positions is
\[
 P_j(x,y)=\frac{j!}{q^j}[z^j]R(z)^O S(z)^A .
\]
Here \(j!\) distributes the labelled positions among the classes after
the factorials in \(r_\ell,s_\ell\) have been included. Let
\(P_j^{\rm R}(x)=P_j(x,1)\) be the represented marginal.
Both initial polynomials at \(j=0\) equal one.

For \(1\le k\le j\), put \(c_{j,k}=(j-1)!/((j-k)!\,q^k)\).
Writing \(F=R^OS^A\), compare coefficients of \(z^{j-1}\) in
\(RSF'=(OR'S+ARS')F\). This gives the joint recurrence
\begin{equation}\label{eq:small-joint}
 P_j=\sum_{k=1}^j c_{j,k}P_{j-k}
       \sum_{u+v=k}\big((O+1)u+(A+1)v-j\big)r_u(x)s_v(y),
\end{equation}
where \(u,v\ge0\). For the marginal, \(S(z)|_{y=1}=e^z\).
Applying the same calculation to \(R^Oe^{Az}\), after multiplying its
derivative by \(R\), gives
\begin{equation}\label{eq:small-marginal}
 P_j^{\rm R}=\sum_{k=1}^j c_{j,k}P_{j-k}^{\rm R}
       \big(A r_{k-1}(x)+((O+1)k-j)r_k(x)\big).
\end{equation}

On the profile domain \(j\le n\le\min(O,A)+1\), all multipliers are
nonnegative. The small-deficit program (\texttt{source\_p2} in the
producer input map) retains the terms with shifts \(k\le64\) in both
sums. Write \(i=2\) for the joint recurrence and \(i=1\) for the
marginal. With the mark variables set to one, the bracket for shift
\(k\) is at most \((q+i)i^{k-1}/(k-1)!\). Each coefficient of an exact preceding
probability polynomial is at most one. Multiplying by \(c_{j,k}\)
and summing the resulting exponential-series tail gives the following
upper remainder per retained coefficient for \(j>64\):
\[
 \frac{q+i}{q}\,
 \frac{(i(j-1)/q)^{64}}{64!}\,
 \frac{1}{1-i(j-1)/(65q)} .
\]
The denominator is positive for \(j\le192\); the remainder is zero for
\(j\le64\). The program propagates earlier upper enclosures through
each positive recurrence before adding its corresponding remainder.
Represented mark degrees above 96 cannot feed a smaller degree and
are unused in the Hall and weak sums for \(K\le96\).
The joint coefficients of \(P_n\) give the Hall term, and those of
\(P_n^{\rm R}\) give the no-strong and weak bounds.

\section{The small deficits and auxiliary availability estimates}\label{the-small-deficits-and-auxiliary-availability-estimates}\label{app:small}

We first prove the two-star packing bounds for deficits two and three,
then estimate how many stars are available after sampling. Throughout,
a good two-star has its centre in the clique and two distinct good
outside neighbours as endpoints. Those endpoints may be adjacent.

\subsection{A two-star packing at deficit two}\label{quarter-supply-at-deficit-two}

Put \(D=750\) and \(q=D-1\). Let \(G\) be a finite simple graph of
maximum degree at most \(D\) that has no proper \(q\)-colouring, and
let \(C\) be a maximal clique of size \(k=D-2\). For \(z\notin C\),
write \(d_C(z)=|N_G(z)\cap C|\). Call \(z\) good if \(d_C(z)\le5\)
and bad otherwise. We assume the following properties.

\begin{enumerate}
\def\labelenumi{\arabic{enumi}.}
\tightlist
\item
  Every \(c\in C\) has degree \(D-1\) or \(D\), and at least two good outside neighbours.
\item
  Adding edges to \(G-C\) gives a properly \(q\)-colourable graph
  whenever the added edges have maximum degree at most five and the
  resulting graph has maximum degree at most \(D\).
\item
  Two distinct vertices outside \(C\) have at most three common neighbours in \(C\).
\end{enumerate}

\begin{lemma}\label{app:lemma:SQ}

Under these assumptions, there exist at least

\begin{equation}\label{src:Q0}
\boxed{\lfloor D/4\rfloor}
\end{equation}

vertex-disjoint good two-stars.

\end{lemma}

For a designated clique these hypotheses follow from
Proposition~\ref{prop:decomposition}(ii)--(iv) and \eqref{eq:critical-range}.
We first exclude certain configurations by showing that they would allow
a colouring of \(G\). We then use the resulting restrictions to assign
weights in the twin graph.

\paragraph{Using Hall's condition to exclude a configuration}\label{hall-reduction-and-a-control-gadget}

Each centre has two or three outside neighbours. After a proper colouring outside \(C\), let \(F_c\) be the set of their colours and let \(L_c=[q]\setminus F_c\). Then \(|F_c|\le3\) and \(|L_c|\ge k-2\). Thus Hall's criterion is satisfied if

\begin{equation}\label{src:QH}
\left|\bigcap_{c\in X}F_c\right|\le2
\quad\text{for every }|X|=k-1,
\qquad
\left|\bigcap_{c\in C}F_c\right|\le1.
\end{equation}

The minimum list size gives Hall's condition for every smaller subset.
An outside colouring satisfying \eqref{src:QH} would therefore extend
to \(G\), a contradiction.

\begingroup\samepage
\begin{lemma}\label{app:lemma:Q1}

There cannot be distinct outside vertices \(b,z_1,z_2\) and distinct centres \(c_0,c_1,c_2\) satisfying

\[
d_C(b)\ge5,\quad d_C(z_1)\ge2,\quad d_C(z_2)\ge3,
\]

\begin{equation}\label{src:Q1}
c_0b\notin E(G),\quad
c_1b,c_1z_1,c_2b,c_2z_2\in E(G),\quad
c_1z_2\notin E(G).
\end{equation}

\end{lemma}
\endgroup
\begin{proof}

In \(G-C\), complete \(b,z_1,z_2\) to a triangle. Add edges from \(b\) to every outside neighbour of \(c_0\). Add edges from \(z_2\) to every outside neighbour of \(c_1\) other than \(b,z_1\); there is at most one such neighbour. Only missing edges are added.

The two nonadjacencies in \eqref{src:Q1} rule out loops. The vertices
\(b,z_1,z_2\) gain at most five, two and three neighbours, respectively.
Each has lost at least that many neighbours in \(C\). Every other
vertex gains at most one edge for each neighbour it lost in
\(\{c_0,c_1\}\). If \(z_1\) or \(z_2\) was an outside neighbour of
\(c_0\), its edge to \(b\) was already counted in the triangle. The new
graph therefore has maximum degree at most \(D\), and the added edges
have maximum degree at most five. By assumption~2 it has a proper
\(q\)-colouring.

Write \(\beta,\alpha,\gamma\) for the distinct colours on \(b,z_1,z_2\). Then

\[
\beta\notin F_{c_0},\qquad
\{\beta,\alpha\}\subseteq F_{c_1},\quad\gamma\notin F_{c_1},
\qquad\{\beta,\gamma\}\subseteq F_{c_2}.
\]

The three forbidden sets are pairwise distinct and each has size at most
three, so any two intersect in at most two colours. Every \((k-1)\)-set
of centres contains at least two of \(c_0,c_1,c_2\), giving the first
condition in \eqref{src:QH}. Moreover, \(F_{c_1}\cap F_{c_2}\)
contains \(\beta\), which is absent from \(F_{c_0}\). Intersecting with
\(F_{c_0}\) leaves at most one colour, giving the second condition and
hence the contradiction.

\end{proof}

\paragraph{Low-degree centres and bad-neighbour groups}\label{low-degree-centres-and-bad-neighbour-groups}

Call a centre \emph{low} if its degree is \(D-1\) and \emph{high}
if its degree is \(D\). A low centre has exactly two outside neighbours,
both good.

\begin{lemma}\label{app:lemma:Q2}

If there are at least two low centres, every outside neighbour of every low centre has actual attachment at most two.

\end{lemma}\begin{proof}

Suppose a low centre \(c\) has an outside neighbour \(b\) with
\(d_C(b)\ge3\). Since \(C\) is maximal, choose \(c_0\in C\) nonadjacent
to \(b\). In \(G-C\), join \(b\) to the outside neighbours of \(c_0\).
There are at most three, and \(b\) lost at least three neighbours when
\(C\) was deleted. Each other endpoint of a new edge lost its neighbour
\(c_0\). The total maximum degree is therefore still at most \(D\),
and the new edges have maximum degree at most three. Assumption~2 gives
a proper colouring.

The colour of \(b\) belongs to \(F_c\), a set of size at most two, and does not belong to \(F_{c_0}\). Hence the intersection over all centres has size at most one. Every \((k-1)\)-subset retains at least one of the two low centres, so its forbidden intersection has size at most two. This is \eqref{src:QH}, a contradiction.

\end{proof}

A high centre without a bad neighbour has exactly three good neighbours. Let \(H\) be the set of all these centres. Every other high centre has exactly one bad neighbour and two good neighbours. Write \(N_C(z)=N_G(z)\cap C\) and put \(C_b=N_C(b)\) for each
bad outside vertex \(b\). The sets \(C_b\) are pairwise disjoint: a centre has at most three outside neighbours and at least two good ones.

An outside endpoint is \emph{private} if it has just one neighbour in
\(C\). A centre in \(C_b\) is private if both its good neighbours are
private endpoints. The two-star at such a centre shares no endpoint
with a star at any other centre.

\begin{lemma}\label{app:lemma:Q3}

Suppose that a bad vertex \(b\) and a good vertex \(z\) with \(d_C(z)\ge3\)
have a common centre. Every nonprivate centre of \(C_b\) is adjacent to
\(z\). There are at most three such nonprivate centres, and at least as
many private centres in \(C_b\).
\end{lemma}
\begin{proof}
Choose \(c_2\in C_b\cap N_C(z)\) and, by maximality of \(C\), a centre
\(c_0\) nonadjacent to \(b\). If a nonprivate centre \(c_1\in C_b\)
misses \(z\), it has a good neighbour \(z_1\) with \(d_C(z_1)\ge2\).
Lemma~\ref{app:lemma:Q1} applies with \(z_2=z\). Thus all nonprivate
centres in \(C_b\) are common neighbours of \(b,z\). Assumption~3
bounds their number by three, whereas \(|C_b|\ge6\).
\end{proof}

Call \(b\) type I if it satisfies the hypothesis of
Lemma~\ref{app:lemma:Q3}, and type II otherwise. Both good neighbours
of every centre in a type-II group have attachment at most two.

\paragraph{The remaining incidence graph}
In each type-I group take the two-star at every private centre.
These \(s\) stars are disjoint and their endpoints meet no remaining
centre. Discard the \(r\) nonprivate centres in these groups.
Lemma~\ref{app:lemma:Q3} gives \(r\le s\). If the clique has exactly one
low centre, discard it as well, and put \(e=1\) in that case and \(e=0\)
otherwise. Let \(L\) consist of the type-II centres and the remaining
low centres. Thus
\[
 k=s+r+e+|H|+|L|.
\]
Form the incidence graph \(B\) from all three good edges at each \(H\)-centre and both good
edges at each \(L\)-centre. The two endpoints at an \(L\)-centre form its
\(L\)-pair, counted once for that centre. Both have attachment at most
two by the type-II condition and Lemma~\ref{app:lemma:Q2}.

Replace each remaining centre \(c\) by adjacent twins, each joined
to its neighbours in \(B\). Give each external edge at either twin weight
\begin{equation}\label{src:QW}
 x_c=\begin{cases}1/12,&c\in H,\\1/6,&c\in L,\end{cases}
\end{equation}
and give the internal edge weight \(1-d_B(c)x_c\). The vertex constraints
require incident weight at most one. Each twin has incident weight one.
At an endpoint of attachment at most two the weight is at most \(2\cdot2/6\).
Every other endpoint meets only \(H\)-centres and has weight at most
\(2\cdot5/12\).

\paragraph{Complete odd-set verification}\label{complete-odd-set-verification}

The odd-set inequality for an odd vertex set \(W\) bounds its internal
edge weight by \((|W|-1)/2\). Suppose \(W\) contains both clones at
centres in \(A\), one at centres in \(A_1\), and endpoints \(T\). Put
\(a=|A|\), \(b=|A_1|\), \(t=|T|\) and \(d_T(c)=|N_B(c)\cap T|\).
The internal weight is \(a+Q\), where

\begin{equation}\label{src:QO}
Q=\sum_{c\in A}x_c(2d_T(c)-d_B(c))
 +\sum_{c\in A_1}x_cd_T(c).
\end{equation}

We must prove \(Q\le(b+t-1)/2\).

If \(b\ge1\), the inequality \(d_T(c)\le d_B(c)\) gives

\[
Q\le \sum_{c\in A\cup A_1}x_cd_T(c)
\le t/2\le(b+t-1)/2,
\]

where the bound \(t/2\) follows by summing half the incident weight at
each endpoint in \(T\). We may therefore assume \(b=0\), so \(t\) is
odd. When \(t=1\), every summand in \eqref{src:QO} is nonpositive.
For odd \(t\ge7\), the sharper bound on endpoint weights gives

\[
Q\le5t/12\le(t-1)/2.
\]

Only \(t=3,5\) remain. Call an endpoint low if it has at most two
neighbours in \(C\). Half the weight incident to a low endpoint is at
most \(1/3\), and the corresponding bound for any other endpoint is
\(5/12\).

\paragraph{\texorpdfstring{The case \(t=5\)}{The case t=5}}\label{the-case-t5}

If \(T\) contains a low endpoint, then \(Q\le(4\cdot5+4)/12=2\).
Otherwise only centres in \(H\) can contribute positively. Let \(u,v\)
count centres in \(A\cap H\) with three and two neighbours in \(T\),
respectively. We have \(12Q\le3u+v\) and \(3u+2v\le25\). If \(v\ge1\)
this gives \(3u+v\le24\), and if \(v=0\) it gives \(u\le8\).
Thus \(Q\le2\).

\paragraph{\texorpdfstring{The case \(t=3\) with a low endpoint}{The case t=3 with a low endpoint}}\label{the-case-t3-with-a-low-endpoint}

If exactly one endpoint is low, no \(L\)-pair is contained in \(T\). Nonpositive contributions can be ignored, and each remaining positive contribution is at most \(1/12\) per incidence. Consequently \(Q\le(5+5+2)/12=1\).

If exactly two endpoints are low, let \(\ell\) count the \(L\)-pairs
in \(T\) with centre in \(A\). They contribute \(\ell/3\) and use
\(2\ell\) of the four available incidences at the low endpoints.
Every other positive summand comes from an \(H\)-centre using at least
one remaining low incidence and contributes at most \(1/4\). Thus
\(Q\le\ell/3+(4-2\ell)/4\le1\).

If all three endpoints are low, their half incident weights give \(Q\le3\cdot(1/3)=1\).

\paragraph{\texorpdfstring{The case \(t=3\) with no low endpoint}{The case t=3 with no low endpoint}}\label{the-case-t3-with-no-low-endpoint}

Again only \(H\)-centres can contribute positively. With \(u,v\)
defined as for \(t=5\), these centres contribute \(3u+v\) incidences
with the three unordered pairs in \(T\). A centre adjacent to all three
endpoints contains three pairs, and a centre adjacent to two contains one.
Assumption~3 gives \(12Q\le3u+v\le9\), so \(Q\le3/4\).

Thus the mixed weight vector satisfies every odd-set inequality.

\paragraph{From the twin graph to a packing}\label{from-the-weighted-twin-graph-to-an-actual-packing}

The Tutte--Berge formula now turns the weight bound into a matching bound.
In any finite graph, suppose nonnegative edge weights satisfy all vertex
and odd-set constraints. Choose a set \(S\) attaining the Tutte--Berge
minimum. The edges meeting \(S\) have total weight at most \(|S|\).
In each even component of the graph after deleting \(S\), the internal
edge weight is at most half its order. In each odd component it is at
most half its order minus one half. The sum of these upper bounds is the
Tutte--Berge value, so the total edge weight is at most the matching number.

By the twin-graph correspondence in Section~\ref{sec:rank}, a bipartite
incidence graph with \(m\) centres has maximum two-star packing size equal
to the matching number of its twin graph minus \(m\).

For \(m=|H|+|L|\) remaining centres, the total weight is
\(m+|H|/4+|L|/3\). Adding the \(s\) already selected private stars,
we obtain a packing of size at least
\[
 s+\frac{|H|}{4}+\frac{|L|}{3}
 =\frac{k-e}{4}+\frac{3s-r}{4}+\frac{|L|}{12}
 \ge\frac{k-1}{4}.
\]
Since \(k=748\), integrality gives the required 187 stars.

\hfill\(\square\)

\subsection{A two-star packing at deficit three}\label{quarter-supply-at-deficit-three}

Set \(D=750\), \(p=3\), \(q=D-1\), \(k=D-p\) and \(r=p+3\). Let
\(G\) be a finite simple graph of maximum degree at most \(D\) with no
proper \(q\)-colouring, and let \(C\) be a maximal clique of order
\(k\). For \(z\notin C\), write \(d_C(z)=|N_G(z)\cap C|\). Call
\(z\) good if \(d_C(z)\le r\) and bad otherwise. We assume the
following two properties.

\begin{enumerate}
\def\labelenumi{\arabic{enumi}.}
\tightlist
\item
  Every centre \(c\in C\) has degree \(D-1\) or \(D\), and at least \(p\) good outside neighbours.
\item
  Adding edges to \(G-C\) gives a properly \(q\)-colourable graph
  whenever the added edges have maximum degree at most \(r=p+3\)
  and the resulting graph has maximum degree at most \(D\).
\end{enumerate}

Let \(h(C)\) be the maximum number of vertex-disjoint good two-stars
centred in \(C\). Each consists of a centre and two distinct good outside
neighbours, which may be adjacent.

\begin{lemma}\label{app:lemma:SP}

Under these hypotheses, \(C\) has at least

\begin{equation}\label{src:P0}
\boxed{\left\lceil\frac{D-p}{4}-\frac{3}{4(p+3)}\right\rceil}
\end{equation}

vertex-disjoint good two-stars. For \(p=3\), this is \(h(C)\ge\lfloor D/4\rfloor=187\).

\end{lemma}

The proof uses the maximality of \(C\), the absence of a \(q\)-colouring
of \(G\), and the colouring property after adding edges to \(G-C\).
Together these restrict the incidence graph, in which every good
endpoint has at most \(p+3\) neighbours in \(C\). For a designated
clique the structural decomposition permits added edges of maximum
degree six, since \(749-6\ge444\).

\paragraph{Hall's condition at deficit three}\label{the-two-hall-conditions-and-the-general-control-gadget}

A low centre, of degree \(D-1\), has exactly \(p\) outside neighbours, all good. A high centre has \(p+1\) outside neighbours and at most one bad neighbour.

For any proper colouring outside \(C\), let \(F_c\) be the set of colours on its outside neighbours and \(L_c=[q]\setminus F_c\). Then \(|F_c|\le p+1\) and \(|L_c|\ge k-2\). Hall's theorem extends the colouring to the clique if

\begin{equation}\label{src:PH}
\left|\bigcap_{c\in X}F_c\right|\le p
\quad (|X|=k-1),
\qquad
\left|\bigcap_{c\in C}F_c\right|\le p-1.
\end{equation}

The minimum list size already gives Hall's condition for smaller sets.
We use the two remaining conditions to exclude the following configuration.

\begin{lemma}\label{app:lemma:P2}

There cannot be distinct outside vertices \(b,z_1,z_2\) and distinct centres \(c_0,c_1,c_2\) such that

\[
d_C(b)\ge p+3,\quad d_C(z_1)\ge2,\quad d_C(z_2)\ge p+1,
\]

\begin{equation}\label{src:PG}
c_0b\notin E(G),\quad
c_1b,c_1z_1,c_2b,c_2z_2\in E(G),\quad
c_1z_2\notin E(G).
\end{equation}

\end{lemma}\begin{proof}

In \(G-C\), complete \(b,z_1,z_2\) to a triangle. Join \(b\) to every outside neighbour of \(c_0\), at most \(p+1\). Join \(z_2\) to every outside neighbour of \(c_1\) other than \(b,z_1\), at most \(p-1\). Only missing edges are added.

Since \(b\notin N(c_0)\) and \(z_2\notin N(c_1)\), these operations
create no loops. The vertices \(b,z_1,z_2\) gain at most \(p+3,2,p+1\)
neighbours, respectively. Each lost at least that many neighbours in
\(C\). Any other vertex gains at most one edge for each neighbour it
lost in \(\{c_0,c_1\}\). If a distinguished vertex is also one of the
specified outside neighbours, the edge to be added is already counted
in the triangle. Thus the resulting graph has maximum degree at most
\(D\) and the added edges have maximum degree at most \(p+3\).
Assumption~2 gives a proper colouring.

Let the colours on \(b,z_1,z_2\) be \(\beta,\alpha,\gamma\), all distinct. Then

\[
\beta\notin F_{c_0},\qquad
\{\beta,\alpha\}\subseteq F_{c_1},\quad\gamma\notin F_{c_1},
\qquad\{\beta,\gamma\}\subseteq F_{c_2}.
\]

The three forbidden sets are pairwise distinct. Each has size at most
\(p+1\), so any two share at most \(p\) colours. Every \((k-1)\)-subset
of centres contains at least two of \(c_0,c_1,c_2\), giving the first
inequality in \eqref{src:PH}. Moreover, \(F_{c_1}\cap F_{c_2}\)
contains \(\beta\), which is absent from \(F_{c_0}\). Intersecting with
\(F_{c_0}\) leaves at most \(p-1\) colours. Hall's theorem now gives a
\(q\)-colouring of \(G\), a contradiction.

\end{proof}

\paragraph{Low centres and stars with private endpoints}\label{low-centres-and-private-star-credit}

\begin{lemma}\label{app:lemma:P3}

If there are at least two low centres, every outside neighbour of every low centre has actual attachment at most \(p\).

\end{lemma}\begin{proof}

Suppose a low centre \(c\) has a neighbour \(b\) with \(d_C(b)\ge p+1\).
Since \(C\) is maximal, some \(c_0\in C\) is nonadjacent to \(b\). In
\(G-C\), join \(b\) to the outside neighbours of \(c_0\), of which there
are at most \(p+1\). The vertex \(b\) lost at least \(p+1\) neighbours
in \(C\), and each other endpoint of a new edge lost \(c_0\). The maximum
degree is still at most \(D\). Since the added edges have maximum degree
at most \(p+1\), assumption~2 gives a \(q\)-colouring.

The colour of \(b\) is in \(F_c\), a set of size at most \(p\), and is absent from \(F_{c_0}\). Thus the full forbidden intersection has size at most \(p-1\). Every \((k-1)\)-subset still contains a low centre, so its forbidden intersection has size at most \(p\). This is \eqref{src:PH}, a contradiction.

\end{proof}

For a bad outside vertex \(b\), let \(C_b=N_C(b)=N_G(b)\cap C\). These groups are pairwise disjoint, since every centre has at most one bad outside neighbour. A centre in \(C_b\) has exactly \(p\) good neighbours. Call the centre private if all these good neighbours have actual attachment one; write \(P_b\) for those centres and \(B_b=C_b\setminus P_b\).

\begin{lemma}\label{app:lemma:P4}

Suppose some centre in \(C_b\) has a good neighbour \(z\) with \(d_C(z)\ge p+1\). Then every nonprivate centre of \(C_b\) is adjacent to \(z\), and

\begin{equation}\label{src:PC}
|B_b|\le p+3,\qquad |P_b|\ge1,
\qquad |B_b|\le(p+3)|P_b|.
\end{equation}

\end{lemma}\begin{proof}

Choose \(c_2\in C_b\cap N_C(z)\) and a centre \(c_0\) nonadjacent to \(b\), using maximality. If a nonprivate \(c_1\in C_b\) were nonadjacent to \(z\), it would have a good neighbour \(z_1\) with attachment at least two. This neighbour is different from \(z\), and \(c_1\ne c_2\). Lemma~\ref{app:lemma:P2} applies with \(z_2=z\), since the bad vertex has attachment at least \(p+4\). Thus all nonprivate centres lie in \(N_C(z)\), giving their bound \(p+3\). The group has at least \(p+4\) centres, so at least one is private, and \eqref{src:PC} follows.

\end{proof}

Call a bad group type I when Lemma~\ref{app:lemma:P4} applies, and type II otherwise. For a type-II group every good endpoint at every centre has actual attachment at most \(p\).

At every private centre of every type-I group, select any two of its \(p\) private good endpoints. All these pairs are disjoint from each other and from the endpoints of every remaining centre. Let their number be \(s\), and delete these centres and selected private pairs from the packing instance. Retain every nonprivate centre of type-I groups, and let their number be \(b\). Summing \eqref{src:PC} gives

\begin{equation}\label{src:PB}
b\le(p+3)s.
\end{equation}

\paragraph{Edge weights in the twin graph}\label{a-mixed-clone-vector}

All centres outside the selected private stars remain. We divide them
into three classes.

\begin{itemize}
\tightlist
\item
  \(H\): high centres with no bad neighbour; retain all \(p+1\) good incidences.
\item
  \(L\): type-II bad-neighbour centres, together with all low centres when there are at least two low centres; retain their \(p\) good incidences. All these endpoints have actual attachment at most \(p\).
\item
  \(R\): nonprivate centres of type-I groups, together with the unique low centre when exactly one exists; retain their \(p\) good incidences.
\end{itemize}

Let \(e\in\{0,1\}\) record the unique-low-centre exception, so \(|R|=b+e\). These classes exhaust the centres not already used in the private packing, and

\begin{equation}\label{src:PD}
k=s+|H|+|L|+b+e.
\end{equation}

Make the twin graph by replacing each remaining centre with two adjacent clones, each incident to all its retained good neighbours. Give each external edge at a centre \(c\) weight

\begin{equation}\label{src:PW}
x_c=
\begin{cases}
1/[4(p+1)],&c\in H,\\
1/(4p),&c\in L,\\
1/[4(p+3)],&c\in R.
\end{cases}
\end{equation}

Give the internal edge at \(c\) weight \(1-d(c)x_c\), where \(d(c)\)
is the number of its retained good neighbours. This number is \(p+1\)
on \(H\) and \(p\) on \(L\cup R\). The weights are positive, and the
weight incident to each twin sums to one.

Define half the incident weight at an endpoint \(z\) by

\[
v_z=\sum_{c\sim z}x_c.
\]

All actual attachments refer to the original clique. If \(d_C(z)\le p\), then \(v_z\le p/(4p)=1/4\). Otherwise no \(L\)-centre is adjacent to \(z\), and

\begin{equation}\label{src:PV}
v_z\le\frac{p+3}{4(p+1)}\le\frac38.
\end{equation}

Thus the actual weighted endpoint degree, \(2v_z\), is at most one.

\paragraph{Odd-set verification}\label{odd-set-verification}

For an odd set of twin-graph vertices, let \(A\) be its centres with both clones present, \(A_1\) those with one clone, and \(T\) its endpoint set. Put \(a=|A|\), \(a_1=|A_1|\), and \(t=|T|\). Its internal edge weight is \(a+Q\), where

\begin{equation}\label{src:PO}
Q=\sum_{c\in A}x_c(2d_T(c)-d(c))
 +\sum_{c\in A_1}x_cd_T(c).
\end{equation}

We need \(Q\le(a_1+t-1)/2\). When \(a_1\ge1\),

\[
Q\le\sum_{c\in A\cup A_1}x_cd_T(c)
\le\sum_{z\in T}v_z\le t/2\le(a_1+t-1)/2.
\]

It remains to consider \(a_1=0\), so \(t\) is odd. At \(t=1\), every
summand is nonpositive because every remaining centre has at least
three retained neighbours. For \(t\ge5\), \eqref{src:PV} gives

\[
Q\le3t/8\le(t-1)/2.
\]

It remains to check \(t=3\). Negative and zero summands may be discarded in an upper bound.

\paragraph{Deficit three, three endpoints}\label{deficit-three-three-endpoints}

At deficit three the edge weights for centres in \(H,L,R\) are
\(1/16,1/12,1/24\), respectively, and their degrees in the incidence
graph are four, three and three. If the three endpoints in \(T\) each
have attachment at most three, then \(v_z\le1/4\) at each endpoint,
giving \(Q\le3/4\).

Otherwise at least one endpoint in \(T\) has attachment greater than
three, so no \(L\)-centre can meet all of \(T\). The table lists all
positive contributions to \(Q\) in this case.

\begin{longtable}[]{@{}
  >{\raggedright\arraybackslash}p{(\linewidth - 6\tabcolsep) * \real{0.2000}}
  >{\raggedleft\arraybackslash}p{(\linewidth - 6\tabcolsep) * \real{0.2667}}
  >{\raggedleft\arraybackslash}p{(\linewidth - 6\tabcolsep) * \real{0.2667}}
  >{\raggedleft\arraybackslash}p{(\linewidth - 6\tabcolsep) * \real{0.2667}}@{}}
\toprule\noalign{}
\begin{minipage}[b]{\linewidth}\raggedright
Centre type
\end{minipage} & \begin{minipage}[b]{\linewidth}\raggedleft
Incidences into \(T\)
\end{minipage} & \begin{minipage}[b]{\linewidth}\raggedleft
Contribution to \(Q\)
\end{minipage} & \begin{minipage}[b]{\linewidth}\raggedleft
Contribution per incidence
\end{minipage} \\
\midrule\noalign{}
\endhead
\bottomrule\noalign{}
\endlastfoot
\(H\) & 3 & \(1/8\) & \(1/24\) \\
\(L\) & 2 & \(1/12\) & \(1/24\) \\
\(R\) & 2 & \(1/24\) & \(1/48\) \\
\(R\) & 3 & \(1/8\) & \(1/24\) \\
\end{longtable}

Every positive contribution is at most \(1/24\) per incidence. The three good endpoints have total actual attachment at most \(3\cdot6=18\). Consequently \(Q\le18/24=3/4<1\).

Thus the mixed weight vector satisfies every odd-set inequality.

\paragraph{Adding the private stars to the packing}\label{actual-packing-and-the-private-centre-budget}

The weighted matching argument above and the twin-graph correspondence
in Section~\ref{sec:rank} give a two-star packing of size at least the
total weight minus the number of centres.

The total weight of \eqref{src:PW}, with \(m=|H|+|L|+b+e\), is

\[
m+\frac{|H|+|L|}{4}+\frac{p(b+e)}{4(p+3)}.
\]

Subtracting \(m\) gives a lower bound on the number of disjoint
two-stars in the remaining incidence graph. We can add the \(s\)
private stars because their endpoints meet no remaining centre. By
\eqref{src:PB} and \eqref{src:PD},

\[
\begin{aligned}
h(C)
&\ge s+\frac{|H|+|L|}{4}+\frac{p(b+e)}{4(p+3)}\\
&=\frac{k}{4}+\frac{3s}{4}-\frac{3b}{4(p+3)}-
   \frac{3e}{4(p+3)}\\
&\ge\frac{k}{4}-\frac{3}{4(p+3)}.
\end{aligned}
\]

Integrality proves \eqref{src:P0}. For \(p=3\) the last expression is \(D/4-7/8\), whose ceiling is \(\lfloor D/4\rfloor\). This identity follows by writing \(D=4a+s\), \(0\le s\le3\). This completes Lemma~\ref{app:lemma:SP}.

\hfill\(\square\)

\subsection{The deficit-one probability bound}
\label{app:p1-probability}
Here \(m=q=749\) and the 373 fixed two-stars have 746 endpoints.
Use the conservative two-star colour data with \(P=5\),
\(\alpha=q/(q-2)\), \(\kappa=3/(q-2)\), and storage limit 192.
At deficit one, every endpoint has attachment one and the representative
exterior degree is at most two, so these floors remain valid.
They do not depend on the clique order, except through the actual
occupancy \(o\). The window \(375\le o\le558\) satisfies
\(o>192\) and \(192\le\min(o,q-o)+1\).

For the two-star experiment with rank-independent floors in
Appendix~\ref{app:colour-laws}, define
\[
 Q_{o,K}=\E_{\rm aux}
   \left[\ind_{\{\forall i:r_i\le1\}}(1/16)^J\right].
\]
Equivalently, \(Q_{o,K}=\sum_j c_j16^{-j}\), where the grouping
coefficients for \(P_2\) use the full joint distribution of the two
mark counts. At deficit one the exponential-product bound is
\(16Q_{o,K}\). We obtain it from the deficit-two data. The field
\texttt{hall\_upper\_hex} in the stored prefix file gives an upper
bound \(H^{(2)}_{o,K}\) for \(256Q_{o,K}\), including every positive
recurrence remainder. Because this bound is stored before taking the
minimum with one, we may divide it by 16 to bound the deficit-one
term. This division must precede the cap at one.

For \(K\ge2\), form two complete bounds
\begin{align*}
 U^{\rm H}_{o,K}&=\min\{1,H^{(2)}_{o,K}/16+WU_{o,K}\},\\
 U^{\rm NS}_{o,K}&=\min\{1,
       \min(H^{(2)}_{o,K}/16,NS_{o,K})+WU_{o,K}\}.
\end{align*}
Here \(NS\) and \(WU\) are upper bounds from
\eqref{eq:NS}--\eqref{eq:WU} for the full law. They include the recurrence
remainders and every term of the weak-star union bound. For \(K=0,1\)
use one. Take the running minimum for each curve and extend it constantly
after \(B=96\). The archive calls these curves \texttt{hall\_only} and
\texttt{hall\_NS}. Apply the smoothed rank-moment bound with a separate
fixed exponential envelope for each curve. At each occupancy, take the
smaller of the two resulting conditional bounds and then average over the
occupancy distribution. The \texttt{p1} entry of \texttt{PROFILE\_MAP.json}
identifies the prefix, curve and smoothed occupancy files used here.

The exact moment identity and smoothing inequality are
\eqref{eq:p1-exact-moment}--\eqref{eq:p1-smoothed} in the main text.

\subsection{Compatible selection in a uniform set}
\label{app:selection}
Let \(C\) be an \(m\)-set and \(J\) a graph on some of its labels.
For each \(c\in V(J)\), let \(D_c\) be its closed neighbourhood, of size
\(\ell_c\), and let \(\Omega_c\subseteq C\setminus D_c\) be a forbidden
set of size \(k_c\). Given a uniform \(h\)-subset \(A\), call \(c\)
eligible if \(c\in A\) and \(A\cap\Omega_c=\varnothing\). The expected
size of a largest independent set of eligible vertices is at least
\begin{equation}\label{eq:selection}
 \sum_{c\in V(J)}
 \frac{A_{k_c}(h)-A_{k_c+\ell_c}(h)}{\ell_c},
 \qquad A_j(h)=\frac{\binom{m-j}{h}}{\binom mh}.
\end{equation}
Generate \(A\) as the first \(h\) labels of a uniform ordering of \(C\).
Select \(c\) if its forbidden set is avoided and it is first in \(D_c\),
with that first label in \(A\). Adjacent vertices cannot both be selected.
The event that \(A\) avoids \(\Omega_c\) and meets \(D_c\) has probability
\(A_{k_c}-A_{k_c+\ell_c}\), and is invariant under permutations of \(D_c\).
Each label is equally likely to be first. This proves \eqref{eq:selection}.

We will use the identity
\begin{equation}\label{eq:selection-monotonicity}
 \frac{A_k-A_{k+\ell}}{\ell}
 =\frac1\ell\sum_{j=k}^{k+\ell-1}
       \frac{\binom{m-j-1}{h-1}}{\binom mh}.
\end{equation}
For \(h>0\) the summands are nonnegative and decrease with \(j\).
The quotient therefore decreases in both \(k,\ell\), on \(k+\ell\le m\),
and lies in \([0,h/m]\). At \(h=0\) it is zero.

\subsection{The original-pool priority coefficient}
Suppose the fixed \(d\)-star family has \(T\ge T_{\min}\ge d\) centres,
and each unused centre has at least \(d\) neighbours in its original
endpoint set. Choose one incident \(d\)-set at each unused centre.
Form a graph on the unused centres by joining two centres if their
chosen \(d\)-sets intersect. Even centres with the same chosen set are
distinct vertices of this graph. Each endpoint is adjacent to its owner
and to the centre that chose it, so it meets at most \(P-2\) other
unused centres. Each closed neighbourhood therefore has size at most
\[
 \ell=\min\{1+d(P-2),\,m-T_{\min}\}.
\]
The owners of the chosen endpoints form a forbidden set of size at most
\(d\), disjoint from the conflict set. Equations
\eqref{eq:selection}--\eqref{eq:selection-monotonicity} give
\[
 \E R_{E_0}(A)\ge hT/m+(m-T)\delta_h,\qquad
 \delta_h=\frac{A_d(h)-A_{d+\ell}(h)}{\ell}.
\]
When \(\ell=0\), only \(T=m\) occurs and the rank is exactly \(h\);
set \(\delta_h=0\). For \(\ell>0\), \(d+\ell\le m\) and
\(0\le\delta_h\le h/m\).
The selected stars can be added to all the available original stars.
The random ordering bounds the maximum packing size. In the probability
argument we may still choose a maximum packing by the fixed graph order.

For the four-star range 12--40, this is combined with
\eqref{eq:affine-mean}, using \(P=p+3\), \(d=4\), and the lower
endpoint \(T_{\min}=T_0(p)\) in the main text.
For the two-star original pool at \(p=3\), take \(m=747\), \(P=6\),
\(d=2\), \(T_{\min}=187\), and \(\ell=9\).
Its coefficient is
\[
 \begin{aligned}
 \beta_h^{\rm old}=\max\Bigl\{&
 \frac{\bigl[(1-h/m)(2/5)
    \E\min\{\Hypergeom(m-1,h,5),3\}-h/m\bigr]_+}{6},\\
 &\frac{A_2(h)-A_{11}(h)}9\Bigr\}.
 \end{aligned}
\]
At \(h=0,m\) the coefficient is zero. The positive-part term comes from
\eqref{eq:affine-mean} with \(a_0=2,s=5,b=3\). Since the enlarged set
contains \(E_0\), the coefficient also gives a lower bound for the
expected maximum packing size there. The Local Lemma allowance still
uses the full enlarged set.

\subsection{The deficit-two mean}
For \(p=2\), \(m=748\), choose a maximum good two-star packing of size
\(T\ge187\), with centre set \(I\) and endpoint set \(E_0\). Let \(J\)
consist of the centres with at least two good neighbours in \(E_0\), so
\(I\subseteq J\). Each centre outside \(J\) has at least two good
neighbours in total. At most one can lie outside \(E_0\), since two
would extend the packing. It therefore has exactly one good neighbour
in \(E_0\) and one outside. Adjoin the latter endpoint for every such
centre. We can now choose a pair of endpoints at every centre, and the
enlarged set satisfies \(|E'|\le2T+m-|J|\le m+T\).

Choose an old--old pair at every centre of \(J\setminus I\), and the
old--new pair at every centre outside \(J\).
A pair of original endpoints has at most two owners and a closed
neighbourhood of size at most \(1+3+3=7\) in the conflict graph. A pair
with one added endpoint has one owner and a closed neighbourhood of
size at most \(1+3+4=8\), since the added endpoint has no original
owner. In both cases the owners lie outside the closed neighbourhood.
Formula~\eqref{eq:selection} yields
\[
 \E R_{E'}(A)\ge hT/m+
 (|J|-T)\frac{A_2-A_9}{7}
 +(m-|J|)\frac{A_1-A_9}{8}.
\]
The last quotient averages the first eight terms of a decreasing
sequence in \eqref{eq:selection-monotonicity}, while the preceding
quotient omits its first term. Hence
\[
 \frac{A_1-A_9}{8}\ge\frac{A_2-A_9}{7}=:\delta_h,\qquad
 \E R_{E'}(A)\ge hT/m+(m-T)\delta_h.
\]
The same coefficient works for every \(T\) in the range and every
sample size. At \(h=0,m\) it is zero. At \(h=1\) it is \(1/m\),
giving the exact rank one.

The paragraph on the deficit-three mean in Section~\ref{sec:p3-mean}
derives the bound for free pairs and combines the three coefficients.

\section{Finite arithmetic and reproducibility}
\label{app:verification}

This appendix describes the finite checks in
Proposition~\ref{prop:finite}. The source programs, exact inputs and
reference outputs are in the \texttt{verification} directory of the
supplement archived at
\href{https://doi.org/10.5281/zenodo.23878887}{10.5281/zenodo.23878887}.
The release is \texttt{BK750-paper-v009-supplement-20261102}.
Its entry points are \texttt{README.md}, \texttt{REPRODUCE.md},
\texttt{PROFILE\_MAP.json} and \texttt{data-index.json}. The data
index gives the stored path for each source key and the SHA-256 hashes
of both its stored and decoded forms. The profile map identifies the
inputs and outputs summarised below. The programs evaluate the events,
recurrences and support bounds specified in the paper.

\subsection{Exact inputs and directed enclosures}
Integer quantities, binomial and Stirling numbers, sample probabilities
and rank-mean coefficients are evaluated exactly.
A binary floating-point endpoint is decoded as its exact dyadic rational
when making a final comparison. Survival floors are rounded down.
Coefficients of the positive recurrences, grouping sums and omitted
positive tails are rounded up. Allowances of the form
\(x_0\exp(-\Lambda)\) are rounded down.

The exponential enclosures use a Taylor series after reducing the
argument to a fixed small interval, with a geometric bound on the
remaining positive tail, followed by directed squaring and reciprocation.
A separate alternating-series enclosure checks the allowance direction.
Factorial roots are bounded by rational \(\ell_j,u_j\) satisfying
\[
 \ell_j^j\le j!\le u_j^j,\qquad
 \ell_j>u_{j-1},\qquad 2\ell_j>u_{j-1}+u_{j+1}.
\]
The finite table through 1502 therefore proves the monotonicity and
concavity used in the interpolated permanent bound.

The numerical kernels are evaluated with strict directed rounding and
without reassociation or contraction of floating-point operations.
The \(2^{-40}\) mark probabilities and their complements are exactly
representable. We include the coefficient remainders from the preceding
appendices. Normalisation checks serve as diagnostics. Each final
inequality must hold strictly between the rational bounds.

\subsection{Coverage of the finite checks}
For each root degree \(n=749,750\), the neighbourhood computation uses
one total-weight bound for \(c<250\) and a separate bound for each
\(250\le c\le634\). This gives 772 inequalities. Each uses all
matching sizes \(0\le r\le\lfloor n/2\rfloor\), for a total of
289886 matching records. The computation reconstructs the coefficients
for restoring the root, the class-size remainders and the bounds on
factorial roots. The largest ratio of an upper bound to its allowance
is less than 0.916.

Table~\ref{tab:cases} covers every clique deficit \(1\le p\le115\).
The structural arguments give the initial packing bounds. We check
their numerical values and verify the integer thresholds defined by
polynomial inequalities at consecutive integers. On each interval where the
mean bound is affine, we use the same coefficients, failure curves and
exponential envelopes from Section~\ref{sec:rank-envelopes} throughout.
We check the envelope inequality at every prefix. For a fixed family
of stars we average over the full hypergeometric distribution. At
occupancies outside the stored range we use the failure bound one.

There are 187 final comparisons. For each of the 75 high deficits we
take the larger of the two bounds obtained by averaging over occupancies,
using the same endpoint bound of 2000 in both cases. For \(4\le p\le10\)
we reconstruct the coefficient \(\beta_h\) in \eqref{eq:affine-mean}
for every \(0\le h\le m\) and check both ends of each interval, giving
14 comparisons.
For the sensitivity comparisons in Section~\ref{sec:closure}, use
the same stored envelopes, occupancy weights and variance terms.
Evaluate \eqref{eq:outer-envelope} at the specified \(T\) with endpoint
bound \(m+T\); for the comparison without recovery replace
\(\mu_h(T)\) by \(hT/m\) alone. The six generic left endpoints for
\(p=5,\ldots,10\) are \(187,207,223,237,247,257\).
These additional comparisons show the effect of changing the initial
packing bound or omitting the recovery of stars. The other 98 final
comparisons use the specified curves and the occupancy weights for
\(m=750-p\). The two-star and four-star profiles use their respective bounds for
colour classes, endpoints and representatives. Throughout the
calculation we use the actual clique order and deficit.

\paragraph{The deficit-two comparisons.}
The stored comparisons for deficit two use the scalar reference function
\(e_{\rm ref}(T)\) below. The endpoint-set bound in the dependency estimate
is \(e(T)=m+T\), with \(m=748\). We obtain the required allowance by
bounding the difference between these functions. Put
\[
 e_{\rm ref}(T)=\min\{m+T,14T-3m+5\},\qquad
 c=\eta_0q^3+\eta_Wq.
\]
For \(187\le T\le374\),
\[
 0\le e(T)-e_{\rm ref}(T)\le556,\qquad 556c<1.
\]
Let \(U_T\) be the stored upper bound on the failure probability and
\(L_T^{\rm ref}\) the stored lower bound on the reference allowance.
The data give \(3U_T<L_T^{\rm ref}\) at \(T=187,229,230,374\).
The exponent in the allowance is affine in \(e\). Since
\(\exp((e(T)-e_{\rm ref}(T))c)<\exp(1)<3\), these comparisons imply
\[
 U_T<x_0\exp[-\eta_0J_3(p,m+T)-\eta_WJ_1(p,m+T)].
\]
The largest of the four outward-rounded values of \(3U_T/L_T^{\rm ref}\)
is \(0.453440593437\). The comparison for 375 fixed stars uses their
750 endpoints and gives a ratio below \(0.064\). Applying
\eqref{eq:endpoints} separately on 187--229 and 230--374 covers every
smaller packing size in the stated range.

\begin{center}\small
\begin{tabular}{@{}lp{99mm}@{}}
\toprule
Map entry & Named input profiles and final comparison outputs\\
\midrule
A comparisons & The 772-row \texttt{A\_activity\_comparisons.json} ledger.\\
p1 & \texttt{p2\_input}, \texttt{p2\_prefix}; both complete p1 curves,
smoothed occupancy rows and \texttt{activity\_p1p2}.\\
p2 & The same inner profile; fixed p2 envelope and endpoints
187, 229, 230, 374, 375; \texttt{activity\_p1p2}.\\
p3 & \texttt{LOW}, \texttt{env\_3}, \texttt{coeff\_3}, \texttt{result\_3}.\\
p4--10 & \texttt{LOW}/\texttt{HIGH};
\texttt{certificates/unified-low/comparisons.json}.\\
p11 & \texttt{P20}, \texttt{env\_11}; the p11 rows of
\texttt{expanded-dummy-owner} selected-comparison file.\\
p12--40 & \texttt{P20}/\texttt{P40}, the corresponding
\texttt{env\_p}, \texttt{result\_p} and ordinary-250 records.\\
p41--115 & The unified high \texttt{exact-input.json} and complete curves;
\texttt{certificates/unified-high/comparisons.json}.\\
\bottomrule
\end{tabular}
\end{center}
The abbreviations in the table are source keys in
\texttt{code/portable-F/INPUTS.json}. The profile map gives their paths
within the archive. The two unified directories each have a
\texttt{FILES.json} manifest. Historical names in the source keys refer
to files within the archive.

\subsection{Programs and their roles}
The archive includes the following checks.
\begin{enumerate}
 \item \path{code/portable-A/verify_a_global.py}\par
 This program reconstructs all
 neighbourhood coefficients and matching bounds, the 772 final inequalities,
 the factorial-root brackets, structural guards, and the \(W\)-bound.
 \item \path{code/portable-F/verify_f.py} reconstructs the full
 sequence of calculations for the conditional failure curves and checks
 262 reference comparisons. The present proof uses 98 of these
 comparisons. The remaining comparisons required for the proof are
 obtained below. The full set of reference checks also validates the
 inputs for the curves.
 \item \path{code/refinement/unify_low.py} checks the seven common
 two-star supplies, all mean coefficients and curve envelopes, and the 14
 full-support endpoint comparisons.
 \item \path{code/refinement/unify_high.py} reconstructs the single
 filtered profile, its complete laws and grouping coefficients, and the 75
 common-support comparisons. Its recurrence uses the positive kernel
 of the reference profile, with \(p_*=115\). Thus the upper bound on the
 number of outside neighbours of a clique representative is
 \(p_*+1=116\), as in Section~\ref{sec:closure}.
 \item \path{code/refinement/verify_p2_support_bridge.py} checks
 the four scalar inequalities for deficit two, the comparison with
 375 fixed stars and the bound on the change in the dependency exponent.
 \item \path{code/portable-F-producer/replay_f_producers.py}
 reconstructs the exact inputs and shared weights, checks small exact
 fixtures, and replays seven reference profiles. It reproduces 31 scientific
 output files and five fixtures byte for byte.
\end{enumerate}

The directories \path{certificates/unified-low} and
\path{certificates/unified-high} contain the exact comparisons and
\texttt{FILES.json} manifests. Reproduction checks their five and sixteen
declared scientific outputs, respectively, byte for byte. The archive
also records interval coverage and all profile parameters.

In the published ZIP, \path{DEPOSIT.md} describes the release layout.
The manuscript is distributed separately. Run the commands below
from the extracted \path{verification/} directory, one job at a time,
using fresh output directories. The verifiers create their shared lock in
\path{/private/tmp}, which must exist and be writable. Python~3.9 or later
is required; the producer and unified high-deficit runs also require a
C++17 compiler. Compiler options and SDK selection are documented in
\path{REPRODUCE.md}. The principal commands are:
\begin{verbatim}
python3 -B code/portable-F/verify_f.py \
  --packet-root . --out /tmp/bk750-F
python3 -B code/portable-F-producer/replay_f_producers.py \
  --packet-root . --out /tmp/bk750-producers
python3 -B code/refinement/unify_low.py \
  --packet-root . --out /tmp/bk750-low
python3 -B code/refinement/unify_high.py \
  --packet-root . --out /tmp/bk750-high
python3 -B code/refinement/verify_p2_support_bridge.py \
  --packet-root .
\end{verbatim}
The neighbourhood and global checks also require the source and input
hashes listed in the archive's \texttt{REPRODUCE.md}. The execution
contracts in the archive give the time and memory limits for each
program.

Replaying the producer reproduces the implementation's reference
outputs. The exact checks of the inputs, weights, test cases and
coefficient transformations evaluate the inequalities derived from
the recurrences and probability bounds above. Together these checks
give the computer-assisted bounds used in the proof.

\end{document}